\documentclass[a4,12pt]{article}%
\usepackage{graphicx}
\usepackage{amsmath}%
\usepackage{amsfonts}%
\usepackage{amssymb}

\begin{document}

\begin{center}

 {\large{\textbf{Representation of symmetric shift registers}}  \vspace{ 1 mm}}

 Jan S\o reng, April\,06, 2026.

 Faculty of  Science and Technology, Norwegian
University of Life Sciences.   \vspace{ 0mm}

\end{center}

\noindent The objective of this work is to establish a mathematical framework for the study of symmetric shift registers over the field
GF(2). The present paper gives a new approach where the symmetric shift registers are represented by associated systems of nonlinear difference equations. Arithmetical progressions will play a central part. This approach clarifies the underlying structures and makes it easier to determine the minimal periods of the sequences generated by the symmetric shift registers. Key words: Shift registers, nonlinear difference equations, periods, arithmetical progressions, GF(2).\vspace{1mm}

\noindent {\large{\textbf{1. Introduction.}}}\vspace{1mm}

\noindent Symmetric shift registers are nonlinear difference
equations with certain symmetry properties. Kjell Kjeldsen studied
some classes of such registers in his paper [3] published in 1976.
In the papers [6], [7] and [8], completed around 1980, the cycle
structure of every symmetric shift register was determined. However,
these papers contain a lot of difficult combinatorial arguments and
most of the ideas and concepts are hidden in complicated proofs.
This serves as a motivation of the present work. The objective of
the present paper is to give a completely new approach by constructing a theory
which clarifies the ideas and simplifies the proofs. Thus it can be
viewed as a complement to the papers [6], [7] and [8].\vspace{1mm}

\noindent   There exist some results of other types
of nonlinear registers.   See for instance [4] by Johannes Mykkeltveit and [5] by
Johannes Mykkeltveit, Man-Keung Siu and Po Tong. In [2] Tor Helleseth  gives a survey of nonlinear shift registers and open problems.\vspace{1mm}

\noindent The symmetric shift registers of length $n$ are symmetric
difference equations which generate periodic sequences of infinite
length by starting with $n$ bits. If $S(x_2, \cdots , x_n)$ is a
symmetric polynomial, the corresponding shift register is
$a_{\,n+1}=a_{\,1}+ S(a_{\,2}, \cdots , a_{\,n})$\,. In  [6] we
proved that every symmetric polynomial $S$ in the variables
$x_{2}, \cdots , x_{n}$ has the form\vspace{1mm}

\hspace{8mm}$S(x_2, \cdots , x_n)=\sum_{j\,=\,0}^{\,n\,-\,1} c_j\, E_j (x_2, \cdots
, x_n) \mbox{\ \ where \ }  c_j \in \left\{0,1 \right\}$\vspace{1mm}

\noindent and $E_j$ is the elementary symmetric polynomial of order $j$ which
is defined in Section 5.  Examples indicate that the minimal periods
are not so large.\vspace{0.8mm}

\noindent  {\footnotesize{\emph{$^1$The author retired from the University in 2012.}}}

\noindent For instance, the minimal periods of the examples
in this paper have values less than $n^3$. In [6] the cycle
structure problem was reduced to the \mbox{case
$S=E_{k}+E_{k+1}+\cdots \,+\,E_{k+p}$\, where $0 \le k \le k+p<
n$.}\vspace{1mm}

\noindent The main problem in [8] was to determine the minimal
periods for each $A \in \{0,1 \}^n$. However, some parts of the deductions in [8] were
only sketched. This paper contains a rigorous new approach to  the problem.\vspace{1mm}

\noindent If $A \in \{0,1\}^n$, we denote the infinite sequence
generated by the shift register by $A^{\infty}$. The sequence
 $A^{\infty}$ seems rather chaotic. But by describing the dynamics with appropriate concepts the underlying structure is revealed.  We introduce certain invariants, and by an inductive reduction process the periods can be determined. Arithmetical progressions will play a central role in  each step of this process. The first two parts of this paper contain the formulation of the results,  examples and visualization of the underlying structure on an example. In Part 3 and 4 the most important concepts used  in the proofs are introduced. Part 5 describes the main lines of the proofs. The remaining parts contain proofs of  results that will be used in this connection. In the end of this paper we have included an index.\vspace{1.5mm}

\noindent{\large{\textbf{2. Acknowledgement.}}}\vspace{1mm}

\noindent I am particularly grateful to Professor Helge Tverberg,
University of Bergen, Norway, who has encouraged me, worked
through different versions of the manuscript and given a great
number of important suggestions. I am also
very grateful to Dr. Kjell Kjeldsen, Headquarter of Defence, Norway,
who presented the problem to me in 1974 and encouraged me to
continue his work [3] on symmetric shift registers. Furthermore I would
like to  thank Professor John Wyller, Norwegian University of Life Sciences, for very helpful discussions. Moreover,
I  am very grateful to my sons Eirik, Martin and P\aa l  for comments and for verifying the results by  computer programs.
Moreover, I owe a special thank to my wife Kari for encouragement
and support.\vspace{2mm}

\hspace{45mm} {\large{\textbf{PART 1.}}  \vspace{ 1 mm}}

\noindent We will in this part introduce necessary notation and definitions.  In particular, we define the symmetric shift registers we will study, in Section 5.  We  also describe how strings can be represented by vectors. Using  vector representations makes it easier to formulate and prove the results.\vspace{0.4mm}

\noindent Suppose $V$ is a vector in $M^*$ where $M^*$ is defined in Section 7.
Then we associate a contraction vector and a distance vector to $V$.
These vectors will give necessary information about the inner structure of the given vector $V$.\vspace{12mm}

\noindent \textbf{{\large 3. Terminology.}}\vspace{1mm}

\noindent Suppose $A=(a_1, \cdots ,a_m)$ where $a_i \in \{0,1\}$ for $1 \le i \le m$\hspace{0.25mm}.
We write  $A=a_{\hspace{0.1mm}1}\, \cdots \,a_{\hspace{0.25mm}m}$ instead of
$A=(\hspace{0.1mm}a_{\hspace{0.1mm}1},\, \cdots ,a_{\hspace{0.25mm}m}\hspace{0.1mm})$\, and call
$A$ \mbox{a string}. The length of $A$ is denoted by $length(A)$\hspace{0.25mm}. That means, $length(a_{\hspace{0.1mm}1}\, \cdots \,a_{\hspace{0.25mm}m})=m$\hspace{0.25mm}.\vspace{1mm}

\noindent We will only consider substrings of $A$ with
adjacent elements of the form $P=a_{\hspace{0.1mm}i}\, \cdots a_{\hspace{0.1mm}j}$\, where
$1 \le i \le j \le m$\hspace{0.25mm}.\vspace{1mm}

\noindent $P\subset A$ means that $P$ is a substring
of $A$\hspace{0.25mm}, i.e. $A=RPS$ where $R$ and/or $S$ may be the empty string. The empty string
is denoted by
$\O$.\vspace{1mm}

\noindent Moreover, $A^{\,\prime}$ denotes the
string $A\,' =a^{\,\prime}_{\hspace{0.1mm}1} \cdots a^{\,\prime} _{\hspace{0.25mm}m}$ where
$a^{\,\prime}_{\hspace{0.1mm}j} =1$ if $a_{\hspace{0.1mm}j} = 0$, \mbox{and $a^{\,\prime} _{\hspace{0.1mm}j}
=0$} if $a_{\hspace{0.1mm}j} =1$, for $1\leq j\leq m$\hspace{0.25mm}. That means,
$a^{\,\prime}_j=1-a_{\hspace{0.1mm}j}$ for $1\leq j\leq m$.\vspace{1mm}

\noindent We define the weight of 1 and 0 as $w(1)=1$ \ and \ $w(0)=0$. The weight of $a_{\hspace{0.1mm}1}\, \cdots
a_{\hspace{0.25mm}m}$ is equal to the number of 1's in $a_{\hspace{0.1mm}1} \,\cdots
\,a_{\hspace{0.25mm}m} $\,, and will be \mbox{denoted $w( a_{\hspace{0.1mm}1} \,\cdots \,a_{\hspace{0.25mm}m}
)$\hspace{0.25mm}.} That means,\vspace{1mm}

\noindent (3.1)\hspace{14mm} $w(a_{\hspace{0.1mm}1} \cdots a_{\hspace{0.25mm}m})=w(a_{\hspace{0.1mm}1}) + \cdots + w(a_{\hspace{0.25mm}m})=a_{\hspace{0.1mm}1} + \cdots + a_{\hspace{0.25mm}m}$\hspace{0.23mm}.\vspace{1mm}

\noindent The positive weights of 1 and 0 are defined as $\overline{w}(1)=1$ and $\overline{w}(0)=-1$\hspace{0.25mm}.\vspace{0.5mm}

\noindent The positive weight of $A=a_{\hspace{0.1mm}1}\, \cdots
a_{\hspace{0.25mm}m}$ is  $\overline{w}(a_{\hspace{0.1mm}1} \cdots a_{\hspace{0.25mm}m})=\overline{w}(a_{\hspace{0.1mm}1}) + \cdots + \overline{w}(a_{\hspace{0.25mm}m})$\hspace{0.25mm}.
That means, \mbox{$\overline{w}(A)$ = the number of ones in $A$ $-$ the
number of zeros in $A\hspace{0.25mm}.$}  If $A=\O$, then $w(A)=\overline{w}(A)=0$\hspace{0.25mm}. If $a=0$ or $a=1$\hspace{0.25mm}, \mbox{then $\overline{w}(a^{\,\prime})=-\,\overline{w}(a)$\hspace{0.25mm}.}\vspace{0.5mm}

\noindent Moreover, $0_j$ denotes $j$ consecutive
0's and $1_j$ denotes $j$ consecutive 1's\hspace{0.25mm}. \mbox{If $j=0$\hspace{0.2mm},} then $0_j$ and
$1_j$ denote the empty string.\vspace{1mm}

\noindent If $A^{\,\infty}=(a_1, a_2, \cdots )$ where $a_i \in \{0,1\}$ for $i \ge 1$, then $A^{\,\infty}$ is regarded as an infinite string
and we write  $A^{\,\infty}=a_1 a_2 \cdots $\,. We define substrings as in the finite case. We call $r \ge 1$ a period of $A^{\,\infty}$
if $a_{r+i}=a_i$ for $i \ge 1$\hspace{0.25mm}.\vspace{1mm}

\noindent If variables, parameters or coordinates of vectors are not specified, they have non-negative integer values. Moreover, $\#$ means "the number of\,".\vspace{2mm}

\noindent \textbf{{\large 4. Vectors.}}\vspace{0.8mm}

\noindent We will mainly use uppercase letters to denote vectors, and lowercase letters to denote the coordinates.
However, $I$ and $J$ will denote integers.\vspace{1mm}

\noindent Suppose  $V=(v_1, \cdots, v_t)$ is a vector where $t >0$\,. If $t$ is odd  or even, we call $V$ an odd or even vector respectively. The number of coordinates in $V$ is denoted by $\#V$. That means, $\#V=t$\hspace{0.25mm}.\vspace{0.5mm}

\noindent Moreover, we define $sum(V)=v_{\hspace{0.1mm}1}+v_{\hspace{0.1mm}2} + \cdots +v_{\hspace{0.1mm}t}$\hspace{0.3mm}.   If $\alpha$ is an integer, we let\vspace{0.5mm}

\hspace{3.5mm} $V+\alpha=(v_{\hspace{0.1mm}1}+\alpha, \cdots, v_{\hspace{0.1mm}t}+\alpha)$ and $V-\alpha=(v_{\hspace{0.1mm}1}-\alpha, \cdots, v_{\hspace{0.1mm}t}-\alpha)$\hspace{0.25mm}.\vspace{10.5mm}

\noindent If $B=(v_r, \cdots, v_s)$ where $1 \le r \le s \le t$\hspace{0.25mm}, then
we write $B \subset V$ and call $B$

\noindent a sub-vector of $V$.\vspace{0.5mm}

\noindent The empty vector is denoted $\O$ or $(\O)$. If $V=\O$, we let $sum(V)=0$\hspace{0.25mm}.\vspace{1mm}

\noindent If $V_1$\hspace{0.25mm}, $V_2$\hspace{0.25mm}, $\cdots$  are vectors, then $(V_1,V_2, \cdots )$ is the vector consisting successively of the coordinates of $V_1$\,, $V_2$\,, $\cdots$\,. In particular, $(V_1)=V_1$\hspace{0.25mm}.\vspace{1mm}

\noindent If  $V^{\,\infty}=(v_1,v_2, \cdots )$\hspace{0.25mm},  then we define $sum(V^{\,\infty},j)=v_1+ \cdots +v_j$   for $j \ge 1$\hspace{0.25mm}.\vspace{1mm}

\noindent If $V=(v_1, \cdots, v_t)$\hspace{0.2mm}, we let $sum(V,j)=v_1+ \cdots +v_j$   for $1 \le j \le t$\hspace{0.15mm}.\vspace{1mm}

\noindent If $V$ is a vector, then the extension $V^*$ of $V$ is constructed by adding one to the last coordinate of $V$.
 For example,  the extension of $V=(1,1,3,2,2,4)$ is given by $V^*=(1,1,3,2,2,4+1)=(1,1,3,2,2,5)$\hspace{0.25mm}.\vspace{0.4mm}

 \noindent We call $r>0$ a vector period of $V^{\,\infty}=(v_1,v_2, \cdots)$ if $v_{r+i}=v_i$ for $i \ge 1$.\vspace{3mm}

 \noindent {\large{\textbf{5. Symmetric shift registers.}}}\vspace{1mm}

\noindent Let $k$, $p$ and $n$ be integers such that $0 \le k \le k+p <
n$\hspace{0.25mm}. In this paper we study the symmetric shift register
$\theta:\{0,1\}^n \rightarrow \{0,1\}^n $ with respect to $k$, $p$ and $n$ defined by
$\theta (a_{\hspace{0.1mm}1} \cdots \,a_{\hspace{0.2mm}n} )\;=\;a_{\hspace{0.2mm}2}
\,\cdots a_{\hspace{0.2mm}n+1} $ where \vspace{0.2mm}

\noindent \hspace{8.2mm} $a_{\hspace{0.2mm}n+1} =a^{\,\prime}_{1} $ if $k\leq a_{\hspace{0.2mm}2}
+\;\cdots \;+a_{\hspace{0.2mm}n} \leq k+p$, and $a_{\hspace{0.2mm}n+1} =a_{\hspace{0.1mm}1} $
otherwise,\vspace{0.2mm}

\noindent where  $a^{\,\prime} _{\hspace{0.1mm}1} =0$
 if $a_{\hspace{0.1mm}1} =1$, and $a^{\,\prime} _{\hspace{0.1mm}1} =1$
 if $a_{\hspace{0.1mm}1} =0$.\vspace{0.5mm}

\noindent The elementary symmetric polynomial $E_j$ of order $j$ has
the property\vspace{0.5mm}

\hspace{26mm}$E_j(x_{2}\,,\cdots,\,x_{n})=1\;\Leftrightarrow
\;x_{2}+\cdots+x_{n}=j\hspace{0.15mm}.$ \vspace{0.5mm}

\noindent We let $S=E_{k}+E_{k+1}+\cdots \,+\,E_{k+p}$\hspace{0.34mm}.
Since\vspace{0.5mm}

\hspace{16.2mm}$S(x_{2} ,\;\cdots \;,x_{n} )=1\;\Leftrightarrow \;k\leq
x_{2} +\;\cdots \;+x_{n} \leq \;k+p$\hspace{0.1mm},\vspace{0.5mm}

\noindent the symmetric shift register $\theta $ corresponds to the
difference equation given by
$a_{\hspace{0.2mm}n+1}=a_{\hspace{0.1mm}1}+S(a_{\hspace{0.2mm}2},\cdots,a_{\hspace{0.2mm}n})$ (\,mod 2\,).

\noindent If $A^{\,\infty}=a_1a_2, \cdots$ is the infinite sequence generated from  $A=a_1 \cdots  a_n$  by the symmetric shift register, then
for each $i \ge 0$ we get that

\noindent \hspace{22.7mm} $a_{\hspace{0.2mm}n+i+1} =a^{\,\prime}_{i+1} $ if $k\leq a_{\hspace{0.2mm}i+2}
+\;\cdots \;+a_{\hspace{0.2mm}i+n} \leq k+p$,

\noindent and $a_{\hspace{0.2mm}n+i+1} =a_{\hspace{0.1mm}i+1} $
otherwise. The main problem is  to determine  the \mbox{least $r \ge 1$} such that $a_{r+i}=a_i$ for $i \ge 1$.\vspace{1.5mm}

\noindent {\large{\textbf{6. Infinite vector representations.}}}\vspace{1mm}

\noindent Let $A^{\,\infty}=a_1a_2a_3 \cdots$ be generated from  $A=a_1 \cdots  a_n$ by the symmetric shift register with parameters $k$, $p$ and $n$\hspace{0.2mm} where  $A$ starts with 1. It is easily seen there exists $i \ge 1$ such that  $a_i=0$. Since $A^{\,\infty}$ is periodic, then
 $A^{\,\infty}$ contains an infinite number of bits equal to 0 and an infinite number of bits
 equal \mbox{to 1.} Hence, $A^{\,\infty}$ can be decomposed as\vspace{0.3mm}

 \noindent (6.1)\hspace{20mm} $A^{\,\infty}=1_{q_1}0_{q_2}1_{q_3}0_{q_4}  \cdots$  where $q_i>0$ for $i \ge 1$\hspace{0.1mm}.\vspace{0.5mm}

  \noindent Then $V(A^{\,\infty})=(q_1,q_2, \cdots )$ is called
 the vector representation of $A^{\,\infty}$.
 Next, we suppose $w(A)=k+p+1$.  Then we get according to Observation 45.3 and Proposition 45.4  that\vspace{0.5mm}

 \noindent (6.2)
 \hspace{27mm}$V(A^{\,\infty})$ has an even vector period,\vspace{0.5mm}

 \noindent (6.3) \hspace{6mm}$sum(V(A^{\,\infty}),j)$ is the
minimal period of $A^{\,\infty}$ if  $j$ is the  least

  \hspace{10mm}even vector period of $V(A^{\,\infty})$.\vspace{2mm}

\noindent {\large{\textbf{7. Sets of vectors.}}}\vspace{1mm}

\noindent Let $M=\{(v_1, \cdots, v_{J+1}): J \ge 1$ is odd, $v_i \ge 1$ for $1 \le i \le J$, and $v_{J+1} \ge 0$\}.\vspace{0.5mm}

 \noindent Moreover, let $M^*=\{(v_1, \cdots, v_{J+1}) \in M: v_1>1\}$ and\vspace{0.4mm}

 \noindent  \hspace{4.5mm}$ M_p=\{(v_1, \cdots ,v_{J+1}) \in M : v_1+v_3+ v_5+ \cdots +v_J \ge p+1\}$ for $p \ge 0$.\vspace{0.5mm}

\noindent In the formulation of the results the  string $A=a_1 \cdots a_n$ always \mbox{starts with 1.}
Suppose $A=a_1 \cdots a_n$ where $a_1=1$.
 Then $A $ has an even  decomposition of the form $1_{v_1} 0_{v_2}1_{v_3}0_{v_4} \cdots 1_{v_{J}}0_{v_{J+1}}$
 where $v_i >0$ \mbox{for $1 \le  i \le J $}, $v_{J+1}\ge 0$ \mbox{and $J \ge 1$ is odd.} We call the vector $V(A)=(v_1, \cdots,v_{J},v_{J+1})$ the even vector representation \mbox{of $A$\hspace{0.25mm}.} We note \mbox{that $V(A) \in M$.}\vspace{0.5mm}

 \noindent If $V=(v_1, \cdots, v_{J+1}) \in M$, we let $A(V)=1_{v_1}0_{v_2}1_{v_3} \cdots 1_{v_J}0_{v_{J+1}}$. We note that $A(V)$ starts with 1. The next observations are trivial.\vspace{2mm}

\noindent {{\textbf{Observation 7.1.}}} a) If $A=A(V)$ where  $V \in M$, then $V=V(A)$.\vspace{0.3mm}

\noindent b) If $V=V(A)$ where $A=a_1 \cdots a_n$ starts with $a_1=1$, then $A=A(V)$.\vspace{2mm}

\noindent {{\textbf{Observation 7.2.}}}  Suppose $A(V)=1_{v_1}0_{v_2}1_{v_3} \cdots 1_{v_J}0_{v_{J+1}}$
where $V \in M$. Then $V=(v_1, \cdots ,v_{J+1})$.\vspace{2mm}

\noindent \noindent {{\textbf{Example 7.3.}}} If  $A=110011100$\hspace{0.35mm}, then we  decompose $A=1_20_21_30_2$\hspace{0.35mm} and get that
$V(A)=(2,2,3,2) \in M$ is the even  vector representation of $A$\hspace{0.25mm}.\vspace{1.5mm}

\noindent \noindent {{\textbf{Example 7.4.}}}   If $A=111100111$\hspace{0.25mm}, then we  decompose $A=1_40_21_30_0$\hspace{0.35mm} and get that
$V(A)=(4,2,3,0)\in M$ is the even  vector representation of $A$\hspace{0.25mm}.\vspace{1.5mm}

\noindent \noindent {{\textbf{Example 7.5.}}} If  $V=(3,4,2,0)$, then $A(V)=1_30_41_20_{\hspace{0.2mm}0}=111000011$.\vspace{1.5mm}

  \noindent \noindent {{\textbf{Example 7.6.}}} If  $V=(2,3,1,3)$, then
$A(V)=1_20_31_10_3=110001000$.\vspace{2mm}

\noindent {\large{\textbf{8. The contractions of vectors.}}}\vspace{0.5mm}

\noindent The function $\pi$ introduced in Section 13 contracts vectors. The definition
of $\pi$ is based on the distance measure and the distance function defined in Section 9 and 10.
In Section 11 and 12  we define proper odd components and component decompositions.  By using these decompositions we will  give a more intuitive and less algorithmic characterization of the function $\pi$ in Section 13.
Next, we  describe  briefly how the function $\pi$ will be used.\vspace{0.5mm}

\noindent  Suppose  $A^{\,\infty}$ is  generated from  $A$  by the symmetric shift register with parameters $k$, $p$ and $n$. If  $A$ satisfies certain assumptions,  the minimal period
\mbox{of $A^{\,\infty}$}
can be determined by an inductive process. Here follows a brief overview of this process.\vspace{0.5mm}

\noindent Let  $Q_p$ be the even vector representation of  $A$.  If $p>0$,  let\vspace{0.4mm}

\noindent  (8.1) \hspace{17.8mm}$Q_{p-1}=\pi(Q_p), Q_{p-2}=\pi(Q_{p-1}), \cdots, Q_0=\pi(Q_1)$.\vspace{0.6mm}

\noindent Then we use the structure of $Q_0$ to find  parameters $j_0$ and $\zeta_0$. Based on  the structure of $Q_1$ and the parameters $j_0$ and $\zeta_0$ we determine $j_1$ and $\zeta_1$. In this way we continue by induction until we find parameters $j_p$ and $\zeta_p$. Then $\zeta_p$ will be the minimal period of $A^{\,\infty}$.\vspace{0.5mm}

\noindent In Section 14 we define the set $M^+_p$, and  according to Observation 14.5 we get that  $Q_{p-1}, \cdots, Q_0$  in (8.1) are  well-defined if $Q_p \in M^+_p$.\vspace{2mm}

\noindent \textbf{{\large 9. The distance measure.}}\vspace{1mm}

\noindent If $v$ is an integer,  let $v^-=v-1$. If $V$ is a vector, let $\delta(V)=sum(V)-\#V$. In particular,  if $V=\O$, then  $\delta(V)=0$\hspace{0.25mm}. We note that\vspace{0.2mm}

\noindent    \hspace{5.3mm} $\delta(V)=sum(V)-s=v_1^-+ \cdots+ v_s^-$ if $V=(v_1, \cdots,v_{s})$ where $s \ge 1$\hspace{0.25mm}.\vspace{0.5mm}

\noindent If $V=(G_1, \cdots, G_r)$ where $r \ge 1$, then $\delta(V)=\delta(G_1)+ \cdots+ \delta(G_r)$\hspace{0.25mm}.\vspace{1.2mm}

\noindent  Suppose  $V=(v_1, \cdots,v_{s})$ where $v_j \ge 1$ for $1 \le j \le s$.
 Then $\delta(V) \ge 0$\,.
 If  in addition one of the coordinates of $V$ are larger than 1, then $\delta(V)>0$\hspace{0.35mm}.\vspace{1.6mm}

 \noindent {\textbf{Example 9.1.}} If $G_1=(5)$\hspace{0.25mm}, then $\delta(G_1)=5-1=4$\hspace{0.25mm}.\vspace{0.3mm}

 \noindent If $G_2=(3,1,2,1,5,1,4)$\hspace{0.25mm}, then $\delta(G_2)=sum(G_2)-7=17-7=10$\hspace{0.25mm}.\vspace{0.3mm}

 \noindent If $G_3=(4,1,3,1,1)$\hspace{0.25mm}, then $\delta(G_3)=sum(G_3)-5=10-5=5$\hspace{0.25mm}.\vspace{1.6mm}

\noindent {\large{\textbf{10. The distance function.}}}\vspace{0.5mm}

\noindent Suppose  $V=(v_1, \cdots,v_{J},v_{J+1}) \in M$. Then $v_r \ge 1$ for $1 \le r \le J$, \mbox{and $v_{J+1} \ge 0$.}
The distance function $\tau$ of $V$ is given by\vspace{0.4mm}

\noindent (10.1) \hspace{1mm}$\tau(0)=0$ and $\tau(r)=\delta(v_{1}, \cdots, v_r)=v_{1}^-+ \cdots+ v_{r}^-$ for $1 \le r \le J+1$.\vspace{0.4mm}

\noindent We note that  $\tau(r+1)=\tau(r)+v^-_{r+1}$ for $0 \le r \le  J$, and $\delta(V)=\tau(J+1)$,\vspace{0.4mm}

\noindent (10.2)\hspace{16mm}$\tau(s)=\tau(r)+\delta(v_{r+1}, \cdots, v_s)$
 if $0 \le r < s \le J+1$,

\noindent (10.3) \hspace{20mm} $\tau(s)>\tau(r)$  if $0 \le r < s \le J$ and $v_{r+1}>1$,\vspace{0.3mm}

\noindent (10.4) \hspace{6.3mm} $\delta(V)+1=v_{1}^-+ \cdots+v_{J}^-+ v_{J+1}^-+1=\tau(J)+v_{J+1} \ge \tau(J)$,\vspace{0.3mm}

\noindent (10.5)  \hspace{12.7mm} $0 < \tau(1) \le \tau(2) \le \cdots \le \tau(J) \le
\delta(V)+1$ if $v_1>1$.\vspace{2mm}

 \noindent \textbf{{\large 11. Odd components.}}\vspace{1mm}

 \noindent    $G=(g_{1}, \cdots , g_{2t+1})$ is called an odd component if  $t \ge 0$ and $g_{2i}=1$ \mbox{for $1 \le i \le t$\hspace{0.25mm}.}
 For example,   $(2)$\hspace{0.35mm}, $(3,1,5)$ $(3,1,2,1,5,1,4)$   are odd components.\vspace{0.8mm}

\noindent In the remaining part of this section we suppose   $V=(v_1, \cdots, v_J,v_{J+1}) \in M^*$. Then $J$ is odd, $v_1 >1$, $v_i \ge 1$ for $2 \le i \le J$, \mbox{and $v_{J+1} \ge 0$.}\vspace{0.5mm}

 \noindent  Suppose $G  \subset V$ is an odd component that ends $V$,   is succeeded by only one coordinate or  is succeeded by a coordinate $>1$. Then $G$ is called a proper odd component  in $V$.\vspace{0.5mm}

 \noindent  Suppose $0 \le r \le J$. Let $t \ge 0$ be maximal such that $v_{r+2i}=1$ for $1 \le i \le t$, and $r+2t \le J$. Then we define\vspace{0.4mm}

 \noindent (11.1)\hspace{8mm} $t_{max}(r)=t$ and $next(r)=r+2t+1=r+2t_{max}(r)+1$.\vspace{0.5mm}

\noindent When it is clear from the context  proper means proper in $V$. The functions $t_{max}$ and $next$ depend on $V$, but this  will always be clear from the context.\vspace{1.4mm}

 \noindent {{\textbf{Observation 11.1.}}} Suppose $G=(v_{r+1}, \cdots, v_{r+2t+1}) \subset V$. Then $G$ is a proper odd component if and only if $t=t_{max}(r)$.\vspace{0.5mm}

 \noindent {{\textbf{Proof.}}} If $r+2t \ge J-1$, then $G$ ends $V$ or is succeeded by one coordinate. Hence,
 $G$ is a proper odd component $\Leftrightarrow$ $G$ is an  odd component $\Leftrightarrow$ $v_{r+2i}=1$ for $1 \le i \le t$
 $\Leftrightarrow$ $t=t_{max}(r)$. Otherwise, $G$ is a proper odd component $\Leftrightarrow$ $v_{r+2i}=1$ for $1 \le i \le t$, \mbox{and $v_{r+2t+2}>1$}
 $\Leftrightarrow$ $t=t_{max}(r)$.\vspace{2mm}

 \noindent {{\textbf{Observation 11.2.}}} Suppose $G=(v_{r+1}, \cdots, v_{s}) \subset V$. Then $G$ is a proper odd component if and only if $s=next(r)$.\vspace{0.5mm}

 \noindent {{\textbf{Proof.}}} Suppose $G$ is a proper odd component. Then there exists $t \ge 0$ such that $G=(v_{r+1}, \cdots, v_{r+2t+1})$ and
  $s=r+2t+1$. By Observation 11.1 we get that
  $t=t_{max}(r)$ and $s=r+2t+1=next(r)$.
 Next, suppose $s=next(r)$. Then $s=r+2t+1$ and $G=(v_{r+1}, \cdots, v_{r+2t+1})$ where $t=t_{max}(r)$. Hence,
  Observation 11.1 implies that $G$ is a proper odd component.\vspace{2mm}

 \noindent {{\textbf{Observation 11.3.}}} Suppose $0 \le r \le J$. Then $G=(v_{r+1}, \cdots, v_{s})$ where $s=next(r)$, is the unique proper odd component starting with $v_{r+1}$.\vspace{0.5mm}

 \noindent {{\textbf{Proof.}}} Follows from Observation 11.2.\vspace{2mm}

 \noindent {{\textbf{Observation 11.4.}}} \noindent a) $(v_J)$ and $(v_{J+1})$ are proper odd components.

 \noindent b) Suppose  $v_{i+1}>1$ where $1 \le i < J$. Then $G=(v_i)$ is a proper

\noindent\hspace{5.3mm}odd component.\vspace{0.5mm}

\noindent {{\textbf{Proof.}}}  a) is true since $(v_J)$ is an odd component succeeded by only one coordinate, and
$(v_{J+1})$ is an odd component ending $V$.\vspace{0.5mm}

\noindent b)  is true since $G=(v_i)$ is an odd component succeeded by a coordinate $>1$.\vspace{2mm}

\noindent \textbf{{\large 12. The component decomposition.}}\vspace{1mm}

 \noindent Suppose  $V=(v_1, \cdots,v_{J},v_{J+1}) \in M^*$.
 If $V=(G_1, \cdots,G_{I+1})$ where   $G_j$ is a proper odd component  for $1 \le j \le I+1$\hspace{0.25mm},
 we call $(G_1, \cdots,G_{I+1})$ the component decomposition of $V$ and $G_1, \cdots, G_{I+1}$  the components of the decomposition. By  Observation 11.3 it is easily seen that
  $V$ has a unique such decomposition.
 In the  examples  the components of the decompositions are embraced.\vspace{1.5mm}

\noindent {\textbf{Example 12.1.}} If $V=(2,1,2,2,2,3,1,3,3,2,1,1,2,3)$\hspace{0.25mm}, then the component decomposition  is $((2,1,2),(2),(2),(3,1,3),(3),(2,1,1),(2),(3))$.\vspace{1.5mm}

\noindent {\textbf{Example 12.2.}}  Suppose $V=(3,3,3,3,2,2,1,998)$\hspace{0.25mm}. Then the component decomposition  is
$((3),(3),(3),(3),(2),(2,1,998))$.\vspace{1.4mm}

\noindent {\large{\textbf{13. The  contraction vector.}}}\vspace{1mm}

 \noindent  Let $V=(v_1, \cdots, v_{J+1}) \in M^*$,
 and choose\vspace{0.2mm}

\noindent\hspace{1.5mm} $r_0=0 < r_1 < \cdots <  r_{I+1}=J+1$ such that $r_{j+1}=next(r_j)$ for $0 \le j \le I$.\vspace{0.4mm}

 \noindent We call $r_0, \cdots, r_{I+1}$ the $r$-indexes of $V$. Let $\tau$ be the distance function of $V$ and $\alpha=\delta(V)+1=\tau(J+1)+1$. We define the contraction vector $\pi(V)$ \mbox{of $V$} in the following way:\vspace{0.5mm}

\noindent (13.1) $\pi(V)=(v_1^*, \cdots, v_I^*,v_{I+1}^*)$  where $v^*_{j+1}=\tau(r_{j+1})-\tau(r_j)$  for $0 \le j < I$,

\noindent \hspace{12mm}and $v^*_{I+1}=\tau(r_{I+1})-\tau(r_I)+1=\tau(J+1)-\tau(r_I)+1=\alpha-\tau(r_I)$.\vspace{0.5mm}

\noindent In Section 81  we  prove that $\pi(V) \in M$.\vspace{2mm}

\noindent {{\textbf{Observation 13.1.}}   Suppose $(G_1, \cdots, G_{I+1})$ is the component decomposition of $V$. Then   $\pi(V)=(\delta(G_{1}), \cdots ,\delta(G_{I}),\delta(G_{I+1})+1)$.\vspace{0.5mm}

 \noindent {{\textbf{Proof.}}}  First, we choose indexes $r_0=0 < r_1 < \cdots <  r_{I+1}=J+1$ such \mbox{that $G_{j+1}=(v_{r_j+1}, \cdots, v_{r_{j+1}})$  for $0 \le j \le I$.}
 Since $G_{j+1}$ is a proper odd  component, then Observation 11.2 implies that $r_{j+1}=next(r_j)$ \mbox{for $0 \le j \le I$.} Hence, $r_0 ,r_1, \cdots, r_{I+1}$ are the $r$-indexes of $V$. Then according to (10.2) \mbox{and (13.1)} we get that
  $v^*_{j+1}=\tau(r_{j+1})-\tau(r_j)=\delta(G_{j+1})$  \mbox{for $0 \le j < I$,}
\mbox{and $v^*_{I+1}=\tau(r_{I+1})-\tau(r_I)+1=\delta(G_{I+1})+1$.}\vspace{2mm}

\noindent {\textbf{Observation 13.2.}} Suppose $V=(v_1, \cdots, v_J,v_{J+1}) \in M^*$ where $v_i >1$ \mbox{for $1 \le i \le J$.} Then $\pi(V)=(v_1-1,v_2-1 \cdots, v_J-1,v_{J+1})$.\vspace{0.5mm}

\noindent {\textbf{Proof.}} By Observation 11.4 we get that $V=((v_1), \cdots, (v_J),(v_{J+1}))$
where the components of the component decomposition of $V$ are embraced. Hence,\vspace{0.3mm}

 \hspace{2.6mm} $\pi(V)=(\delta(v_1), \cdots, \delta(v_J),\delta(v_{J+1})+1)=(v_1-1, \cdots, v_J-1,v_{J+1})$\vspace{0.5mm}

 \noindent where the first equality follows from Observation 13.1.\vspace{1.5mm}

 \noindent In the  examples the components of the  decompositions  are embraced.\vspace{1.5mm}

\noindent {\textbf{Example 13.3.}} If $V=((2,1,2),(2),(2),(3,1,3),(3),(2,1,1),(2),(3))$, then\vspace{0.33mm}

 \noindent \hspace{10mm} $\pi(V) =(\delta(2,1,2), \delta(2), \cdots, \delta(2),\delta(3)+1)=(2,1,1,4,2,1,1,3)$.\vspace{2mm}

 \noindent {\textbf{Example 13.4.}} Suppose $V=((3),(3),(3),(3),(2),(2,1,998))$\hspace{0.25mm}.
 Then\vspace{0.4mm}

 \noindent \hspace{5mm} $\pi(V) =(\delta(3), \delta(3), \delta(3), \delta(3), \delta(2), \delta(2,1,998)+1)=(2,2,2,2,1,999)$.\vspace{2mm}

 \noindent {\textbf{Example 13.5.}} Suppose $V=((2),(2),(2),(2,1,999))$. Then\vspace{0.4mm}

 \noindent \hspace{15mm} $\pi(V) =(\delta(2), \delta(2), \delta(2),  \delta(2,1,999)+1)=(1,1,1,1000)$.\vspace{2mm}

\noindent {\textbf{Example 13.6.}} Suppose  $V=((3),(4),(2),(4,1,0))$\,. Then\vspace{0.3mm}

\noindent \hspace{22mm} $\pi(V) =(\delta(3), \delta(4),\delta(2), \delta(4,1,0)+1)=(2,3,1,3)$.\vspace{2mm}

\noindent {\textbf{Example 13.7.}} Suppose  $V=((2),(3,1,3))$\,. Then

\hspace{30mm}$\pi(V) =(\delta(2),  \delta(3,1,3)+1)=(1,5)$.\vspace{2mm}

\noindent {\textbf{Example 13.8.}} If $V=(3,2,2,5,3,2,2,3)$, then Observation 13.2 implies that $\pi(V) =(3-1,2-1,2-1, \cdots, 2-1,3)=(2,1,1,4,2,1,1,3)$.\vspace{2mm}

 \noindent {\large{\textbf{14. Admissible $p$\,-\,values and the set $M_p^+$.}}}\vspace{1mm}

\noindent \noindent Suppose $V=(v_1, \cdots, v_{J+1}) \in M$ and $p\ge 0$, Let $\rho_0=0$. For $0 \le i \le J$  let \vspace{0.5mm}

\noindent (14.1)\hspace{5.7mm} $\rho_{i+1}=\rho_i+v_{i+1}$ if $i$ is even, and $\rho_{i+1}=\rho_i-v_{i+1}$ if $i$ is odd.\vspace{0.5mm}

\noindent Then $\rho_0=0$, $\rho_1=\rho_0+v_1=v_1, \rho_2=v_1-v_2,
\rho_3=v_1-v_2+v_3$, $\cdots$.  We \mbox{call $\rho_0, \rho_1, \cdots, \rho_{J+1}$} the alternating parameters of $V$. Moreover,  we call $p$  an admissible $p$\,-\,value of $V$  if there exists  $t$ such that\vspace{0.5mm}

\noindent  (14.2)\hspace{5.7mm}$1 \le t \le J$ where $t$ is odd, $\rho_i>0$ for $1 \le i \le t$, and $\rho_t \ge p+1$.\vspace{0.4mm}

\noindent We also call $(v_1, \cdots, v_t)$ an admissible start vector of $V$ with respect to $p$
\mbox{if (14.2) is fulfilled.} Next, we define\vspace{0.5mm}

\noindent (14.3) \hspace{6mm}$M^+_p=\{V \in M: p$ is an admissible $p$\,-\,value of $V$\} for  $p \ge 0$.\vspace{0.5mm}

 \noindent If $V \in M$ and $p \ge 0$, we note that  $V \in M^+_p$ if and only if there exists an admissible start vector of $V$ with respect to $p$. \vspace{0.4mm}

\noindent  For example, we get that  $p=2$ is an admissible $p$\,-\,value of $V=(2,1,2,4,1,2)$ since
 $\rho_1=2, \rho_2=2-1=1,
  \rho_3=2-1+2=3$. Hence, $V \in M_2^+$.\vspace{0.5mm}

\noindent In Section 84  we  prove  the following result:\vspace{0.7mm}

\noindent (14.4) \hspace{18.5mm}If $V \in M^+_p$  where $p > 0$\hspace{0.25mm}, then  $\pi(V) \in M^+_{p-1}$\hspace{0.25mm}.\vspace{2mm}

\noindent {\textbf{Observation 14.1.}} Suppose $p \ge 0$ and $V \in M$. If  $V$ starts with a coordinate larger than $p$, then $\rho_1 \ge p+1$ and $V \in M_p^+$.\vspace{0.5mm}

\noindent {\textbf{Proof.}} This is trivial.\vspace{2mm}

\noindent {\textbf{Observation 14.2.}} $M_0^+=M$.\vspace{0.5mm}

\noindent {\textbf{Proof.}}  Let $p=0$. If $V=(v_1, \cdots, v_{J+1}) \in M$,  then $v_1 \ge 1 \ge p+1$. Hence,  Observation 14.1 implies that $V\in M_0^+$. By (14.3) we get that $M_0^+ \subset M$.\vspace{20mm}

\noindent {\textbf{Observation 14.3.}}  If  $V=(v_1, \cdots ,v_{J+1}) \in M_p^+$ where $p \ge 0$, then $V \in M_p$.\vspace{0.5mm}

\noindent {\textbf{Proof.}}  If  $(v_1,v_2, \cdots ,v_J,v_{J+1}) \in M_p^+$ where $p \ge 0$\hspace{0.2mm},  then there exists  an odd index $t$ such that $1 \le t \le J$
and $\rho_t\ge p+1$\hspace{0.12mm}. Hence,\vspace{0.5mm}

\noindent \hspace{24mm} $v_1+v_3+v_5+ \cdots +v_J \ge v_1+v_3+v_5+ \cdots +v_t$\vspace{0.5mm}

\noindent \hspace{24mm} $ \ge v_1-v_2+v_3-v_4 + v_5- \cdots +v_{\hspace{0.2mm}t} = \rho_t \ge p+1$\hspace{0.2mm}.\vspace{2mm}

\noindent {\textbf{Observation 14.4.}} Suppose $V \in M_p^+$ where $p >0$. Then $V \in M^*$.\vspace{0.5mm}

\noindent {\textbf{Proof.}} Suppose $V=(v_1, \cdots, v_{J+1})$. Since $V \in M_p^+ \subset M$, it is sufficient to prove that $v_1>1$. Suppose $v_1=1$. Then $\rho_1=v_1=1$.
If $J >1$, then $v_2 >0$ and $\rho_2=\rho_1-v_2=v_1-v_2 \le 0$. Since $p >0$, then $p$ is not an admissible $p$\,-\,value of $V$. Hence,  $V \not\in M_p^+$. This is a contradiction.\vspace{2mm}

\noindent {\textbf{Observation 14.5.}} Suppose $V \in M_p^+$ and $p \ge 0$. Let $V_p=V$ \mbox{and $V_{i-1}=\pi(V_i)$} for $1 \le i \le p$.
Then  $V_i \in M_{i}^+$  \mbox{for $0 \le i \le p$.}\vspace{0.5mm}

\noindent {\textbf{Proof.}} If $i=p$\hspace{0.2mm}, this is trivial. Next, suppose $V_i \in M_{i}^+$  where $1 \le i  \le p$\hspace{0.2mm}. Since $i>0$\hspace{0.2mm}, then (14.4) implies that  $V_{i-1}=\pi(V_i) \in M_{i-1}^+$.\vspace{2mm}

\noindent {\large{\textbf{15. The distance vector.}}}\vspace{1mm}

\noindent We will define the distance vector $D(V)$ of $V=(v_1, \cdots, v_J,v_{J+1}) \in M^*$. Then $J \ge 1$ is odd, $v_1>1$, $v_i \ge 1$ \mbox{for $2 \le i \le J$,} and $v_{J+1} \ge 0$\hspace{0.2mm}.\vspace{0.5mm}

 \noindent First, we suppose  $v_i>1$ for $1 \le i \le J$. Then we   let $D(V)=\O$.
 Otherwise,  there exists $i$ such that $1 < i \le J$ and $v_i=1$.  Hence, we can in this case choose
 indexes $c_0, \cdots, c_{\gamma}$ such that

  \noindent (15.1) \hspace{34mm} $c_0=0< c_1 < \cdots <c_{\gamma} \le  J$,
 \vspace{0.5mm}

\noindent (15.2) \hspace{1mm}$c_{i+1}$ is the least index $>c_i+1$ such that $v_{c_{i+1}}=1$ for $0 \le i < \gamma$,\vspace{0.3mm}

 \noindent (15.3)  \hspace{34mm}$v_c>1$ for $c_{\gamma}+1 < c \le J$.\vspace{0.6mm}

 \noindent Since $c_0=0$, then we get from (15.2) that

 \noindent (15.4) \hspace{20mm}$c_{1}$ is the least index $>1$ such that $v_{c_{1}}=1$.\vspace{0.3mm}

\noindent We call
$c_0, \cdots , c_{\gamma}$  the $c$\,-\,indexes of $V$\hspace{0.2mm}.
Then we define\vspace{0.5mm}

\noindent (15.5) \hspace{20.5mm} $D(V)=(\tau(c_1), \cdots, \tau(c_{\gamma}))=(d_1, \cdots, d_{\gamma})$\vspace{0.4mm}

\noindent where $d_i=\tau(c_i)$ for $1 \le i \le \gamma$. Since $v_1>1$, then we get by  (10.5)  that\vspace{0.5mm}

\noindent (15.6) \hspace{10mm}$0 < d_1 \le d_2 \le \cdots \le d_{\gamma} \le \tau(J) \le \alpha$ where $\alpha=\delta(V)+1$\,.\vspace{0.5mm}

\noindent   In the following examples   the coordinates corresponding to the  positive $c$\,-\,indexes are overlined.\vspace{1mm}

 \noindent {\textbf{Example 15.1.}} Suppose
$V=(2,\overline{1},{2},{2},{2},3,\overline{1},{3},{3},2,\overline{1},{1},
{2},3)$. Then $c_1=2$, $c_2=7$ and $c_3=11$ are the  positive $c$\,-\,indexes
of $V$, and
 $\alpha=\delta(V)+1=15$.  Hence, $D(V)=(\tau(c_1),\tau(c_2),\tau(c_3))=(\tau(2),\tau(7),\tau(11))=(1,6,11)$

\noindent since $\tau(2)=\delta(2,1)=1$, $\tau(7)=\delta(2,1,2,{2},{2},3,{1})=6$ and
$\tau(11)=11$\hspace{0.25mm}.\vspace{2mm}

\noindent {\textbf{Example 15.2.}} Suppose
$V=({3},{3},{3},{3},{2},2,\overline{1},998)$. Then  $c_1=7$   is the only  positive $c$\,-\,index of $V$ and  $\alpha=\delta(V)+1=1008$\hspace{0.25mm}.
 Moreover, we get that
 $D(V)=(\tau(c_1))=(10)$  since $\tau(c_1)=\tau(7)=\delta({3},{3},{3},{3},{2},2,{1})=10$\hspace{0.35mm}.\vspace{2mm}

\noindent {\textbf{Example 15.3.}} Suppose
 $V=({2},{2},{2},{2},\overline{1},999)$  Then  $c_1=5$   is the only  positive $c$\,-\,index of $V$, and $\alpha=\delta(V)+1=1003$\hspace{0.25mm}.
 Moreover, we get \mbox{that
 $D(V)=(\tau(c_1))=(4)$} since $\tau(c_1)=\tau(5)=\delta({2},{2},{2},{2},{1})=4$\hspace{0.25mm}.\vspace{2mm}

 \noindent {\textbf{Example 15.4.}} Suppose
$V=(3,4,2,4,\overline{1},0)$. Then  $c_1=5$   is the only  positive $c$\,-\,index of $V$, $\alpha=\delta(V)+1=9$
\mbox{and $D(V)=(\tau(c_1))=(\tau(5))=(9)$.}\vspace{1.5mm}

 \noindent {\textbf{Example 15.5.}} Suppose
$V=(2,3,\overline{1},3)$. Then  $c_1=3$   is the only  positive $c$\,-\,index of $V$, $\alpha=\delta(V)+1=6$ \mbox{and  $D(V)=\tau(c_1)=(\tau(3))=(3)$.}\vspace{2mm}

\noindent {\textbf{Example 15.6.}} If $V=(3,2,2,5,3,2,2,3)$, then $D(V)=\O$.\vspace{2mm}

\noindent {\textbf{Example 15.7.}} If $V=(3,4,2,3)$, then $D(V)=\O$.\vspace{2mm}

\noindent \textbf{{\large  {16. Progression parameters.}}}\vspace{1mm}

 \noindent In this section we suppose $\alpha>0$ and\vspace{0.5mm}

\hspace{12mm}$D=(d_1,d_{\hspace{0.2mm}2}, \cdots ,d_{\gamma})$ where\vspace{0.4mm}
 $ 0<d_1 \le d_{\hspace{0.2mm}2} \le \cdots \le d_{\gamma} \le \alpha$.\vspace{0.5mm}

\noindent We note that $\gamma=\#D$ = the number of coordinates of $D$. Let $E=(D,D+\alpha)$. That means, $E=(e_1, \cdots, e_{2\gamma})$ where $e_{i}=d_i$ and $e_{\gamma+i}=d_i+\alpha$ for $1 \le i \le \gamma$.\vspace{0.4mm}

\noindent Let $F=\{m >0: m$ factor of $gcd(\alpha,\gamma)\}$ where $gcd(\alpha,\gamma)$ is the greatest common divisor of $\alpha$ and $\gamma$.\vspace{0.8mm}

\noindent Suppose $m \in F$.  Let $\beta=\frac{\alpha}{m}$  and $r=\frac{\gamma}{m}$. If
$e_{r+i}=d_i + \beta$ for $1 \le i \le \gamma$,  \mbox{then $m$} is a called a progression coefficient
of $D$ with respect to $\alpha$.

\noindent Alternatively, $m$ is called a progression coefficient
of $D$ with respect to $\alpha$ \mbox{if
$(e_{r+1}, \cdots, e_{r+\gamma})=(d_1, \cdots, d_{\gamma})+\beta$.}\vspace{0.8mm}

\noindent By Observation 16.1 we get that $m=1$ is a progression coefficient
of $D$ with respect to $\alpha$. Let $m^*$ be the maximal progression coefficient
of $D$ with respect to $\alpha$. Then $\alpha^*=\frac{\alpha}{m^*}$ and  $\gamma^* =\frac{\gamma}{m^*}$ are called the least progression parameters of $D$ with respect to $\alpha$.\vspace{0.5mm}

 \noindent Let $m_1>m_2> \cdots > m_j=1$ be the factors in $F$ in descending order. Then we test successively
 if $m_i$ is a  progression coefficient
of $D$ with respect to $\alpha$ for $i =1, i=2, \cdots$, until we find $i$ such that\vspace{0.2mm}

\hspace{6mm}$(e_{r+1}, \cdots, e_{r+\gamma})=(d_1, \cdots, d_{\gamma})+\beta$ where $\beta=\frac{\alpha}{m_i}$ and  $r =\frac{\gamma}{m_i}$.\vspace{0.65mm}

\noindent Then $m^*=m_i$ is the maximal progression coefficient
of $D$ with respect to $\alpha$, and $\alpha^*=\frac{\alpha}{m^*}=\beta$ and  $\gamma^* =\frac{\gamma}{m^*}=r$.\vspace{20mm}

\noindent \textbf{{{Observation 16.1.}}} \mbox{$m=1$ is a progression coefficient
of $D$ with respect to $\alpha$.}\vspace{0.5mm}

\noindent \textbf{{{Proof.}}} If $m=1$, then $\beta=\frac{\alpha}{m}=\alpha$  and $r=\frac{\gamma}{m}=\gamma$. It is sufficient to prove that $e_{r+i}=d_i + \beta$ for $1 \le i \le \gamma$. This is true since\vspace{0.5mm}

\hspace{22mm}$e_{r+i}=e_{\gamma+i}=d_i+\alpha=d_i+\beta$ for $1 \le i \le \gamma$.\vspace{2mm}

\noindent {\textbf{Observation 16.2.}}  Suppose $gcd(\alpha,\gamma)=1$.

\noindent a) $m^*=1$ is the maximal
progression coefficient
of $D$ with respect to $\alpha$.\vspace{0.5mm}

\noindent b) $\alpha^*=\frac{\alpha}{1}=\alpha$  and  $\gamma^*=\frac{\gamma}{1}=\gamma$ are the least
progression parameters

\noindent \hspace{5mm}of $D$ with respect to $\alpha$.\vspace{0.6mm}

\noindent \textbf{{{Proof.}}}  If $gcd(\alpha,\gamma)=1$, then $F=\{m >0: m$ factor of $gcd(\alpha,\gamma)\}=\{1\}$. By Observation 16.1 we get that $m^*=1$ is the only
progression coefficient
\mbox{of $D$} with respect to $\alpha$. Hence, a) is true.  Moreover, b) follows from a).\vspace{2mm}

\noindent In particular, if $\gamma=1$, then $gcd(\alpha,\gamma)=1$ and  Observation 16.2 b) implies \mbox{that
$\alpha^*=\frac{\alpha}{1}=\alpha$  and  $\gamma^*=\frac{\gamma}{1}=\gamma=1$.}\vspace{2mm}

\noindent \textbf{{{Example 16.3.}}}  Let $D=(d_1, \cdots, d_6)=(2,4,8,10,14,16)$ and $\alpha=18$\hspace{0.2mm}. Then $\gamma=\#D=6$ and $gcd(\alpha,\gamma)=6$. In this case $F=\{6,3,1\}$ and\vspace{0.5mm}

\noindent \hspace{3mm} $E=(D,D+\alpha)=(e_1, \cdots, e_{12})=(2,4,8,10,14,16, 20,22,26,28,32,34)$.\vspace{0.5mm}


\noindent Suppose $m=6$. Let $\beta=\frac{\alpha}{m}=3$ and  $r =\frac{\gamma}{m}=1$. Then\vspace{0.5mm}

\noindent \hspace{4mm} $(e_{r+1}, \cdots, e_{r+6})=(e_2, \cdots, e_7)=(4,8,10,14,16, 20) \ne (d_1, \cdots, d_{6})+\beta$.\vspace{0.5mm}

\noindent Suppose $m=3$. Let $\beta=\frac{\alpha}{m}=6$ and $r =\frac{\gamma}{m}=2$. Then\vspace{0.5mm}

\noindent \hspace{3.9mm}$(e_{r+1}, \cdots, e_{r+6})=(e_3, \cdots, e_8)=(8,10,14,16, 20,22) = (d_1, \cdots, d_6)+\beta$.\vspace{0.5mm}

\noindent Hence, $m^*=3$ is the maximal progression coefficient of $D$ with respect to $\alpha$, and $\alpha^*=\frac{\alpha}{m^*}=6$  and  $\gamma^*=\frac{\gamma}{m^*}=2$.\vspace{2mm}

\noindent {\textbf{Example 16.4.}} Let $D=(d_1,d_2,d_3)=(1,6,11)$ and $\alpha=15$\hspace{0.25mm}. Then we get that \mbox{$\gamma=\#D=3$} and
$E=(D,D+\alpha)=(e_1, \cdots, e_{6})=(1,6,11, 16,21,26)$.  The positive factors of $gcd(\alpha,\gamma)=3$ are $3$  and $1$\hspace{0.15mm}.
 In this case $F=\{3,1\}$. \noindent  \vspace{0.5mm}

 \noindent Suppose $m=3$. Let $\beta=\frac{\alpha}{m}=\frac{15}{3}=5$ and
 $r =\frac{\gamma}{m}=\frac{3}{3}=1$. Then\vspace{0.5mm}

\noindent \hspace{16mm}$(e_{r+1}, e_{r+2}, e_{r+3})=(e_2, e_3,e_4)=(6,11,16)
= (d_1, d_2, d_{3})+\beta$.\vspace{0.5mm}

\noindent We conclude that $m^*=3$ is the maximal progression coefficient of $D$ with respect to $\alpha$, and $\alpha^*=\frac{\alpha}{m^*}=5$  and  $\gamma^*=\frac{\gamma}{m^*}=1$.\vspace{2mm}

\noindent {\large{\textbf{17. Progression parameters of distance vectors.}}}\vspace{1mm}

\noindent {\textbf{Example 17.1.}} Suppose $Q_1=(2,{1},{2},{2},{2},3,{1},{3},{3},2,{1},{1},
{2},3)$. Then we get by Example 15.1  that $D(Q_1)=(1,6,11)$ \mbox{and $\alpha_1=\delta(Q_1)+1=15$.}\vspace{0.5mm}

\noindent  By Example 16.4 we also get  that $\alpha_1^*=5$ and $\gamma_1^*=1$ are the least progression parameters
of $D(Q_1)$ with respect to $\alpha_1$.\vspace{2mm}

\noindent {\textbf{Example 17.2.}} Suppose $Q_2=({3},{3},{3},{3},{2},2,{1},998)$. Then we get according to Example 15.2   that $D(Q_2)=(10)$ \mbox{and $\alpha_2=\delta(Q_2)+1=1008$}.
\mbox{Then $\gamma_2=\#D(Q_2)=1$.} Observation 16.2 implies that  $\alpha_2^*=\alpha_2=1008$
\mbox{and $\gamma_2^*=1$} are the least progression parameters
of $D(Q_2)$ \mbox{with respect to $\alpha_2$.}\vspace{2mm}

\noindent {\textbf{Example 17.3.}} Let $Q_1=({2},{2},{2},{2},{1},999)$. Then we get by Example 15.3  that $D(Q_1)=(4)$ \mbox{and $\alpha_1=\delta(Q_1)+1=1003$}.
Then $\gamma_1=\#D(Q_1)=1$. Hence, Observation 16.2 implies that  $\alpha_1^*=\alpha_1=1003$ and $\gamma_1^*=1$ are the least progression parameters
of $D(Q_1)$ with respect to $\alpha_1$.\vspace{2mm}

\noindent {\textbf{Example 17.4.}} Let $Q_2=(3,4,2,4,{1},0)$. Then we get by Example 15.4   that $D(Q_2)=(9)$ and $\alpha_2=\delta(Q_2)+1=9$.
Then $\gamma_2=\#D(Q_2)=1$. Hence, Observation 16.2 implies that  $\alpha_2^*=\alpha_2=9$ and $\gamma_2^*=1$ are the least progression parameters
\mbox{of $D(Q_2)$} with respect to $\alpha_2$.\vspace{2mm}

\noindent {\textbf{Example 17.5.}} Suppose $Q_1=(2,3,{1},3)$. Then  Example 15.5   implies that $D(Q_1)=(3)$ \mbox{and $\alpha_1=\delta(Q_1)+1=6$}.
Then $\gamma_1=\#D(Q_1)=1$. Hence, Observation 16.2 implies that  $\alpha_1^*=\alpha_1=6$ and $\gamma_1^*=1$ are the least progression parameters
of $D(Q_1)$ with respect to $\alpha_1$.\vspace{2mm}

\hspace{45mm} {\large{\textbf{PART 2.}}  \vspace{ 1 mm}}

\noindent In Section 19 we describe the main case and explain how every case can be reduced to this case. In Section 22 we describe the determination of the minimal periods  in the main case.
In Section $20$ and $21$  we define
parameters. Section $24, 25$ and  $26$ contain examples.  In Section 27 we define and plot weight parameters.
These plots bear a resemblance to plots in soliton theory.\vspace{2mm}

\noindent {\large{\textbf{18. Notation.}}}\vspace{1mm}

\noindent Suppose $Q=(q_1, \cdots, q_J,e_0) \in M_p$ where $p \ge 0$. Let\vspace{0.5mm}

 $A=A(Q)=1_{q_1}0_{q_2} \cdots 1_{q_J}0_{e_0}$,  $k=w(A)-(p+1)$ and $n=length(A)$.\vspace{0.5mm}

\noindent Then we define   $A^{\,\infty}_p(Q)=A^{\,\infty}$ where $A^{\,\infty}$ is  generated from $A=A(Q)$  by the symmetric shift register with respect to
 $k$, $p$ and $n$. We note that $k \ge 0$ since $w(A)= q_1+q_3+ \cdots+q_J \ge p+1$,
 $w(A)=k+p+1$ \mbox{and $A$ starts with 1.} We call $A^{\,\infty}_p(Q)$ the shift symmetric bit string generated by $Q$ and $p$.\vspace{0.7mm}

\noindent  Since $w(A)=k+p+1$, then we get according to (6.2) and (6.3) that\vspace{0.4mm}

 \noindent (18.1)
 \hspace{28mm}$V(A^{\,\infty})$ has an even vector period,\vspace{0.5mm}

 \noindent (18.2) \hspace{6mm}$sum(V(A^{\,\infty}),j)$ is the
minimal period of $A^{\,\infty}$ if  $j$ is the  least

  \hspace{11.7mm}even vector period of $V(A^{\,\infty})$.\vspace{0.5mm}

\noindent By Observation 14.3 we also get that $A^{\,\infty}_p(Q)$ is well-defined if $Q \in M_p^+$.\vspace{2mm}

 \noindent {{\textbf{Example 18.1.}}} Suppose  $Q=(3,4,2,4,1,0)$ and $p=2$\hspace{0.25mm}.
 Then $Q \in M_p$ and\vspace{0.5mm}

\hspace{20mm} $A=A(Q)=1_30_41_20_41_10_{\hspace{0.2mm}0}=11100001100001$\hspace{0.25mm}.

\noindent Let $k=w(A)-(p+1)=6-(2+1)=3$ and $n=length(A)=14$.
Then $A^{\,\infty}_2(Q)=A^{\,\infty}$ where $A^{\,\infty}$ is  generated from $A$  by the symmetric shift register with respect to
 $k=3$, $p=2$ and $n=14$.\vspace{2mm}

\noindent {\large{\textbf{19. The main case.}}}\vspace{1mm}

 \noindent Let $A^{\,\infty}$ be generated from $A $  by the  symmetric shift register with respect \mbox{to $k$, $p$ and $n$.} We note that $n=length(A)$. Suppose\vspace{0.5mm}

 \noindent (19.1) \hspace{12mm}$w(A)=k+p+1$\hspace{0.2mm}, $A$ starts with 1 and
 $V(A) \in M_p^+$.\vspace{0.5mm}

 \noindent   This is the main case.   Let $Q=V(A)$. According to Observation 7.1 we get that $A=A(Q)$. Hence, $A^{\,\infty}$ is generated from $A=A(Q) $  by the  symmetric shift register with respect \mbox{to $k=w(A)-(p+1)$,} $p$ and $n=length(A)$. By Section 18 we get that\vspace{0.5mm}

\noindent (19.2)\hspace{24mm} $A^{\,\infty}=A_p^{\,\infty}(Q)$ where $Q=V(A)\in M_p^+$.\vspace{0.5mm}

\noindent  In Part 12  we describe how every case can be reduced to the case where (19.1) is satisfied.  Therefore by  (19.2) every case can be reduced to $A^{\,\infty}=A_p^{\,\infty}(Q)$ where $Q\in M_p^+$. In Section 22 we describe how the minimal periods of $A_p^{\,\infty}(Q)$ can be determined for $Q\in M_p^+$ and $p \ge 0$.\vspace{2mm}

\noindent \textbf{{\large 20. Cyclic parameters.}}\vspace{1mm}

\noindent If $V=(v_1, \cdots, v_{r})$ where $r >0$, let $\psi(V)=(v_{2}, \cdots, v_{r},v_1)$\hspace{0.25mm}. It is easily seen that  $\psi^{\hspace{0.1mm}r}(V)=V$ and  $\psi^{\hspace{0.2mm}j}(V)=(v_{j+1}, \cdots, v_{r},v_1, \cdots, v_{j})$ for $1 \le j<r$\hspace{0.05mm}.\vspace{0.5mm}

 \noindent We suppose $V=(v_1, \cdots, v_{r})$ where $r >0$ is even. Since $\psi^r(V)=V$,
 there exists a  least even positive integer $j \le r$ such that $\psi^{\hspace{0.2mm}j}(V)=V$.\vspace{0.5mm}

 \noindent   Let   $\zeta=sum(V,j)=v_1+ \cdots +v_j$\hspace{0.25mm}. We call $j$ \mbox{and $\zeta$}  the cyclic parameters \mbox{of $V$.} It is easily  proved that $j$ is an even positive factor of $r=\#V$.\vspace{1.5mm}

\noindent {\textbf{Example 20.1.}} Let $V=(2,1,1,4,2,1,1,4)$\hspace{0.25mm}.
Since $j$ must be an even positive factor of $\#V=8$, there are three possibilities, $j=2$\hspace{0.25mm}, $j=4$ \mbox{or $j=8$\hspace{0.25mm}.}
We get that $\psi^2(V)=(1,4,2,1,1,4,2,1) \ne V$ and $\psi^4(V)=V$.\vspace{0.5mm}

\noindent  Then  $j=4$ is the least even positive integer such that $\psi^{\hspace{0.2mm}j}(V)=V$. Hence,   $j=4$ and $\zeta=sum(V,4)=2+1+1+4=8$ are  the cyclic parameters of $V$.\vspace{1mm}

\noindent {\textbf{Example 20.2.}} Suppose $V=(1,1,1,1001)$\hspace{0.25mm}.
Since $j$ must be an even positive factor of $\#V=4$\hspace{0.25mm}, there are two possibilities, $j=2$ or $j=4$.  We note that $\psi^2(V)=(1,1001,1,1) \ne V$. Hence,  $j=4$ is the least even positive integer such that $\psi^{\hspace{0.2mm}j}(V)=V$.  Then  $j=4$ and $\zeta=sum(V,4)=1004$ are  the cyclic parameters of $V$.\vspace{1.4mm}

\noindent {\textbf{Example 20.3.}} Suppose $V=(1,6)$. Since $\#V=2$, then   $j=2$ is the least even positive integer such that $\psi^{\hspace{0.2mm}j}(V)=V$. Hence,\vspace{0.3mm}

\noindent\hspace{6mm}$j=2$ and $\zeta=sum(V,2)=1+6=7$
are  the cyclic parameters of $V$.\vspace{2mm}

\noindent {\large{\textbf{21. Dynamical parameters.}}}\vspace{0.5mm}

\noindent If $Q=(q_1, \cdots,q_J,e_0)$, then $Q^*=(q_1, \cdots,q_J,e_0+1)$ is the extension of $Q$.\vspace{0.5mm}

\noindent
  Suppose $Q_p \in M^+_p$ where $p \ge 0$. If $p>0$,  let\vspace{0.4mm}

\noindent   \hspace{27mm}$Q_{p-1}=\pi(Q_p), Q_{p-2}=\pi(Q_{p-1}), \cdots, Q_0=\pi(Q_1)$.\vspace{0.6mm}

\noindent By Observation 14.5 we get that $Q_{p-1}, Q_{p-2}, \cdots, Q_0$ are well-defined.
Next, we  define the dynamical parameters $j_0, j_1, \cdots, j_p$ and $\zeta_0,\zeta_1, \cdots ,\zeta_p$ \mbox{of $Q_p$} with respect to $p$ by an inductive process as follows. Let\vspace{0.8mm}

\noindent (21.1) \hspace{3mm} $j_0$ and $\zeta_{\hspace{0.2mm}0}$ be the cyclic parameters of the extension $Q_0^*$   of $Q_0$.\vspace{0.5mm}

\noindent Suppose $j_{i-1}$ and $\zeta_{i-1}$ where $1 \le i \le p$, are determined. Then $j_{i}$ and $\zeta_{i}$ are  determined as follows:\vspace{0.8mm}

\noindent (21.2)
\hspace{2mm}If $D(Q_i)=\O$, let $j_{i}=j_{i-1}$ and $\zeta_i=\zeta_{i-1}+j_{i}=\zeta_{i-1}+j_{i-1}$.\vspace{0.8mm}

\noindent  (21.3) \hspace{2mm}If $D(Q_i)\ne \O$, let $\alpha^*_i$ and $\gamma^*_i$ be the least progression parameters

\hspace{8mm}of $D(Q_i)$ with respect to $\alpha_i=\delta(Q_i)+1$,\hspace{10mm}

\hspace{8mm}let $x$ and $y$ be the least positive integer solution
of $x\alpha^*_i=y\zeta_{i-1}$,\hspace{10mm}

 \hspace{8mm}let $j_i=2x\gamma_i^*+yj_{i-1}$ and  $\zeta_i=y\zeta_{i-1}+j_i$.\vspace{0.5mm}

 \noindent If $p=0$, then we note that $j_0$ and $\zeta_0$ are the only dynamical parameters \mbox{of $Q_p$} with respect to $p$.\vspace{2mm}

 \noindent \textbf{{\large 22. The determination of the periods.}}\vspace{1mm}

\noindent
Suppose $Q_p \in M^+_p$ where $p \ge 0$.
If $p>0$,  let\vspace{0.4mm}

\noindent (22.1)  \hspace{14mm}$Q_{p-1}=\pi(Q_p), Q_{p-2}=\pi(Q_{p-1}), \cdots, Q_0=\pi(Q_1)$.\vspace{0.4mm}

\noindent Let $j_0, j_1, \cdots, j_p$ and $\zeta_0,\zeta_1, \cdots, \zeta_p$ be the dynamical parameters of $Q_p$ with respect to $p$.  By Observation 14.5 we get that $Q_i \in M_i^+$ for $0 \le i \le p$. Suppose $A_i^{\,\infty}=A^{\,\infty}_i(Q_i)$ where  $0 \le i \le p$. By Proposition 41.3 we get that\vspace{0.4mm}

\noindent (22.2)\hspace{19.5mm}$j_i$ is the least even vector period of $V(A_i^{\,\infty})$,\vspace{0.5mm}

\noindent (22.3) \hspace{29.3mm}$\zeta_i$ is the minimal period  of $A_i^{\,\infty}$.\vspace{2mm}

\noindent \textbf{{\large 23. Solving equations.}}\vspace{0.5mm}

\noindent  As previously, $gcd(\alpha,\beta)$ denotes the greatest common divisor of $\alpha$ and $\beta$.
 Suppose we will  determine the least positive integers $x$ and $y$ such that $x \cdot 30 = y \cdot 65$\hspace{0.25mm}.
We divide both sides of the equation by $5=gcd(30,65)$ and get $x \cdot 6 = y \cdot 13$\hspace{0.25mm}. The answer is $x=13$ and $y=6$\hspace{0.25mm}.
 It is easily seen that this method can be generalized as described by the following observations.\vspace{2mm}

\noindent {\textbf{Observation 23.1.}} If $\alpha>0$ and $\beta>0$\hspace{0.35mm}, let  $g=gcd(\alpha,\beta)$\hspace{0.25mm}. Then $x=\frac{\beta}{g}$ \mbox{and $y=\frac{\alpha}{g}$}  are the least positive integers  satisfying  $x\alpha=y\beta$\hspace{0.25mm}.\vspace{1.6mm}

\noindent {\textbf{Observation 23.2.}} If $gcd(\alpha,\beta)=1$ where $\alpha>0$ and $\beta>0$\hspace{0.35mm}, then $x=\beta$ and $y=\alpha$ are the least positive integers satisfying  $x\alpha=y\beta$\hspace{0.25mm}.\vspace{1.6mm}

\noindent {\textbf{Example 23.3.}} Suppose $\alpha=1008$ and $\beta=1013032$.\vspace{0.01mm}
 Then\vspace{0.4mm}

\noindent \hspace{6mm}$x=\frac{\beta}{g}=126629$ and $y=\frac{\alpha}{g}=126$ where $g=gcd(1008,1013032)=8$,\vspace{0.4mm}

\noindent are the least positive integers  satisfying the equation $x\alpha=y\beta$\hspace{0.25mm}.\vspace{2mm}

 \noindent {\large{\textbf{24. Example in the case  $p=0$.}}}\vspace{0.5mm}

 \noindent  Suppose  $Q_0=(2,1,1,4,2,1,1,3)$ and $p=0$.  We note that\vspace{0.6mm}

\noindent\hspace{3mm}$Q_0^*=(2,1,1,4,2,1,1,3+1)=(2,1,1,4,2,1,1,4)$ is the extension of $Q_0$.\vspace{0.3mm}

\noindent By Example 20.1  we get that
$j_0=4$ and $\zeta_0=8$ are the cyclic parameters of $Q_0^*$.
Hence, $j_0=4$ and $\zeta_0=8$ are the  dynamical parameters of $Q_0$ with respect to $p=0$.
 Let
$A^{\,\infty}_0=A^{\,\infty}_0(Q_0)$. By (22.2)  we get that $j_0=4$ is the least even vector period of $V(A_0^{\,\infty})$. By (22.3) we get   $\zeta_0=8$ is the minimal period of $A_0^{\,\infty}$. We note that  $A^{\,\infty}_0=A^{\,\infty}_0(Q_0)$ is  generated from\vspace{0.4mm}

\hspace{15.8mm}$A_0=A(Q_0)=1_2 0_1 1_1 0_4 1_2 0_1 1_1 0_3=110100001101000$\vspace{0.4mm}

\noindent by the symmetric shift register with respect to\vspace{0.4mm}

\noindent \hspace{3mm} $k=w(A_0)-(p+1)=6-(0+1)=5$, $p=0$ and $n=length(A_0)=15$.\vspace{2mm}

\noindent \textbf{{\large 25. Examples in the case  $p =1$.}}\vspace{1mm}

  \noindent {\textbf{Example 25.1.}} Suppose $Q_1=(3,2,2,5,3,2,2,3)$ and $p=1$.   According to Example 13.8  we get that
$Q_0=\pi(Q_1)=(2,1,1,4,2,1,1,3)$.

\noindent The extension of $Q_0$ is $Q_0^*=(2,1,1,4,2,1,1,4)$\hspace{0.25mm}. By Example 20.1 we get that $j_0=4$ and $\zeta_{\hspace{0.2mm}0}=8$ are the cyclic parameters of $Q_0^*$.\vspace{0.5mm}

\noindent By Example 15.6 we get that  $D(Q_1)=\O$.  Then (21.2) implies that\vspace{0.3mm}

\hspace{23.5mm}$j_1=j_0=4$ and $\zeta_1=\zeta_0+j_1=8+4=12$.\vspace{1mm}

 \noindent Let $A_0^{\,\infty}=A^{\,\infty}_0(Q_0)$ and $A_1^{\,\infty}=A^{\,\infty}_1(Q_1)$.  By (22.2)  we get that $j_0=4$ and $j_1=4$ are the least even vector periods \mbox{of
 $V(A_0^{\,\infty})$} and $V(A_1^{\,\infty})$. By (22.3)  we get that $\zeta_0=8$ and $\zeta_1=12$ are the minimal periods of
 $A_0^{\,\infty}$ and $A_1^{\,\infty}$.\vspace{2mm}

 \noindent {\textbf{Example 25.2.}} Suppose $Q_1=(2,1,2,2,2,3,1,3,3,2,1,1,2,3)$ and $p=1$.   By Example 13.3 we get that
$Q_0=\pi(Q_1)=(2,1,1,4,2,1,1,3)$.

\noindent The extension of $Q_0$ is $Q_0^*=(2,1,1,4,2,1,1,4)$\hspace{0.25mm}. By Example 20.1 we get that $j_0=4$ and $\zeta_{\hspace{0.2mm}0}=8$ are the cyclic parameters of $Q_0^*$.\vspace{0.7mm}

\noindent By Example 17.1 we get that $D(Q_1) \ne \O$, and  $\alpha^*_1=5$ and $\gamma^*_1=1$ are the least progression parameters  of $D(Q_1)$ with respect to $\alpha_1=\delta(Q_1)+1$. The equation $x\alpha_1^*=y\zeta_0$ is equivalent to $5x=8y$. Hence, the least positive integer  solution is $x=8$ and $y=5$. By (21.3) we get that\vspace{0.5mm}

\noindent $j_1=2x\gamma_1^*+yj_{0}=2\cdot 8 \cdot 1+5 \cdot 4=36$ and $\zeta_1=y\zeta_{0}+j_1=5\cdot 8+36=76$.\vspace{0.5mm}

\noindent Let $A_0^{\,\infty}=A^{\,\infty}_0(Q_0)$ and $A_1^{\,\infty}=A^{\,\infty}_1(Q_1)$. By (22.2)  we get that $j_0=4$ and $j_1=36$ are the least even vector periods \mbox{of
 $V(A_0^{\,\infty})$} \mbox{and $V(A_1^{\,\infty})$,} and \mbox{by (22.3)} that  $\zeta_0=8$ and $\zeta_1=76$ are the minimal periods of
 $A_0^{\,\infty}$ \mbox{and $A_1^{\,\infty}$.}\vspace{2mm}

 \noindent \textbf{{\large 26. Examples in the  case $p =2$.}}\vspace{1mm}

 \noindent {\textbf{Example 26.1.}} Suppose  $Q_2=(3,4,2,4,1,0)$ and $p=2$\hspace{0.25mm}.
By  Example 13.6 and 13.7  we get that
$Q_1=\pi(Q_2)=(2,3,1,3)$ \mbox{and $Q_0=\pi(Q_1)=(1,5)$\hspace{0.25mm}.}\vspace{0.6mm}

\noindent  The extension of $Q_0$ \mbox{is $Q_0^*=(1,6)$\hspace{0.15mm}.} By Example 20.3 we get that
$j_{\hspace{0.2mm}0}=2$ \mbox{and $\zeta_{\hspace{0.3mm}0}=7$} are the cyclic parameters of $Q_0^*$\hspace{0.25mm}.\vspace{0.5mm}

\noindent By Example 17.5   we get that $D(Q_1) \ne \O$, and $\alpha_1^*=6$ and $\gamma_1^*=1$ are the least progression parameters  of $D(Q_1)$ with respect to $\alpha_1=\delta(Q_1)+1$.\vspace{0.5mm}

\noindent The equation $x\alpha_1^*=y\zeta_0$ is equivalent to $6x=7y$. The least positive integer solution is $x=7$ and $y=6$. By (21.3) we get that\vspace{0.5mm}

\hspace{25mm}$j_1=2x\gamma_1^*+yj_{0}=2\cdot 7 \cdot 1+6 \cdot 2=26$,

\hspace{25.4mm}$\zeta_1=y\zeta_{0}+j_1=6\cdot 7+26=68$.\vspace{0.5mm}

\noindent By Example  17.4  we get that $D(Q_2) \ne \O$, and $\alpha_2^*=9$ and $\gamma_2^*=1$ are the least progression parameters  of $D(Q_2)$ with respect to $\alpha_2=\delta(Q_2)+1$. The equation $x\alpha_2^*=y\zeta_1$ is equivalent to $9x=68y$. Since $gcd(9,68)=1$, then the least positive integer solution is $x=68$ and $y=9$. By (21.3) we get that\vspace{0.5mm}

\hspace{25mm}$j_2=2x\gamma_2^*+yj_{1}=2\cdot 68 \cdot 1+9 \cdot 26=370$,

\hspace{25.4mm}$\zeta_2=y\zeta_{1}+j_2=9\cdot 68+370=982$.\vspace{0.5mm}

\noindent
$A_2^{\,\infty}=A^{\,\infty}_2(Q_2)$ is generated from $A_2=A(Q_2)=1_30_41_20_41_10_{\hspace{0.2mm}0}=11100001100001$ by the symmetric shift register  with respect to\vspace{0.3mm}

 \hspace{10mm} $k_{\hspace{0.2mm}2}=w(A_2)-(2+1)=3$, $p_2 =2$
and $n_2=length(A_2)=14$.\vspace{0.5mm}

\noindent $A_1^{\,\infty}=A^{\,\infty}_1(Q_1)$ is generated from $A_1=A(Q_1)=1_20_31_10_3=110001000$ by the symmetric shift register  with respect to\vspace{0.1mm}

\hspace{8mm} $k_{\hspace{0.2mm}1}=w(A_1)-(1+1)=1$, $p_1 =1$
and $n_1=length(A_1)=9$.\vspace{0.4mm}

\noindent $A_0^{\,\infty}=A^{\,\infty}_0(Q_0)$ is generated from $A_0=A(Q_0)=1_10_{\hspace{0.2mm}5}=100000$ by the symmetric shift register  with respect to\vspace{0.1mm}

\hspace{8mm} $k_{\hspace{0.2mm}0}=w(A_0)-(0+1)=0$, $p_0 =0$
and $n_0=length(A_0)=6$.\vspace{0.5mm}

\noindent By (22.2)  we get that $j_0=2$,  $j_1=26$ and $j_2=370$ are the least even vector periods \mbox{of
 $V(A_0^{\,\infty})$}, $V(A_1^{\,\infty})$  and $V(A_2^{\,\infty})$. Moreover, by (22.3)  we get that $\zeta_0=7$, $\zeta_1=68$ and  $\zeta_2=982$ are the minimal periods of
 $A_0^{\,\infty}$, $A_1^{\,\infty}$  and $A_2^{\,\infty}$.\vspace{20mm}

\noindent {\textbf{Example 26.2.}} Suppose  $Q_2=(3,3,3,3,2,2,1,998)$ and $p=2$\hspace{0.25mm}.
 According to  Example 13.4 and 13.5  we get that\vspace{0.3mm}

\hspace{4.8mm}$Q_1=\pi(Q_2)=(2,2,2,2,1,999)$ and $Q_0=\pi(Q_1)=(1,1,1,1000)$\hspace{0.25mm}.\vspace{0.5mm}

\noindent  The extension of $Q_0$ \mbox{is $Q_0^*=(1,1,1,1001)$\hspace{0.15mm}.} By Example 20.2 we get that
$j_{\hspace{0.2mm}0}=4$ \mbox{and $\zeta_{\hspace{0.3mm}0}=1004$} are the cyclic parameters of $Q_0^*$\hspace{0.25mm}.\vspace{0.5mm}

\noindent By Example 17.3   we get that $D(Q_1) \ne \O$, and $\alpha_1^*=1003$ and $\gamma_1^*=1$ are the least progression parameters  of $D(Q_1)$ with respect to $\alpha_1=\delta(Q_1)+1$. The equation $x\alpha_1^*=y\zeta_0$ is equivalent to $1003x=1004y$. The least  positive integer solution is $x=1004$ and $y=1003$. By (21.3) we get that\vspace{0.5mm}

\hspace{20mm}$j_1=2x\gamma_1^*+yj_{0}=2\cdot1004\cdot 1 + 1003 \cdot 4=6020$,

\hspace{20mm}$\zeta_1=y\zeta_{0}+j_1=1003 \cdot 1004 + 6020=1013032$.\vspace{0.5mm}

\noindent By Example 17.2  we get that $D(Q_2) \ne \O$, and $\alpha_2^*=1008$ and $\gamma_2^*=1$ are the least progression parameters  of $D(Q_2)$ with respect to $\alpha_2=\delta(Q_2)+1$. The equation $x\alpha_2^*=y\zeta_1$ is equivalent to $1008x=1013032y$. By Example 23.3   the least  positive integer solution \mbox{is $x=126629$} and $y=126$. By (21.3) we get

\hspace{10mm}$j_2=2x\gamma_2^*+yj_{1}=2\cdot126629\cdot 1 + 126 \cdot 6020=1011778$,

\hspace{10mm}$\zeta_2=y\zeta_{1}+j_2=126 \cdot 1013032 + 1011778=128653810$\hspace{0.35mm}.\vspace{0.5mm}

\noindent Let $A_i^{\,\infty}=A^{\,\infty}_i(Q_i)$  for $0 \le i \le 2$. By (22.2)  we get that $j_0=4$,  $j_1=6020$ and $j_2=1011778$ are the least even vector periods \mbox{of
 $V(A_0^{\,\infty})$}, $V(A_1^{\,\infty})$  and $V(A_2^{\,\infty})$. By (22.3)  we get that $\zeta_0=1004$, $\zeta_1=1013032$ \mbox{and  $\zeta_2=128653810$} are the minimal periods of
 $A_0^{\,\infty}$, $A_1^{\,\infty}$  and $A_2^{\,\infty}$.\vspace{2mm}

\noindent {\textbf{Example 26.3.}} Suppose $A^{\,\infty}$ is generated from $A=111000011000$  by the  symmetric shift register with respect to $k=2$, $p=2$ and $n=12$. We note that $w(A)=k+p+1$ and $A$ starts with 1. By Observation 14.1 we get \mbox{that  $V(A)=(3,4,2,3) \in M_2^+$} and we  conclude that $A$ satisfies (19.1).\vspace{0.5mm}

 \noindent Let $Q_2=V(A)=(3,4,2,3)$. By (19.2) we  get  $A^{\,\infty}=A_2^{\,\infty}(Q_2)$.
 \noindent Let $j_0, j_1, j_2$ and $\zeta_0,\zeta_1, \zeta_2$ be the dynamical parameters of $Q_2$ with respect to $p=2$.  \mbox{By (22.3)} we get that   $\zeta_2$ is the minimal period of $A^{\,\infty}=A_2^{\,\infty}(Q_2)$. We will prove that  $\zeta_2=94$.\vspace{0.5mm}
 By Observation 13.2 and Example 13.7 we get that\vspace{0.05mm}

  \hspace{16.6mm}   $Q_1=\pi(Q_2)=(2,3,1,3)$ and $Q_0=\pi(Q_1)=(1,5)$\hspace{0.25mm}.\vspace{0.5mm}

  \noindent  The extension of $Q_0$ \mbox{is $Q_0^*=(1,6)$\hspace{0.15mm}.} By Example 20.3 we get that
  $j_{\hspace{0.2mm}0}=2$ \mbox{and $\zeta_{\hspace{0.3mm}0}=7$} are the cyclic parameters of $Q_0^*$\hspace{0.25mm}.\vspace{0.5mm}

 \noindent By Example 17.5  we get that $D(Q_1) \ne \O$ and $\alpha_1^*=6$ and $\gamma_1^*=1$ are the least progression parameters  of $D(Q_1)$ with respect to $\alpha_1=\delta(Q_1)+1$.\vspace{0.5mm}

 \noindent The equation $x\alpha_1^*=y\zeta_0$ is equivalent to $6x=7y$. The least  positive integer solution is $x=7$ and $y=6$.\vspace{12mm}

 \noindent By (21.3) we get that\vspace{0.5mm}

\hspace{27mm}$j_1=2x\gamma_1^*+yj_{0}=2\cdot 7 \cdot 1+6 \cdot 2=26$,

\hspace{27mm}$\zeta_1=y\zeta_{0}+j_1=6\cdot 7+26=68$.\vspace{0.5mm}

\noindent  By Example 15.7 we get that $D(Q_2)=\O$. Then (21.2) implies that\vspace{0.5mm}

 \hspace{20mm}$j_2=j_1=26$ and $\zeta_2=\zeta_1+j_2=68+26=94$.\vspace{2mm}

\noindent {\large{\textbf{27. Visualization.}}}\vspace{1mm}

\noindent Suppose  $A^{\,\infty}=a_1a_2 \cdots$ is generated from $A=a_1 \cdots a_n $ where $a_1=1$,
\mbox{by the symmetric shift register} with parameters  $k$, $p$ and $n$.
We also suppose \mbox{that $ w(A) = k+p+1$} and $V(A) \in M_p^+$. Let $A_i=a_{i+1} \cdots a_{i+n}$ for $i \ge 0$.\vspace{0.8mm}

\noindent The string $A^{\,\infty}$ has an internal  structure that can be visualized as waves for small values of $p$\hspace{0.35mm}. This  can be done by plotting the relative  weight of each substring of length $n$\hspace{0.25mm}. Let $w_i=w(A_i)-k=w(a_{i+1} \cdots a_{i+n})-k$ for $i \ge 0$\hspace{0.25mm}. That means, we normalize the weights by subtracting $k$\hspace{0.25mm}. By plotting $w_0\hspace{0.25mm}, w_1, \cdots $  we obtain a wave structure.
 However  by plotting the \mbox{integers $w_i^*=p+1-w_i$} for $j \ge 0$\hspace{0.35mm}, we see this  structure more clearly.
 We call $w_0\hspace{0.25mm},w_1\hspace{0.25mm}, \cdots$  the weight parameters and $w_0^*,w_1^*, \cdots$ the modified weight parameters of $A^{\,\infty}$ with respect \mbox{to $k$\hspace{0.25mm}, $p$ and $n$\hspace{0.2mm}.}
 Next, we will illustrate this by  Example 26.1. We let $A_i^{\infty}$, $p_{\hspace{0.3mm}i}$\hspace{0.2mm}, $k_{\hspace{0.2mm}i}$ and $n_{\hspace{0.2mm}i}$ be as in Example 26.1 for $0 \le i \le 2$.\vspace{0.4mm}

\noindent  In Figure 1 we have plotted the modified weight parameters of $A_2^{\,\infty}$ with respect to $k_{\hspace{0.2mm}2}$\hspace{0.2mm}, $p_{\hspace{0.3mm}2}$ and $n_{\hspace{0.2mm}2}$\hspace{0.2mm}.\vspace{0.7mm}

\begin{picture}(300,60)(10,20)
\put(20,30){\line(1,3){12}}
\put(32,66){\line(1,-3){12}}
\put(44,30){\line(1,0){4}}
\put(48,30){\line(1,3){8}}
\put(56,54){\line(1,-3){8}}
\put(64,30){\line(1,0){8}}
\put(72,30){\line(1,3){4}}
\put(76,42){\line(1,-3){4}}
\put(80,30){\line(1,0){8}}
\put(88,30){\line(1,3){12}}
\put(100,66){\line(1,-3){12}}
\put(112,30){\line(1,3){8}}
\put(120,54){\line(1,-3){8}}
\put(128,30){\line(1,0){4}}
\put(132,30){\line(1,3){4}}
\put(136,42){\line(1,-3){4}}
\put(140,30){\line(1,0){16}}
\put(156,30){\line(1,3){12}}
\put(168,66){\line(1,-3){8}}
\put(176,42){\line(1,3){8}}
\put(184,66){\line(1,-3){8}}
\put(192,42){\line(1,3){4}}
\put(196,54){\line(1,-3){8}}
\put(204,30){\line(1,0){20}}
\put(224,30){\line(1,3){8}}
\put(232,54){\line(1,-3){8}}
\put(240,30){\line(1,3){8}}
\put(248,54){\line(1,-3){4}}
\put(252,42){\line(1,3){8}}
\put(260,66){\line(1,-3){12}}
\put(272,30){\line(1,0){16}}
\put(288,30){\line(1,3){8}}
\put(296,54){\line(1,-3){8}}
\put(304,30){\line(1,3){4}}
\put(308,42){\line(1,-3){4}}
\put(312,30){\line(1,0){4}}
\put(316,30){\line(1,3){12}}
\put(328,66){\line(1,-3){12}}

\put(362,28){\scriptsize{$i$}}
\put(336,22){\scriptsize{$80$}}
\put(7,74){\scriptsize{$w_i^*$}}
\put(12,64){\tiny{$3$}}
\put(18,66){\line(1,0){4}}
\put(12,52){\tiny{$2$}}
\put(18,54){\line(1,0){4}}
\put(12,40){\tiny{$1$}}
\put(18,42){\line(1,0){4}}

\put(20,30){\qbezier[80](0,0),(165,0),(330,0)}

\put(20,30){\vector(0,5){40}}
\put(350,30){\vector(4,0){5}}

\end{picture}

\noindent \hspace{20mm}{\textbf{Figure 1.}} The modified weight parameters of $A_2^{\,\infty}$.\vspace{1.5mm}

\noindent In this figure we see three types of waves with different heights.
 These waves "collide" and "move" with different "velocities". But after the collisions they obtain there original form. These plots bear a resemblance to plots in soliton theory. We refer to [1]. We  do not use these wave
 structures   in the proofs. But there are analogous arithmetical structures that we reduce successively to obtain a simpler problem that we can solve. Afterwards we use an inverse process to find the periods we requested.

\noindent If we plot the modified wave parameters on each reduction step of this  process, we will observe that the waves of minimal heights  disappear on each step. Moreover, the height of the other waves are reduced by one unit.
 In this example  $A_2^{\,\infty}$ is reduced successively  to $A_1^{\,\infty}$ and $A_0^{\,\infty}$.
The next figure  contains the plot of  the modified weight parameters   of $A_1^{\,\infty}$   with respect \mbox{to $k_1$, $p_{\hspace{0.2mm}1}$} and $n_1$. In this figure we  have two types of waves.\vspace{0.7mm}

\begin{picture}(250,50)(10,20)
\put(20,30){\line(1,3){8}}
\put(28,54){\line(1,-3){8}}
\put(36,30){\line(1,0){4}}
\put(40,30){\line(1,3){4}}
\put(44,42){\line(1,-3){4}}
\put(48,30){\line(1,0){12}}
\put(60,30){\line(1,3){8}}
\put(68,54){\line(1,-3){8}}
\put(76,30){\line(1,3){4}}
\put(80,42){\line(1,-3){4}}
\put(84,30){\line(1,0){16}}
\put(100,30){\line(1,3){8}}
\put(108,54){\line(1,-3){4}}
\put(112,42){\line(1,3){4}}
\put(116,54){\line(1,-3){8}}
\put(124,30){\line(1,0){16}}
\put(140,30){\line(1,3){4}}
\put(144,42){\line(1,-3){4}}
\put(148,30){\line(1,3){8}}
\put(156,54){\line(1,-3){8}}
\put(164,30){\line(1,0){12}}
\put(176,30){\line(1,3){4}}
\put(180,42){\line(1,-3){4}}
\put(184,30){\line(1,0){4}}
\put(188,30){\line(1,3){8}}
\put(196,54){\line(1,-3){8}}
\put(204,30){\line(1,0){8}}
\put(212,30){\line(1,3){4}}
\put(216,42){\line(1,-3){4}}
\put(220,30){\line(1,0){8}}
\put(228,30){\line(1,3){8}}
\put(236,54){\line(1,-3){8}}
\put(244,30){\line(1,0){4}}
\put(248,30){\line(1,3){4}}
\put(252,42){\line(1,-3){4}}
\put(256,30){\line(1,0){12}}
\put(268,30){\line(1,3){8}}
\put(276,54){\line(1,-3){8}}
\put(284,30){\line(1,0){16}}
\put(300,30){\line(1,3){8}}
\put(308,54){\line(1,-3){4}}
\put(312,42){\line(1,3){4}}
\put(316,54){\line(1,-3){8}}
\put(324,30){\line(1,0){12}}

\put(362,28){\scriptsize{$i$}}
\put(333,22){\scriptsize{$80$}}
\put(7,61){\scriptsize{$w_i^*$}}

\put(12,52){\tiny{$2$}}
\put(18,54){\line(1,0){4}}
\put(12,40){\tiny{$1$}}
\put(18,42){\line(1,0){4}}

\put(20,30){\qbezier[80](0,0),(165,0),(330,0)}

\put(20,30){\vector(0,5){28}}
\put(350,30){\vector(4,0){5}}

\end{picture}

\noindent \hspace{21mm}{{\textbf{Figure 2.}}}   The modified weight parameters of $A_1^{\,\infty}$.\vspace{2mm}

\noindent The next figure  contains the plot of  the modified weight parameters   of  $A_0^{\,\infty}$
with respect to $k_0$, $p_{\hspace{0.3mm}0}$ and $n_0$\hspace{0.35mm}.
In this figure we  have only one type of waves.\vspace{0.1mm}

\begin{picture}(40,50)(10,20)
\put(20,30){\line(1,3){4}}
\put(24,42){\line(1,-3){4}}
\put(28,30){\line(1,0){24}}

\put(52,30){\line(1,3){4}}
\put(56,42){\line(1,-3){4}}
\put(60,30){\line(1,0){24}}

\put(84,30){\line(1,3){4}}
\put(88,42){\line(1,-3){4}}
\put(92,30){\line(1,0){24}}

\put(116,30){\line(1,3){4}}
\put(120,42){\line(1,-3){4}}
\put(124,30){\line(1,0){24}}

\put(148,30){\line(1,3){4}}
\put(152,42){\line(1,-3){4}}
\put(156,30){\line(1,0){24}}

\put(180,30){\line(1,3){4}}
\put(184,42){\line(1,-3){4}}
\put(188,30){\line(1,0){24}}

\put(212,30){\line(1,3){4}}
\put(216,42){\line(1,-3){4}}
\put(220,30){\line(1,0){24}}

\put(244,30){\line(1,3){4}}
\put(248,42){\line(1,-3){4}}
\put(252,30){\line(1,0){24}}

\put(276,30){\line(1,3){4}}
\put(280,42){\line(1,-3){4}}
\put(284,30){\line(1,0){24}}

\put(308,30){\line(1,3){4}}
\put(312,42){\line(1,-3){4}}
\put(316,30){\line(1,0){24}}

\put(362,28){\scriptsize{$i$}}
\put(333,22){\scriptsize{$80$}}
\put(7,53){\scriptsize{$w_i^*$}}

\put(12,40){\tiny{$1$}}
\put(18,42){\line(1,0){4}}

\put(20,30){\qbezier[80](0,0),(165,0),(330,0)}
\put(20,30){\vector(0,5){25}}
\put(350,30){\vector(4,0){5}}

\end{picture}

\noindent \hspace{21mm}{{\textbf{Figure 3.}}}   The modified weight parameters of $A_0^{\,\infty}$.\vspace{1.7mm}

 \hspace{45mm} {\large{\textbf{PART 3.}}  \vspace{ 0.5 mm}}

 \noindent In Section 28 we define and  derive basic properties of  shift symmetric vectors.
 Suppose $A^{\infty}=A_p^{\infty}(Q)$ where  $Q \in M_p$ and $p \ge 0$.  Then we prove in Part 7 that $V(A^{\,\infty})$ is the shift symmetric vector denoted $C_p^{\infty}(Q)$, we refer to (46.1). Hence, we can use shift symmetric vectors to determine the least even vector period of $V(A^{\,\infty})$.
 In Section 29 and 30 we  derive  periodic properties  \mbox{if  $p=0$.} Moreover, Section 31  contains some   preliminary reduction results.\vspace{2mm}

\noindent {\large{\textbf{28. Shift symmetric vectors.}}}\vspace{1mm}

\noindent  Suppose $p \ge 0$ and $Q=(q_1, \cdots, q_J,e_0) \in M$. Then $q_i \ge 1$ for $1 \le i \le J$, $e_0 \ge 0$ and $J$ is odd. Let $\lambda_0 = p+1$.  For $j \ge 0$ we let

\noindent(28.1)\hspace{12mm} $s_{j+1}=min\{q_{j+1},\lambda_j\}$ and $\lambda_{j+1}=\lambda_j-s_{j+1}$ if $j$ is even,

\noindent(28.2)\hspace{4mm} $s_{j+1}=min\{q_{j+1},p+1-\lambda_j\}$ and $\lambda_{j+1}=\lambda_j+s_{j+1}$ if $j$ is odd,

\noindent(28.3)\hspace{41mm} $e_{j+1}=q_{j+1}-s_{j+1}$,\vspace{0.2mm}

\noindent(28.4)\hspace{41mm} $q_{J+j+1}=e_{j}+ s_{j+1}$\,.\vspace{0.5mm}

\noindent Then   $Q^{\,\infty}=(q_1,q_2, \cdots )$\hspace{0.25mm} is called  the   shift  symmetric  vector  generated by $Q$ with respect to  $p$\hspace{0.2mm}.
The  vectors $(s_1, s_2, \cdots )$, $(e_0,e_1, \cdots )$ and $(\lambda_0, \lambda_1, \cdots )$   are called the associated sequences. In particular, we call $\lambda_0, \lambda_1, \lambda_2,\cdots $ \mbox{the $\lambda$}-parameters of $Q^{\,\infty}$. We let   $C^{\infty}_p(Q)=Q^{\,\infty}$.
That means, $C^{\infty}_p(Q)$ is  the   shift  symmetric  vector  generated by $Q$ with respect to  $p$\hspace{0.2mm}.\vspace{0.5mm}

\noindent Let $j \ge 0$. Then  (28.1), (28.2) and (28.3) imply that\vspace{0.5mm}

 \noindent (28.5) \hspace{34mm} $s_{j+1} \le q_{j+1}$ and $e_{j+1} \ge 0$,

 \noindent (28.6) \hspace{10mm}  $s_{j+1} \le \lambda_j$ if $j$ is even, and $s_{j+1} \le p+1-\lambda_j$ if $j$ is odd.\vspace{0.5mm}

 \noindent Since $e_0 \ge 0$, then (28.5) implies that\vspace{0.3mm}

 \noindent (28.7) \hspace{40.5mm} $e_j \ge 0$ for $j \ge 0$.\vspace{0.2mm}

 \noindent   As previously, $q^-=q-1$ if $q$ is an integer.  By (28.3) and (28.4) we get that\vspace{0.2mm}

\noindent(28.8)\hspace{14mm}  $q_{j+1}=s_{j+1}+e_{j+1}$ and $q_{j+1}^-=s_{j+1}^-+e_{j+1}$ for $j \ge 0$\hspace{0.25mm},\vspace{0.5mm}

\noindent(28.9)\hspace{31.5mm} $q_{J+j+1}^-=e_{j}+ s_{j+1}^-$ for $j \ge 0$\hspace{0.25mm}.\vspace{1.7mm}

 \noindent {\textbf{Observation 28.1.}} Suppose $q_{j+1}\ge 1$ and $1 \le \lambda_j \le p+1$ where $j \ge 0$ is even. Then $s_{j+1}>0$, $q_{J+j+1}>0$ and
 $0 \le \lambda_{j+1} \le p$.\vspace{0.5mm}

 \noindent {\textbf{Proof.}} By (28.1) we get that $1 \le s_{j+1} \le \lambda_j$ and $\lambda_{j+1}=\lambda_j-s_{j+1}$.
 Hence, $0 \le \lambda_{j+1} \le p$. Moreover, (28.4) and (28.7) imply that  $q_{J+j+1}=e_j+s_{j+1}>0$.\vspace{1.7mm}

\noindent {\textbf{Observation 28.2.}} Suppose $q_{j+1}\ge 1$ and $0 \le \lambda_j \le p$ where $j\ge 0$ is odd. \mbox{Then $s_{j+1}>0$, $q_{J+j+1}>0$ and
 $1 \le \lambda_{j+1} \le p+1$.}\vspace{0.5mm}

 \noindent {\textbf{Proof.}} By (28.2) we get that $1 \le s_{j+1} \le p+1-\lambda_j$ and $\lambda_{j+1}=\lambda_j+s_{j+1}$.
 Then $1 \le \lambda_{j+1} \le p+1$, and (28.4) and (28.7) \mbox{imply that  $q_{J+j+1}=e_j+s_{j+1}>0$.}\vspace{2mm}

 \noindent {\textbf{Observation 28.3.}} Let $j \ge 0$.\vspace{0.4mm}

\noindent a) $q_i >0$ for $1 \le i \le J+j$, and $s_{j}>0$ for  $1 \le i \le j$.

\noindent b) $1 \le \lambda_j \le p+1$ if $j$ is even, and $0 \le  \lambda_j \le p$
if $j$ is odd.\vspace{0.7mm}

\noindent {\textbf{Proof.}}  Since $q_i \ge 1$ for $1 \le i \le J$, and $0<\lambda_0 = p+1$, then the results are  true for $j=0$.
 Next, suppose the results are true for $j \ge 0$. Then  $q_{j+1}>0$ \mbox{since $1 \le j+1 \le J+j$.} If $j$ is even, then $1 \le \lambda_j \le p+1$  and \mbox{Observation 28.1} implies that the results are true for $j+1$. If $j$ is odd, then $0 \le \lambda_j \le p$ and Observation 28.2 implies that the results are true for $j+1$.\vspace{1.5mm}

 \noindent {\textbf{Observation 28.4.}} $q_{j+1}>0$ and  $s_{j+1}>0$ for $j \ge 0$.\vspace{0.4mm}

 \noindent {\textbf{Proof.}} This follows from Observation 28.3 a).\vspace{2mm}

\noindent {\large{\textbf{29. Auxiliary results.}}}\vspace{0.5mm}

\noindent Let  $\psi$  be as Section 20. We suppose in this section that $r$ is a vector period \mbox{of $Q^{\,\infty}=(q_1,q_2, \cdots )$\hspace{0.2mm}.} That means, $q_{r+i}=q_i$ for $i \ge 1$. It is easily seen that

\noindent (29.1)\hspace{29.5mm} $q_{\,yr+i}=q_i$  for $y \ge 0$ and $i \ge 1$\hspace{0.2mm}.\vspace{2mm}

\noindent {\textbf{Observation 29.1.}} Suppose  $j \ge 1$. Then\vspace{0.2mm}

 \noindent \hspace{20mm} $q_{j+i}=q_i$ for $1 \le i \le r$ $\Leftrightarrow$  $j$ is a vector period of $Q^{\,\infty}$.\vspace{0.8mm}

 \noindent {\textbf{Proof.}} Suppose $q_{j+i}=q_i$ for $1 \le i \le r$.
 Let $i \ge 1$.
 We choose $y$ and $x$ such that $i=yr+x$ where $y \ge 0$ and $1 \le x \le r$. Then  (29.1) implies that\vspace{0.5mm}

\hspace{21mm} $q_{j+i}=q_{j+yr+x}=q_{j+x}=q_x=q_{yr+x}=q_i$\hspace{0.2mm}.\vspace{0.8mm}

\noindent  Hence, $j$ is a vector period of $Q^{\,\infty}$. The reverse implication \mbox{is trivial.}\vspace{200mm}

 \noindent {\textbf{Observation 29.2.}}  Let  $Q^*=(q_1, \cdots, q_r)$.

\noindent a) If  $0 \le j \le r$, then $\psi^{\,j}(Q^*)=(q_{j+1}, \cdots, q_{j+r})$.

\noindent b) If  $1 \le j \le r$, then $\psi^{\,j}(Q^*)=Q^*$ $\Leftrightarrow$ $j$ is a vector period of $Q^{\,\infty}$.\vspace{1mm}

\noindent {\textbf{Proof.}} a) If $j =0$, this is trivial. Next, suppose $\psi^{\,j}(Q^*)=(q_{j+1}, \cdots, q_{j+r})$ where $0 \le j < r$. Since $q_{j+1}=q_{r+j+1}$, then a) is true for $j+1$ since\vspace{0.3mm}

$\psi^{\,j+1}(Q^*)=(q_{j+2}, \cdots, q_{j+r},q_{j+1})=(q_{(j+1)+1}, \cdots, q_{(j+1)+r-1},q_{(j+1)+r})$.\vspace{0.5mm}

\noindent b) By a) and Observation 29.1 we get that\vspace{0.6mm}

\hspace{17mm}$\psi^{\,j}(Q^*)=Q^*$ $\Leftrightarrow$ $(q_{j+1}, \cdots, q_{j+r})=(q_{1}, \cdots, q_{r})$

\hspace{11mm}$\Leftrightarrow$ $q_{j+i}=q_i$ for $1 \le i \le r$
$\Leftrightarrow$ $j$ is a vector period of $Q^{\,\infty}$.\vspace{2mm}

\noindent {\large{\textbf{30. The  case $p =0$.}}}\vspace{0.5mm}

\noindent In this section we suppose  $Q^{\,\infty}=C_0^{\infty}(Q)$ is the shift symmetric vector generated by  $Q=(q_1, \cdots, q_J,e_0) \in M$ with respect to $p=0$. Then $J$ is odd.  Moreover,
we suppose
$Q^{\,\infty}=(q_1,q_2, \cdots)$ and that $(s_1, s_2, \cdots )$, $(e_0,e_1, \cdots )$ and $(\lambda_0, \lambda_1, \cdots )$  are the associated sequences.
\mbox{Since $p+1=1$,} then  \mbox{Observation 28.3 b)} implies that\vspace{0.4mm}

\noindent (30.1) \hspace{14mm}$ \lambda_j =1$ if $j \ge 0$ is even, and $\lambda_j = 0$
if $j \ge 0$ is odd.\vspace{0.4mm}

\noindent By Observation 28.4 we get that $q_{j+1} \ge 1$. Hence, (28.1), (28.2) and (30.1) imply  \mbox{that $s_{j+1}=1$ for $j \ge 0$.}
Then  we get from  (28.3) and (28.4)  that\vspace{0.3mm}

\noindent (30.2)  \hspace{10mm} $q_{J+j+1}=e_j+s_{j+1}=e_j+1=e_j+s_j=q_j$   for $j \ge 1$,

\noindent (30.3)  \hspace{20mm} $r=J+1$ is an even vector period of $Q^{\,\infty}$.\vspace{0.2mm}

\noindent  In particular, $q_{J+1}=e_0+s_1=e_0+1$\hspace{0.1mm}. Let  $Q^*$ be the extension of $Q$\vspace{0.4mm}.
Then

\noindent (30.4) \hspace{8mm} $Q^*=(q_1, \cdots, q_J,e_0+1)=(q_1, \cdots, q_{J+1})=(q_1, \cdots, q_r)$\hspace{0.2mm}.\vspace{0.4mm}

\noindent Suppose  $1 \le j \le r$, then Observation 29.2 b), (30.3) and (30.4) imply that\vspace{0.4mm}

\noindent (30.5) \hspace{8.5mm} $\psi^{\hspace{0.2mm}j}(Q^*)=Q^*$ if and only if $j$ is a vector period of $Q^{\,\infty}$\vspace{0.5mm}

\noindent where $\psi$ is defined as in Section 20.\vspace{1.5mm}

 \noindent {\textbf{Proposition 30.1.}} Suppose  $j$ and $\zeta$ are the cyclic parameters of $Q^*$.

\noindent a) $j$ is the least even  integer such that
$1 \le j \le r$ and $\psi^{\hspace{0.2mm}j}(Q^*)=Q^*$.

\noindent b)   $sum(Q^{\,\infty},j)=sum(Q^*,j)=q_1+ \cdots+ q_j=\zeta$\hspace{0.2mm}.\vspace{0.3mm}

\noindent c)  $j$ is the least even vector period of $Q^{\,\infty}$.\vspace{0.5mm}

\noindent {\textbf{Proof.}} a) and b) follow from the definition of $j$ and $\zeta$  in Section 20.\vspace{0.5mm}

\noindent  c)  By a) and (30.5) we get that $j$ is an even vector period of $Q^{\,\infty}$. \mbox{Suppose  $i$} is an even vector period of $Q^{\,\infty}$ satisfying $1 < i < j$. Then  $1 \le i < r$, \mbox{and  (30.5)} implies that $\psi^{\hspace{0.2mm}i}(Q^*)=Q^*$.
According to a) this is a contradiction. Hence,
$j$ is the least even vector period of $Q^{\,\infty}$.\vspace{200mm}

\noindent {\large{\textbf{31. Reduction results.}}}\vspace{1mm}

\noindent Suppose $Q^{\,\infty}=C_p^{\infty}(Q)$ where $Q \in M^+_p$,
and $p >0$. According to   (14.4) we get \mbox{that $\pi(Q) \in M_{p-1}^+$.} Let  $Q_*^{\,\infty}=C_{p-1}^{\,\infty}(\pi(Q))$.\vspace{0.5mm}

\noindent If we know   the least even vector period $j^*$ of $Q_*^{\,\infty}$ and
 $sum(Q_*^{\,\infty},j^*)$, we can determine the least even vector period $j$ of $Q^{\,\infty}$ and
 $sum(Q^{\,\infty},j)$.  We divide the proof of this into two cases.
  The case $D(Q)\ne\O$
 is complicated and is treated in Part 5. In this section we suppose $D(Q)=\O$,\vspace{0.5mm}

\noindent (31.1) \hspace{20mm} $j^*$ is the least even vector period of $Q_*^{\,\infty}$ and

\hspace{25.3mm} $sum(Q_*^{\,\infty},j^*)=\zeta^*$ where $Q_*^{\,\infty}=C_{p-1}^{\,\infty}(\pi(Q))$.\vspace{0.6mm}

\noindent Let  $Q^{\,\infty}=(q_1,q_2, \cdots)$ and $Q=(q_1, \cdots, q_J,e_0)$. We will prove that\vspace{0.5mm}

\noindent (31.2) \hspace{26.5mm}$Q_*^{\,\infty}=(q_1-1,q_2-1,q_3-1 \cdots )$.\vspace{0.5mm}

\noindent Since $Q^{\,\infty}=(q_1,q_2, \cdots)$, then (31.1) and (31.2) imply that\vspace{0.7mm}

\noindent (31.3)\,\hspace{18mm}$j^*$ is the least even vector period of $Q^{\,\infty}$ and

\hspace{11mm}$sum(Q^{\,\infty},j^*)=q_1+ \cdots+q_{j^*}=sum(Q_*^{\,\infty},j^*)+j^*=\zeta^*+j^*$.\vspace{2mm}

\noindent  \noindent {{\textbf{Proof of (31.2).}}}
Let  $(s_1, s_2, \cdots )$, $(e_0,e_1, \cdots )$ and $(\lambda_0, \lambda_1, \cdots )$   be  the associated sequences of $Q^{\,\infty}$. In particular, we get that $\lambda_0=p+1$. Let\vspace{0.3mm}

\noindent   \hspace{19.5mm} $q_{j+1}^*=q_{j+1}-1$, $s_{j+1}^*=s_{j+1}-1$  and $e_j^*=e_j$ for $j \ge 0$,\vspace{0.5mm}

\noindent    \hspace{19.5mm} $\lambda_j^*=\lambda_j-1$ if $j \ge 0$ is even,
  and $\lambda_j^*=\lambda_j$ if $j \ge 1$ is odd.\vspace{0.4mm}

\noindent Let $p^*=p-1$ and $Q^*=\pi(Q)$. Since $D(Q)=\O$,  then $q_i>1$ for $1 \le i \le J$, and $e_0\ge 0$. According to Observation 13.2  we get that

\noindent (31.4) \hspace{12mm}$Q^*=\pi(Q)=(q_1-1, \cdots, q_J-1,e_0)=(q_1^*, \cdots, q_J^*,e_0^*)$.\vspace{0.5mm}

\noindent  If $j \ge 0$ is even, then $\lambda^*_{j+1}=\lambda_{j+1}$ and $\lambda_j-1=\lambda_j^*$, and (28.1)  implies that\vspace{0.5mm}

\noindent (31.5)   \hspace{8mm}$s^*_{j+1}=s_{j+1}-1=min\{q_{j+1}-1,\lambda_j-1\}=min\{q_{j+1}^*,\lambda_j^*\}$,

\hspace{10.5mm}$\lambda^*_{j+1}=\lambda_{j+1}=\lambda_j-s_{j+1}=(\lambda_j-1)-(s_{j+1}-1)=
\lambda_j^*-s_{j+1}^*$.\vspace{0.5mm}

 \noindent   If $j \ge 0$ is odd, then $\lambda^*_{j+1}=\lambda_{j+1}-1$ and $\lambda_j=\lambda_j^*$, and  (28.2) implies that\vspace{0.5mm}

\noindent (31.6)\hspace{1mm} $s^*_{j+1}=s_{j+1}-1=min\{q_{j+1}-1,p-\lambda_j\}=min\{q_{j+1}^*,p^*+1-\lambda_j^*\}$,

 \hspace{20.5mm} $\lambda^*_{j+1}=\lambda_{j+1}-1=\lambda_j+s_{j+1}-1=\lambda_j^*+s_{j+1}^*$.\vspace{0.5mm}

\noindent   If $j \ge 0$, then (28.3) and (28.4) imply that

\noindent (31.7) \hspace{1mm} $e_{j+1}^*=e_{j+1}=q_{j+1}-s_{j+1}=(q_{j+1}-1)-(s_{j+1}-1)=q_{j+1}^*-s_{j+1}^*$,\vspace{0.6mm}

\noindent (31.8) \hspace{15mm}$q_{J+j+1}^*=q_{J+j+1}-1=e_{j}+s_{j+1}-1=e_{j}^*+s_{j+1}^*$.\vspace{0.7mm}

\noindent  We note that $\lambda_0^*=\lambda_{0}-1=p=p^*+1$. Hence, by (31.4), $\cdots$, (31.8) we get that
 $(q_1^*,q_2^*, \cdots )$ is   the   shift  symmetric  vector  generated by $Q^*$ with respect
 \mbox{to  $p^*$.} That means,\vspace{0.6mm}

\noindent \hspace{8.5mm} $Q_*^{\infty}=C_{p-1}^{\infty}(\pi(Q))=C_{p^*}^{\infty}(Q^*)=(q_1^*,q_2^*, \cdots )=(q_1-1,q_2-1, \cdots )$,\vspace{0.9mm}

\noindent and the proof of (31.2) is complete.\vspace{200mm}

\hspace{45mm} {\large{\textbf{PART 4.}}  \vspace{ 0.5 mm}}

\noindent  We will define and study the structure of complete vectors.  These vectors  will play a central role in the forthcoming proofs and results.\vspace{1mm}

\noindent  {\large{\textbf{32. Complete vectors.}}\vspace{ 1 mm}}

\noindent  In this part we suppose\vspace{0.4mm}

 \noindent (32.1) \hspace{12mm}$Q^{\,\infty}=(q_1,q_2, \cdots)$ where  $q_1 > 1$ and $q_i \ge 1$ for $i \ge 2$,\vspace{0.3mm}

 \noindent (32.2) \ $\#\{ i \ge 1: q_i> 1$ and  $i$ even\} = $\#\{ i \ge 1: q_i > 1$ and  $i$ odd\} = $\infty$\hspace{0.2mm},\vspace{0.6mm}

 \noindent (32.3) \hspace{34mm} $\#\{ i \ge 1: q_i=1 \}  = \infty$\hspace{0.2mm}.\vspace{0.5mm}

  \noindent   Then $Q^{\,\infty}$ is called a complete vector.
 For $r \ge 0$ we let   $next(r)=r+2t+1$ and $t_{max}(r)=t$ where
$t \ge 0$ is maximal such that \mbox{$q_{r+2i}=1$ for $1 \le i \le t$.}
 \mbox{According to (32.2)} the functions $next(r)$ and $t_{max}(r)$ are well-defined. Moreover, these functions depend
 on $Q^{\,\infty}$, but this  will always be clear from the context.
 If $r \ge 0$ and $t \ge 0$, we note \mbox{that $next(r)=r+2t_{max}(r)+1$ and}\vspace{0.4mm}

\noindent (32.4) \hspace{3.5mm}$next(r)=r+2t+1$ $\Leftrightarrow$  $q_{r+2i}=1$ for $1 \le i \le t$, and $q_{r+2t+2}>1$.\vspace{0.5mm}

\noindent The distance function  $\tau$  of $Q^{\,\infty}$ is given by\vspace{0.4mm}

\noindent (32.5) \hspace{21mm}$\tau(0)=0$ and $\tau(r)=\delta(q_1, \cdots, q_r)$ for $r \ge 1$.\vspace{0.5mm}

\noindent Then $\tau(r)=q_{1}^-+ \cdots +q_{r}^-$ for $r \ge 1$, and $\tau(r+1)=\tau(r)+q_{r+1}^-$  for $r \ge 0$,\vspace{0.3mm}

\noindent (32.6) \hspace{11mm} $q_1+ \cdots+ q_r=\delta(q_1, \cdots, q_r)+r=\tau(r)+r$  for $r > 0$,

\noindent (32.7) \hspace{23.4mm} $\tau(r+1)=\tau(r)$ if $q_{r+1}=1$ and $r \ge 0$,

\noindent (32.8) \hspace{1.6mm} $\tau(s)=\tau(r)+\delta(q_{r+1}, \cdots,q_{s})
=\tau(r)+q_{r+1}^-+ \cdots +q_{s}^-$ if $0 \le r < s $,

\noindent (32.9) \hspace{23mm} $ \tau(r) < \tau(r+1)$ if $r \ge 0$ and $q_{r+1}>1$.\vspace{0.5mm}

\noindent By (32.1) we get that $q_1 >1$ and  $q_i \ge 1$ for $i \ge 2$. Moreover, by (32.2)  we get that $q_i > 1$ for an infinite number of indexes. Hence,

\noindent (32.10) \hspace{8mm} $\tau(0)=0 < \tau(1) \le \tau(2) \le \cdots$ and $\tau(r) \rightarrow \infty$ if $r \rightarrow \infty$.\vspace{2mm}

\noindent  {\large{\textbf{33. The  contraction vector.}}\vspace{ 1 mm}}

\noindent  We choose integers $r_0, r_1, r_2, \cdots$   such that

 \noindent (33.1)  \hspace{12.5mm} $r_0=0 < r_1 < r_2 <\cdots $ and     $r_{j+1}=next(r_j)$ for $j \ge 0$.\vspace{0.3mm}

 \noindent Let $t_{j+1}=t_{max}(r_j)$ for $j \ge 0$. Then we get that

\noindent (33.2)  \hspace{31.5mm} $r_{j+1}=r_j+2t_{j+1}+1$   for $j \ge 0$,\vspace{0.4mm}

\noindent (33.3)  \hspace{9.5mm} $q_{r_j+2i}=1$ for $1 \le i \le t_{j+1}$, and $q_{r_j+2t_{j+1}+2}>1$ for $j \ge 0$,\vspace{0.3mm}

\noindent (33.4)  \hspace{15.5mm} $r_0, r_2, r_4, \cdots$ are even, and  $r_1, r_3, r_5, \cdots$ are odd.\vspace{0.1mm}

\noindent We call  $r_0, r_1, r_2, \cdots$   the $r$-indexes   and  $t_1, t_2, \cdots$ the $t$-indexes of $Q^{\,\infty}$.\vspace{0.5mm}

\noindent   By (32.1) we get that $q_{r_0+1}=q_1>1$. If $j \ge 0$, then (33.2) and (33.3) imply that
  $q_{r_{j+1}+1}=q_{r_j+2t_{j+1}+2}>1$. Hence,

\noindent (33.5)\hspace{39mm}$q_{r_j+1}>1$ for $j \ge 0$.\vspace{0.5mm}\vspace{0.5mm}

\noindent
 The contraction vector of $Q^{\,\infty}$ is defined as\vspace{0.4mm}

\noindent (33.6)\hspace{6mm} $\pi(Q^{\,\infty})=(q_1^*,q_2^*,  \cdots)$
where $q_{j+1}^*=\tau(r_{j+1})-\tau(r_j)$ for $j \ge 0$\hspace{0.2mm}.\vspace{2mm}

 \noindent  {\large{\textbf{34. The  component decomposition.}}\vspace{ 1 mm}}

\noindent By the component decomposition defined in this section we can give an alternative  characterization of the contraction vector.\vspace{0.5mm}

\noindent Suppose $G  \subset Q^{\,\infty}$ is an odd component  succeeded by a coordinate $>1$.
\mbox{Then $G$} is called a proper odd component  in $Q^{\,\infty}$.\vspace{0.5mm}

\noindent If $Q^{\,\infty}=(G_1,G_2, \cdots )$ where $G_1, G_2, \cdots$ are proper odd components in
 $Q^{\,\infty}$, then  $(G_1,G_2, \cdots )$ is called the component decomposition of $Q^{\,\infty}$.
By Observation 34.3 we get that $Q^{\,\infty}$  has a unique component decomposition. Moreover, if
$Q^{\,\infty}=(G_1,G_2, \cdots )$ where $(G_1,G_2, \cdots )$ is  the component decomposition of $Q^{\,\infty}$, then we will prove that\vspace{0.5mm}

\noindent (34.1)\hspace{34mm} $\pi(Q^{\,\infty})=(\delta(G_1),\delta(G_2), \cdots ).$\vspace{2mm}

\noindent {{\textbf{Observation 34.1.}}} If $G=(q_{r+1}, \cdots, q_{r+2t+1}) \subset Q^{\,\infty}$. Then $G$ is a proper odd component  $\Leftrightarrow$ $q_{r+2i}=1$ for $1 \le i \le t$, \mbox{and $q_{r+2t+2}>1$}
 $\Leftrightarrow$ $t=t_{max}(r)$.\vspace{0.4mm}

 \noindent {{\textbf{Proof.}}} This is trivial.\vspace{2mm}

 \noindent {{\textbf{Observation 34.2.}}} Suppose $G=(q_{r+1}, \cdots, q_{s}) \subset Q^{\,\infty}$. Then $G$ is a proper odd component if and only if $s=next(r)$.\vspace{0.5mm}

 \noindent {{\textbf{Proof.}}} This is proved as Observation 11.2 by using Observation 34.1.\vspace{2mm}

 \noindent {{\textbf{Observation 34.3.}}} Suppose $0 \le r \le J$. Then $G=(v_{r+1}, \cdots, v_{s})$ where $s=next(r)$, is the unique proper odd component starting with $v_{r+1}$.\vspace{0.5mm}

 \noindent {{\textbf{Proof.}}} Follows from Observation 34.2.\vspace{2mm}

 \noindent {{\textbf{Proof of (34.1).}}} Suppose $Q^{\,\infty}=(G_1,G_2, \cdots )$ where $(G_1,G_2, \cdots )$ is  the component decomposition of $Q^{\,\infty}$.\vspace{0.5mm}

 \noindent We choose indexes $r_0=0 < r_1 < r_2 <\cdots  $ such that $G_{j+1}=(q_{r_j+1}, \cdots, q_{r_{j+1}})$  for $j \ge 0$.
 Since $G_{j+1}$ is a proper odd  component in $Q^{\,\infty}$, then Observation 34.2 implies that $r_{j+1}=next(r_j)$ for $j \ge 0$. Hence, $r_0 ,r_1, r_2,\cdots $ are the $r$-indexes of $Q^{\,\infty}$. By (33.6) we get that
  $\pi(Q^{\,\infty})=(q_1^*,q_2^*,  \cdots)$
\mbox{where $q_{j+1}^*=\tau(r_{j+1})-\tau(r_j)$ for $j \ge 0$\hspace{0.2mm}.}\vspace{0.5mm}
 Finally, according to (32.8)  we get that
  $\delta(G_{j+1})=\tau(r_{j+1})-\tau(r_j)=q_{j+1}^*$  for $j \ge 0$.\vspace{2mm}

\noindent  {\large{\textbf{35. The  distance vector.}}\vspace{ 1 mm}}

 \noindent According to (32.3) we can choose $c_0,c_1,c_2, \cdots$ such that

 \noindent (35.1) \hspace{31mm} $c_0=0< c_1 < c_2 < c_3 < \cdots $
 \vspace{0.5mm}

 \noindent (35.2) \hspace{5mm} $c_{i+1}$ is the least index $>c_i+1$ such that $q_{c_{i+1}}=1$ for $i \ge 0$.\vspace{0.5mm}

\noindent We call $c_0, c_1, c_2, \cdots $  the $c$\,-\,indexes of $Q^{\,\infty}$. Let $D(Q^{\,\infty})=D^{\,\infty}$ where

\noindent (35.3)\hspace{20.4mm} $D^{\,\infty}=(d_1,d_2,  \cdots)$
and $d_i=\tau(c_i)$ for $i \ge 1$\hspace{0.2mm}.\vspace{0.5mm}

\noindent We call $D(Q^{\,\infty})$ the distance vector of $Q^{\,\infty}$.
 By (32.10) and (35.1) we get

\noindent (35.4) \hspace{20.4mm}  $0 <d_1 \le d_2 \le \cdots$ and  $d_i \rightarrow \infty$ if $\rightarrow \infty$\,.\vspace{0.5mm}

\noindent We let    $r(D^{\infty},\beta)=\#\{i \ge 1: d_i \le \beta\}$ for $\beta \ge 0$\,.
Suppose $\beta \ge 0$. Then we get from (35.4) that the following statements are true:\vspace{0.4mm}

\noindent (35.5) \hspace{17mm} If $y >0$, then $r(D^{\infty},\beta)=y$ $\Leftrightarrow$  $d_y \le \beta < d_{y+1}$.\vspace{0.3mm}

\noindent (35.6) \hspace{34.5mm} $r(D^{\infty},\beta)=0$ $\Leftrightarrow$  $\beta < d_{1}$.\vspace{0.5mm}

\noindent  We call $\beta>0$ a progression parameter of $D^{\,\infty}$ if\vspace{0.5mm}

\noindent (35.7)\hspace{18mm} $d_{y+i}=d_i+\beta$ for $i \ge 1$, where $y =r(D^{\infty},\beta)$.\vspace{0.3mm}

\noindent For each $r \ge 0$ we let $y(r)=y$ where $y$  is maximal such that $c_{y} \le r$. In particular, since $c_0=0< c_1$, then $y(0)=0$.  According to  (35.1) we get that the following result is true:\vspace{0.4mm}

\noindent (35.8) \hspace{19mm} If $y \ge 0$, then  $y(r)=y$  $\Leftrightarrow$ $c_y \le r < c_{y+1}$.\vspace{1.4mm}

\noindent  {\large{\textbf{36. Auxiliary results.}}\vspace{ 1 mm}}

\noindent {\textbf{Observation 36.1.}} Suppose $y(r)=y$ and $q_{r+1}>1$ where $r \ge 0$.

\noindent a) $c_y \le r < c_{y+1}$, $c_y < r+1 < c_{y+1}$ and $y(r+1)=y$.

\noindent b) Suppose $q_{r+2}=1$. Then $c_{y+1}=r+2$.\vspace{0.5mm}

\noindent {\textbf{Proof.}} a) By (35.8) we get that  $c_y \le r < c_{y+1}$.\vspace{0.5mm}
 Since $q_{c_{y+1}}=1$ and $q_{r+1}>1$, then $c_y < r+1 < c_{y+1}$. Hence, we get by (35.8) that $y(r+1)=y$.\vspace{1mm}

\noindent b) By a) we get that $c_y+1<r+2 \le c_{y+1}$.  Since $c_y+1<r+2$ and $q_{r+2}=1$, then
(35.2) implies that $c_{y+1} \le r+2$. Hence, $c_{y+1}=r+2$.\vspace{2mm}

\noindent {\textbf{Observation 36.2.}}  Suppose $r \ge 0$ and $t \ge 1$.

\noindent a)  If $c_y=r$ and $q_{r+2}=1$, then $c_{y+1}=r+2$.\vspace{0.3mm}

\noindent b) If  $c_{y+1}=r+2$ and $q_{r+2i}=1$ for $1 \le i \le t$, then $c_{y+i}=r+2i$
for $1 \le i \le t$.\vspace{0.6mm}

\noindent {\textbf{Proof.}} a) follows from (35.2), and b) follows from a) by induction.\vspace{2mm}

\noindent {\textbf{Observation 36.3.}} Suppose $\beta=\tau(r)$  and $q_{r+1}>1$ where $r \ge 0$.

\noindent a) If $i \ge 1$, then $d_i \le \beta$ $\Leftrightarrow$ $\tau(c_i) \le \tau(r)$
$\Leftrightarrow$ $c_i \le r$.

\noindent b) $r(D^{\infty},\beta)=y(r)$.\vspace{0.5mm}

\noindent {\textbf{Proof.}} a) follows from (32.9), (32.10)  and (35.3).\vspace{0.3mm}

\noindent b) If $y(r)=0$, then (35.8) and a) imply that
 $c_0 \le r < c_1$ and $d_1>\beta$. Hence, we get from (35.6)  that $r(D^{\infty},\beta)=0$.
If $y(r)=y>0$, then we get by (35.8) and a) that $c_y \le r < c_{y+1}$ and   $d_y \le \beta < d_{y+1}$. Hence, according to (35.5) we get that
$r(D^{\infty},\beta)=y$.\vspace{2mm}

 \noindent  {\large{\textbf{37. Properties of the indexes.}}\vspace{ 1 mm}}

\noindent Let $\beta_j=\tau(r_j)$ and $y_j=y(r_j)$ for $j \ge 0$. Then $\beta_0=\tau(r_0)=\tau(0)=0$
\mbox{and $y_0=y(r_0)=y(0)=0$.}
 By  (33.5)} and Observation  36.3 b) we get that\vspace{0.5mm}

\noindent (37.1)\hspace{15.8mm}$q_{r_j+1}>1$ and  $y_j=y(r_j)=r(D^{\,\infty},\beta_j)$ for $j \ge 0$.\vspace{2mm}

\noindent {\textbf{Observation 37.1.}} Let $j \ge 0$.

\noindent a)   $c_{y_j+i}=r_j+2i$  for $1 \le i \le t_{j+1}$.

\noindent b) $y_{j+1}=y_j+t_{j+1}$.\vspace{0.6mm}

\noindent {\textbf{Proof.}} By (33.3) and (33.5) we get that\vspace{0.3mm}

\noindent (37.2) \hspace{18mm} $q_{r_j+2i}=1$ for $1 \le i \le t_{j+1}$, and $q_{r_j+1}>1$.\vspace{0.5mm}

\noindent a) If $t_{j+1}=0$, this is trivial. Suppose $t_{j+1}>0$. By (37.2) we get
\mbox{that $q_{r_j+2}=1$} and $q_{r_j+1}>1$.
Since  $y(r_j)=y_j$, then Observation 36.1 b) implies that $c_{y_j+1}=r_j+2$.
\mbox{Hence,
(37.2) and Observation 36.2 b) imply that a) is true.}\vspace{0.8mm}

\noindent b)  First,  we suppose $t_{j+1}=0$.  By (33.2)  we get  that $r_j+1=r_{j+1}$.
\mbox{Since $q_{r_j+1}>1$} and $y(r_j)=y_j$, then Observation 36.1 a) implies that\vspace{0.3mm}

 \hspace{22.5mm}$y_{j+1}=y(r_{j+1})=y(r_j+1)=y_j=y_j+t_{j+1}$.\vspace{0.5mm}
 \vspace{0.4mm}

\noindent Next, we suppose $t_{j+1}>0$. By a), (33.2) and (35.2) we get that\vspace{0.4mm}

\noindent \hspace{2mm}$c_{y_j+t_{j+1}}=r_j+2t_{j+1}=r_{j+1}-1< r_{j+1}$ and $c_{y_j+t_{j+1}+1} \ge  c_{y_j+t_{j+1}}+2>r_{j+1}$.\vspace{0.4mm}

\noindent Hence, (35.8) implies that $y_{j+1}=y(r_{j+1})=y_j+t_{j+1}$.\vspace{2mm}

\noindent {\textbf{Observation 37.2.}} Suppose $1 \le i \le t_{j+1}$ where  $j \ge 0$. Then

\noindent \hspace{23mm} $\tau(r_j+2i-1)=\tau(r_j+2i)=\tau(c_{y_j+i})=d_{y_j+i}$.\vspace{0.4mm}

\noindent {\textbf{Proof.}}  By  (33.3) we get that  $q_{r_j+2i}=1$. Hence, the first equality is true. The last equalities follow from Observation 37.1 a) and (35.3).\vspace{2mm}

\noindent {\textbf{Observation 37.3.}}
 a) $y_0=0$ and $y_j=t_1+ \cdots+t_j$ for $j \ge 1$.

  \noindent b) $r_j=2y_j+j$ for $j \ge 0$.

  \noindent {\textbf{Proof.}}  a) We note that $y_0=0$. Hence, we get by using \mbox{Observation 37.1 b)}  that $y_1=y_0+t_1=t_1$.
  Suppose $y_j=t_1+ \cdots+t_j$ where $j >0$. Then  Observation 37.1 b) implies that $y_{j+1}=t_1+ \cdots+t_j+t_{j+1}$.\vspace{0.5mm}

  \noindent b) Since $r_0=y_0=0$, the result is true for $j=0$. If $r_j=2y_j+j$ where $j \ge 0$, then Observation 37.1 b) and (33.2) imply that\vspace{0.2mm}

  \hspace{6.7mm} $r_{j+1}=r_j+2t_{j+1}+1=2y_j+j+2t_{j+1}+1=2y_{j+1}+(j+1)$.\vspace{2mm}

 \noindent {\textbf{Observation 37.4.}}
 a) $y_0=0 \le y_1 \le y_2 \le \cdots$  and $y_j \rightarrow \infty$ if $j \rightarrow \infty$.\vspace{0.6mm}

  \noindent b) $sum(Q^{\,\infty},r_j)=q_1+ \cdots+q_{r_j}=\tau(r_j)+r_j=\beta_j+r_j$ for $j \ge 1$.\vspace{0.5mm}

  \noindent {\textbf{Proof.}}  a) Since $r_0=0 < r_1 < r_2 < \cdots$ and $y_j=y(r_j)$  for $j \ge 0$, the inequalities in a) are true.\vspace{12mm}

 \noindent  Let $i>0$\,. Choose $j$ such that $r_j>c_i$\,. Then we get that $y_j=y(r_j) \ge i$\,. Hence,  $y_j \rightarrow \infty$ if $j \rightarrow \infty$.\vspace{0.5mm}

\noindent  b) follows from (32.6).\vspace{2mm}

\noindent {\textbf{Observation 37.5.}} a) $\beta_0=0< \beta_1< \beta_2 < \beta_3 < \cdots$\,.

\noindent  b) $q_{j+1}^*=\tau(r_{j+1})-\tau(r_j)=\beta_{j+1}-\beta_j>0$  for $j \ge 0$.\vspace{0.3mm}

\noindent  c) $\beta_0=0$ and $\beta_j=q_1^*+ \cdots+q_j^*$ for $j \ge 1$.\vspace{0.5mm}

 \noindent {\textbf{Proof.}} a)   Since $r_j < r_{j+1}$, then (32.9), (32.10) and (33.5) imply \mbox{that $q_{r_j+1}>1$} and $\beta_{j}=\tau(r_{j})<\tau(r_{j}+1) \le \tau(r_{j+1})=\beta_{j+1}$ for $j \ge 0$.\vspace{0.5mm}

 \noindent b) follows from a) and (33.6).\vspace{0.3mm}

 \noindent c) By  b) we get  $\beta_0=0$,\,
 $\beta_1=\beta_0+q_1^*=q_1^*,\,\, \beta_2=\beta_1+q_2^*=q_1^*+q_2^*$, etc.\vspace{2mm}

 \noindent \textbf{{\large 38. Periodic properties.}}\vspace{1mm}

\noindent Let $\pi(Q^{\,\infty})=(q_1^*,q_2^*, \cdots)$ be as in  (33.6).  In this section we suppose $j^*$ is the least even vector period of $\pi(Q^{\,\infty})$ and
$\zeta^*=q_1^*+ \cdots+q_{j^*}^*=sum(\pi(Q^{\,\infty}),j^*)$.
According to Observation 37.5 c) we get that $\beta_{j^*}=\zeta^*$. It is easily seen that\vspace{0.5mm}

\noindent (38.1) \hspace{26.3mm} $q^*_{yj^*+i}=q^*_{i}$ for $i \ge 1$ and $y >0$,\vspace{0.3mm}

\noindent (38.2) \hspace{4.5mm}$j$ is an even vector period of $\pi(Q^{\,\infty})$ $\Leftrightarrow$ $j=yj^*$ where $y >0$.\vspace{2mm}

 \noindent {\textbf{Observation 38.1.}}  a) $\beta_{yj^*}=y\zeta^*$ for $y \ge 1$.\vspace{0.3mm}

 \noindent b) If $y \ge 1$, then $j=yj^*$ $\Leftrightarrow$ $\beta_j=\beta_{yj^*}$
 $\Leftrightarrow$ $\beta_{j}=y\zeta^*$.\vspace{0.5mm}

 \noindent c) $j$ is an even vector period of $\pi(Q^{\,\infty})$ $\Leftrightarrow$  $\beta_j=y\zeta^*$ where $y >0$.\vspace{0.6mm}

 \noindent {\textbf{Proof.}}
a)  is true for $y =1$. Suppose a) is true for $y \ge 1$. By (38.1)  we get that
$q_{yj^*+1}^*+ \cdots + q_{(y+1)j^*}^*=q_{yj^*+1}^*+ \cdots + q_{yj^*+j^*}^*=q_1^*+ \cdots +q_{j^*}^*=\zeta^*$, and Observation 37.5 c) implies that

 $\beta_{(y+1)j^*}=q_1^*+ \cdots+q_{(y+1)j^*}^*=\beta_{yj^*}+q_{yj^*+1}^*+ \cdots + q_{(y+1)j^*}^*=y\zeta^*+\zeta^*$.\vspace{0.6mm}

 \noindent b) follows from a) and Observation 37.5 a).\vspace{0.3mm}

  \noindent   c) follows from b) and (38.2).\vspace{2mm}

 \hspace{45mm} {\large{\textbf{PART 5.}}  \vspace{ 0.5 mm}}

 \noindent We will describe in Section 41 how minimal periods can be determined.
  The deductions are based on the reduction results in Section 31 and 40.
  Moreover, we will use a lot of results that will be proved in the next parts.\vspace{1.5mm}

 \noindent {\large{\textbf{39.  A crucial observation.}}}\vspace{1.5mm}

\noindent {{\textbf{Observation 39.1.}}} Let  $A^{\,\infty}=A^{\infty}_{p}(Q)$ and $Q^{\,\infty}=C^{\infty}_{p}(Q)$ where $Q \in M^+_p$ and $p  \ge 0$. Suppose $j$ is the least even vector period of $Q^{\,\infty}$. Then\vspace{0.5mm}

\hspace{17mm} $j$ is the least even vector period of $V(A^{\,\infty})$

\hspace{17mm} and $sum(Q^{\,\infty},j)$  is the minimal period of $A^{\,\infty}$.\vspace{0.5mm}

  \noindent {{\textbf{Proof.}}} By Observation 14.3 we get that $Q \in M_p$. Moreover, by (46.1) we get that $V(A^{\,\infty})=Q^{\,\infty}$. Hence, $j$ is the least even
vector period of $V(A^{\,\infty})$. Then  according to (18.2) we get that $sum(Q^{\,\infty},j)=sum(V(A^{\,\infty}),j)$ is the minimal period of $A^{\,\infty}$.\vspace{2mm}

\noindent {\large{\textbf{40.  Reduction results - complete case.}}}\vspace{0.5mm}

\noindent The following definition will be used to formulate the results.\vspace{0.5mm}

 \noindent Suppose $\alpha^*$, $\gamma^*$, $j^*$ and $\zeta^*$ are positive integers. Let $x$ and $y$ be the least positive integers satisfying $x\alpha^*=y\zeta^*$. Then we define\vspace{0.7mm}

 \hspace{7mm}$\omega(\alpha^*, \gamma^*, j^*, \zeta^*)=(r,\zeta)$ where $r=2x\gamma^*+yj^*$ and $\zeta=y\zeta^*+r$.\vspace{0.7mm}

\noindent We suppose in this section that\vspace{0.4mm}

\noindent (40.1) \hspace{22.8mm} $Q=(q_1, \cdots,q_J,e_0)\in M^+_p$   and $p> 0$.\vspace{0.4mm}

\noindent By   (14.4) we get that $\pi(Q) \in M_{p-1}^+$. Moreover, we suppose $D(Q) \ne \O$,\vspace{0.7mm}

\noindent (40.2) \hspace{18mm} $j^*$ is the least even vector period of $Q_*^{\,\infty}$ and

\hspace{23.5mm} $sum(Q_*^{\,\infty}),j^*)=\zeta^*$ where $Q_*^{\,\infty}=C_{p-1}^{\,\infty}(\pi(Q))$.\vspace{0.5mm}

\noindent Let
 $\alpha^*$ and $\gamma^*$ be the least progression parameters of $D(Q)$
with respect \mbox{to $\alpha=\delta(Q)+1$.}
Suppose $Q^{\,\infty}=C^{\infty}_{p}(Q)$. Let $(r,\zeta)=\omega(\alpha^*, \gamma^*, j^*, \zeta^*)$.
 We will in the end of this section prove  that\vspace{0.6mm}

\noindent (40.3) \hspace{3mm} $r$ is the least even vector period of $Q^{\,\infty}$
and $sum(Q^{\,\infty},r)=\zeta$.\vspace{0.6mm}

\noindent By Section 58 we get that $Q^{\,\infty}$ is complete. Let $r_0, r_1, \cdots$ be the $r$-indexes, and  $y_j=y(r_j)$ \mbox{and $\beta_j=\tau(r_j)$}  for $j \ge 0$. We also let
$\pi(Q^{\,\infty})$ be the contraction vector and $D^{\,\infty}=D(Q^{\,\infty})$ the distance vector of $Q^{\,\infty}$.\vspace{0.5mm}

\noindent By (71.13) we get  $\pi(Q^{\,\infty})= C_{p-1}^{\infty}(\pi(Q))=Q^{\,\infty}_*$. \mbox{Then (40.2) implies that}\vspace{0.7mm}

 \noindent (40.4)  $j^*$ is the least even vector period of $\pi(Q^{\,\infty})$
 and $sum(\pi(Q^{\,\infty}),j^*)=\zeta^*$.\vspace{0.7mm}

\noindent By   Proposition 67.1 we  get that\vspace{0.5mm}

\noindent  (40.5) \hspace{28.5mm} $r(D^{\,\infty},x\alpha^*)=x\gamma^*$ for $x > 0$,\vspace{0.3mm}

\noindent (40.6)  \hspace{4mm}$\beta$ is a progression parameter of $D^{\,\infty}$   $\Leftrightarrow$ $\beta=x\alpha^*$ where $x >0$.\vspace{0.5mm}

\noindent  If $s \ge 1$ is an even vector period of $Q^{\,\infty}$, then  Proposition 66.2 implies  that\vspace{0.5mm}

 \noindent (40.7)  \hspace{24mm}there exists $j > 0$ such that $s=r_j$\hspace{0.2mm}.\vspace{0.5mm}

 \noindent    If $j \ge 1$, then
  we get from  (48.2)  that\vspace{0.5mm}

 \noindent (40.8) \hspace{6.5mm}$r_j$ is a vector period of $Q^{\,\infty}$ $\Leftrightarrow$    $j$ is a vector period of $\pi(Q^{\,\infty})$

\hspace{10.8mm} and $\beta_j$ is a progression parameter of $D^{\,\infty}$.\vspace{0.5mm}

\noindent If $j \ge 0$, then   (33.4), (37.1), Observation 37.3 b) and 37.4 b) imply that\vspace{0.5mm}

\noindent (40.9)\hspace{26.7mm} $r_j$ is even if and only if $j$ is even,\vspace{0.4mm}

\noindent (40.10)\hspace{7.7mm}$y_j=r(D^{\,\infty},\beta_j)$, $r_j=2y_j+j$ and $sum(Q^{\,\infty},r_j)=\beta_j+r_j$.\vspace{20mm}

\noindent {{\textbf{Observation 40.1.}}} a) If   $y>0$\hspace{0.25mm}, then $j=yj^*$ $\Leftrightarrow$   $\beta_j=y\zeta^*$.

 \noindent   b) $j$ is an  even vector period of $\pi(Q^{\,\infty})$  $\Leftrightarrow$ $j=yj^*$ where $y>0$\hspace{0.25mm}.\vspace{0.2mm}

\noindent c) $j$ is an  even vector period of $\pi(Q^{\,\infty})$  $\Leftrightarrow$ $\beta_j=y\zeta^*$  where  $y>0$\hspace{0.25mm}.\vspace{0.5mm}

 \noindent {{\textbf{Proof.}}} The results follow from (40.4), (38.2),  Observation 38.1 b) and c).\vspace{1.5mm}

 \noindent {{\textbf{Observation 40.2.}}} Suppose $j \ge 1$. Then
 $r_j$ is an even vector period of $Q^{\,\infty}$ $\Leftrightarrow$
 $j$ is an even vector period \mbox{of $\pi(Q^{\,\infty})$}
and $\beta_j$ is a progression
parameter \mbox{of $D^{\,\infty}$} $\Leftrightarrow$
$\beta_j=x\alpha^*=y\zeta^*$ where  $x>0$ and $y>0$.\vspace{0.5mm}

\noindent {{\textbf{Proof.}}}
The first equivalence  follows from (40.8) and (40.9), and the last  from (40.6) and
Observation 40.1 c).\vspace{1mm}

\noindent {{\textbf{Observation 40.3.}}} Suppose $x>0$, $y >0$ and $j \ge 1$.

\noindent a) If $\beta_j=x\alpha^*$, then $y_j=r(D^{\,\infty},\beta_j)=r(D^{\,\infty},x\alpha^*)=x\gamma^*$.

\noindent b) If $\beta_j=x\alpha^*=y\zeta^*$, then $r_j=2x\gamma^*+yj^*$.\vspace{0.5mm}

\noindent {{\textbf{Proof.}}} a) follows from (40.5) and (40.10).\vspace{0.3mm}

\noindent b) Suppose $\beta_j=x\alpha^*$ and $\beta_j=y\zeta^*$. Then we get from a), \mbox{Observation 40.1 a)} and (40.10)  that
$y_j=x\gamma^*$,   $j=yj^*$ and $r_j=2y_j+j=2x\gamma^*+yj^*$.\vspace{2mm}

\noindent {{\textbf{Proof of (40.3).}}} Let
$(r,\zeta)=\omega(\alpha^*, \gamma^*, j^*, \zeta^*)$. Then
 $r=2x\gamma^*+yj^*$ and $\zeta=y\zeta^*+r$
 where $x$ and $y$ are
the least positive integers satisfying $x\alpha^*=y\zeta^*$.  Let $j = yj^*$. Since $j^*$ is even, then $j$ is even. By using \mbox{Observation 40.1 a)} we get that
 $\beta_j=y\zeta^*$\,. Hence, $\beta_j=x\alpha^*=y\zeta^*$.
  Then  Observation 40.2 implies that $r_j$ is an even vector period of $Q^{\,\infty}$. Moreover, $r=2x\gamma^*+yj^*=r_j$ where the last equality follows from Observation 40.3 b). Hence,  $r$ is an even vector period of $Q^{\,\infty}$.\vspace{0.6mm}

\noindent Suppose $s$ is an even vector period of $Q^{\,\infty}$.
 By (40.7)  there exists  $i >0$ such that $s=r_i$, and
 by Observation 40.2 there exist  integers $x^*>0$ \mbox{and $y^*>0$}  such that $\beta_i=x^*\alpha^*=y^*\zeta^*$. By the minimality property of $x$ and $y$  we obtain  $x^* \ge x$ and $y^*\ge y$.
 Hence, we get according to  Observation 40.3 b) that
$s=r_i=2x^*\gamma^*+y^*j^* \ge 2x\gamma^*+yj^*=r$. Then  $r$ is the least  even vector period of $Q^{\,\infty}$. Moreover, (40.10) implies that\vspace{0.5mm}

\hspace{15mm}$sum(Q^{\,\infty},r)=sum(Q^{\,\infty},r_j)=\beta_j+r_j=y\zeta^*+r=\zeta$.\vspace{2mm}

\noindent {\large{\textbf{41.  Determination of the periods.}}}\vspace{0.5mm}

\noindent Suppose   $Q \in M^+_p$ where $p \ge 0$.
  Let  $Q_p=Q$ and $Q_{i-1}=\pi(Q_i)$ \mbox{for $1 \le i \le p$\hspace{0.2mm}.}\vspace{0.5mm}

 \noindent  Observation 14.5 implies that $Q_{i} \in M^+_{i}$   \mbox{for $0 \le i \le p$\hspace{0.25mm}.}  Let
 $Q^{\,\infty}_i=C_{i}^{\infty}(Q_i)$    \mbox{for $0 \le i \le p$.} We note that $Q^{\,\infty}_{i-1}=C_{i-1}^{\infty}(\pi(Q_i))$
 for $1 \le i \le p$. We also let  $j_0, j_1, \cdots, j_p$ and $\zeta_0,\zeta_1, \cdots ,\zeta_p$ be the dynamical parameters \mbox{of $Q_p$} with respect to $p$.\vspace{100mm}

  \noindent {\textbf{Proposition  41.1.}} Suppose $j_{i-1}$ is the least even vector period of $Q_{i-1}^{\,\infty}$ \mbox{and $sum(Q_{i-1}^{\,\infty},j_{i-1})=\zeta_{i-1}$} where $1 \le i \le p$.
 Then $j_{i}$ is the least even vector period of $Q_{i}^{\,\infty}$ and $sum(Q_{i}^{\,\infty},j_{i})=\zeta_{i}$.\vspace{0.5mm}

 \noindent {\textbf{Proof.}}  Let $p_i=i$, $j^*=j_{i-1}$, $\zeta^*=\zeta_{i-1}$ and
 $Q^{\,\infty}_{*}=Q^{\,\infty}_{i-1}$. Next, we note \mbox{that
 $Q^{\,\infty}_{*}=Q^{\,\infty}_{i-1}=C_{i-1}^{\,\infty}(Q_{i-1})
 =C_{p_i-1}^{\,\infty}(\pi(Q_i))$.} Then we  get that\vspace{0.5mm}

  \noindent (41.1) \hspace{18mm}$j^*$ is the least even vector period of $Q^{\,\infty}_{*}$ and

 \hspace{22mm} $sum(Q^{\,\infty}_{*},j^*)=\zeta^*$ where  $Q^{\,\infty}_{*}=C_{p_i-1}^{\,\infty}(\pi(Q_i))$.\vspace{0.9mm}

 \noindent Suppose $D(Q_i)=\O$.
Since $Q^{\,\infty}_i=C_{p_i}^{\infty}(Q_i)$, then  (31.3) and (41.1) imply that
$j^*$ is the least even vector period of $Q_{i}^{\,\infty}$ and $sum(Q_{i}^{\,\infty},j^*)=\zeta^*+j^*$. By (21.2) we get  that $j_i=j_{i-1}=j^*$
and $\zeta_i=\zeta_{i-1}+j_{i}=\zeta^*+j^*$. Hence, $j_{i}$ is the least even vector period of $Q_{i}^{\,\infty}$ and $sum(Q_{i}^{\,\infty},j_{i})=\zeta_{i}$.\vspace{0.8mm}

\noindent  Next, we suppose $D(Q_i)\ne\O$. Let $\alpha_i^*$ and $\gamma_i^*$ be the least progression parameters of $D(Q_i)$
with respect to $\alpha_i=\delta(Q_i)+1$.
 By (21.3) we get \mbox{that $j_i=2x\gamma_i^*+yj_{i-1}$} and  $\zeta_i=y\zeta_{i-1}+j_i$ where
$x$ \mbox{and $y$} are the least positive integers satisfying
$x\alpha^*_i=y\zeta_{i-1}$. Hence,\vspace{0.5mm}

\hspace{22mm}$(j_i,\zeta_i)=\omega(\alpha_i^*, \gamma_i^*, j_{i-1}, \zeta_{i-1})=\omega(\alpha_i^*, \gamma_i^*, j^*, \zeta^*)$.\vspace{0.6mm}

\noindent Since $Q^{\,\infty}_i=C_{p_i}^{\infty}(Q_i)$, then we get from (40.3) and (41.1) that
$j_{i}$ is the least even vector period \mbox{of $Q_{i}^{\,\infty}$} and $sum(Q_{i}^{\,\infty},j_{i})=\zeta_{i}$.\vspace{2mm}

\noindent {\textbf{Proposition  41.2.}} If $0 \le i \le p$, then $j_{i}$ is the least even vector period \mbox{of $Q_{i}^{\,\infty}$} and $sum(Q_{i}^{\,\infty},j_{i})=\zeta_{i}$ for $0 \le i \le p$.\vspace{0.5mm}

 \noindent {\textbf{Proof.}} By (21.1) we get that  $j_0$ and $\xi_{\hspace{0.2mm}0}$ are the cyclic parameters of the extension $Q_0^*$ of $Q_0$.  Observation 14.2 implies that $Q_0 \in M^+_0=M$.\vspace{0.5mm}

  \noindent Since  $Q_0^{\,\infty}=C_{0}^{\infty}(Q_0)$, then Proposition 30.1 b) and c) imply that $j_0$ is the least even vector period \mbox{of $Q_0^{\,\infty}$} and
 $sum(Q_0^{\,\infty},j_0)=\zeta_0$.  Next, we suppose\vspace{0.5mm}

\noindent\hspace{4mm} $j_{i-1}$ is the least even vector period of $Q_{i-1}^{\,\infty}$ and $sum(Q_{i-1}^{\,\infty},j_{i-1})=\zeta_{i-1}$\vspace{0.7mm}

 \noindent where $1 \le i \le p$.
 Then Proposition 41.1 implies that the result is true for $i$.\vspace{2mm}

 \noindent {\textbf{Proposition  41.3.}} Let
$A_i^{\infty}=A^{\infty}_{i}(Q_i)$  where $0 \le i \le p$. Then
 $j_{i}$ is the least even vector period of $V(A_i^{\infty})$ and
 $\zeta_{i}$ is the minimal period of $A_i^{\infty}$.\vspace{0.5mm}

 \noindent {\textbf{Proof.}} By Proposition 41.2 we get that
 $j_{i}$ is the least even vector period \mbox{of $Q_{i}^{\,\infty}$} and $sum(Q_{i}^{\,\infty},j_{i})=\zeta_{i}$.
 Since  $Q^{\,\infty}_i=C_{i}^{\infty}(Q_i)$, then the results follow from
 Observation 39.1.\vspace{2mm}

\hspace{49mm} {\large{\textbf{PART 6.}}  \vspace{ 0.5 mm}}

\noindent We will derive properties of the infinite string $A^{\,\infty}=a_1a_2 \cdots  $  generated from
$A=a_1 \cdots a_n$ by the symmetric shift register with parameters $k$\hspace{0.06mm}, $p$ and $n$.
We suppose $k \le w(A) \le k+p+1$.\vspace{200mm}

\noindent {\large{\textbf{42. Basic properties.}}}\vspace{0.4mm}

\noindent We let $A_r=a_{r+1} \cdots a_{r+n}$ for $r \ge 0$\hspace{0.2mm}. In particular, $A_0=A$\hspace{0.2mm}. If $j \ge 0$\,, then

 \noindent (42.1)  \hspace{14mm}$a_{\,n+j+1} =a^{\,\prime} _{j+1} $ if $k\leq w(a_{\,j+2}
\;\cdots \;a_{\,j+n}) \leq k+p$\hspace{0.2mm},

\hspace{20mm}and $a_{\,n+j+1} =a_{\,j+1} $ otherwise.
 \vspace{0.4mm}

 \noindent  We also let
$w_r=w(a_{r+1} \cdots a_{r+n})-k=w(A_r)-k$ for $r \ge 0$\hspace{0.35mm}. By (3.1) we get that $w_r=a_{r+1}+ \cdots +a_{r+n}-k$ for $r \ge 0$\hspace{0.35mm}.
 If $a_{r+1}=1$ where $r \ge 0$, then

\noindent (42.2)\hspace{12mm} $a_{r+n+1}=a_{r+1}^{\,\prime}$ $\Leftrightarrow$
 $k \le w(a_{r+2} \cdots a_{r+n}) \le k+p$

\hspace{16mm} $\Leftrightarrow$
 $k < w(a_{r+1} \cdots a_{r+n}) \le k+p+1$ $\Leftrightarrow$ $0 < w_r \le p+1$.\vspace{0.5mm}

 \noindent If $a_{r+1}=0$ where $r \ge 0$, then

\noindent (42.3)\hspace{12mm} $a_{r+n+1}=a_{r+1}^{\,\prime}$ $\Leftrightarrow$
 $k \le w(a_{r+2} \cdots a_{r+n}) \le k+p$

\hspace{16mm} $\Leftrightarrow$
 $k \le  w(a_{r+1} \cdots a_{r+n}) < k+p+1$ $\Leftrightarrow$ $0 \le w_r < p+1$.\vspace{0.5mm}

\noindent  Since $w_r=a_{r+1}+ \cdots + a_{r+n}-k$ and  $w_{r+1}=a_{r+2}+ \cdots + a_{r+n+1}-k$\hspace{0.2mm},
then\vspace{0.2mm}

 \noindent (42.4)\hspace{20mm} $w_{r+1}=w_r+a_{r+n+1}-a_{r+1}$  for $r \ge 0$\hspace{0.2mm}.\vspace{0.6mm}

 \noindent Suppose $r \ge 0$\hspace{0.2mm}. Then we will prove that
 \vspace{0.5mm}

\noindent (42.5)  \hspace{3mm}$w_{\hspace{0.25mm}r+1}=w_r-1$  and $a_{\hspace{0.35mm}r+\hspace{0.25mm}n+\hspace{0.25mm}1}=0$
if $0 < w_r \le p+1$ and $a_{\hspace{0.35mm}r+1}=1$\hspace{0.1mm},\vspace{0.5mm}

\noindent (42.6)  \hspace{15mm}$w_{\hspace{0.25mm}r+1}=w_r$  and $a_{\hspace{0.35mm}r+\hspace{0.25mm}n+\hspace{0.25mm}1}=1$ if $ w_r=0$ and $a_{\hspace{0.35mm}r+1}=1$\hspace{0.1mm},\vspace{0.5mm}

\noindent (42.7)  \hspace{3mm}$w_{\hspace{0.25mm}r+1}=w_r+1$  and $a_{\hspace{0.35mm}r+\hspace{0.25mm}n+\hspace{0.25mm}1}=1$
if $0 \le w_r < p+1$ and $a_{\hspace{0.35mm}r+1}=0$\hspace{0.25mm},\vspace{0.2mm}

\noindent (42.8)  \hspace{12mm}$w_{\hspace{0.25mm}r+1}=w_r$  and $a_{\hspace{0.35mm}r+\hspace{0.25mm}n+\hspace{0.25mm}1}=0$ if $ w_r=p+1$ and $a_{\hspace{0.35mm}r+1}=0$\hspace{0.25mm}.\vspace{0.5mm}

\noindent In fact,   (42.5) and (42.6) follow from (42.2) and (42.4).
Moreover, (42.7) \mbox{and (42.8)} follow from (42.3) and (42.4).
By these relations we get that\vspace{0.4mm}

 \noindent (42.9) \hspace{11.4mm} $0 \le w_{r+1}<p+1$ if $0 \le w_r \le p+1$ and $a_{r+1}=1$\hspace{0.2mm},\vspace{0.4mm}

\noindent (42.10)  \hspace{9mm} $0 < w_{r+1} \le p+1$ if $0 \le w_r \le p+1$ and $a_{r+1}=0$\hspace{0.2mm}.\vspace{2mm}

\noindent {{\textbf{Observation 42.1.}}}
 a) $0 \le w_r \le p+1$  if $r \ge 0$\hspace{0.25mm}.\vspace{0.2mm}

 \noindent b) $0 \le w_{r} < p+1$  if $a_r=1$ and  $r \ge 1$\hspace{0.2mm}.\vspace{0.2mm}

 \noindent c) $0 < w_{r} \le  p+1$  if $a_r=0$ and  $r \ge 1$\hspace{0.2mm}.\vspace{0.7mm}

 \noindent {{\textbf{Proof.}}} $0 \le w_0 \le p+1$ since  $k \le w(a_1 \cdots a_n) \le k+p+1$\hspace{0.2mm}. Next, we suppose that $0 \le w_r \le p+1$  where $r \ge 0$.
  By   (42.9) and (42.10) we get that the results are  true for $r+1$\hspace{0.1mm}.\vspace{2mm}

\noindent {{\textbf{Proposition 42.2.}}} Suppose  $a_{r+i}=1$ for $1 \le i \le q$, where $r \ge 0$ and $q \ge 1$\hspace{0.2mm}.
 \mbox{Let  $s=min\{q,w_r\}$\hspace{0.2mm}.}\vspace{0.3mm}

 \noindent a) $w_{r+i}=w_r-i$ for $0 \le i \le s$, and $a_{r+n+i}=0$ for $1 \le i \le s$.

 \noindent b) If $s < q$, then $w_{r+i}=0$  for $s \le i \le q$, and $a_{r+n+i}=1$ for $s < i \le q$.

 \noindent c) $w_{r+q}=w_r-s$ and $a_{r+n+1} \cdots a_{r+n+q}=0_s1_{q-s}$.\vspace{0.5mm}

  \noindent {{\textbf{Proof.}}} a) The result is true for  $i=0$\hspace{0.2mm}.
 Next, we suppose $w_{r+i}=w_r-i$ where $0 \le i < s$. Since $i <s \le w_r$, then $0 < w_r-i =w_{r+i} \le p+1$\hspace{0.2mm}
 where the last inequality follows from Observation 42.1 a). Since $a_{r+i+1}=1$\hspace{0.2mm},
 then we get from (42.5)  that $a_{r+n+i+1}=0$ \mbox{and $w_{r+i+1}=w_{r+i}-1=w_r-(i+1)$\hspace{0.15mm}.}\vspace{0.5mm}

 \noindent b) Suppose $s <q$. Then $s=min\{q,w_r\}=w_r$\hspace{0.2mm}, and according to  a) we get that $w_{r+s}=w_r-s=0$\hspace{0.1mm}.
 Suppose $w_{r+i}=0$  where $s \le i < q$\hspace{0.25mm}.
 Since $a_{r+i+1}=1$, then   (42.6) implies that $a_{r+n+i+1}=1$ and $w_{r+i+1}=0$\hspace{0.2mm}.\vspace{0.5mm}

 \noindent c) If $s=q$\hspace{0.25mm}, then a) implies  $w_{r+q}=w_{r+s}=w_r-s$. If $s <q$, then $s=w_r$ \mbox{and  b)} implies
 $w_{r+q}=0=w_r-s$\hspace{0.2mm}. The last equality follows from a) and b).\vspace{1.7mm}

 \noindent {{\textbf{Proposition 42.3.}}} Suppose  $a_{r+i}=0$ for $1 \le i \le q$, where $r \ge 0$ and $q \ge 1$\hspace{0.2mm}.
 \mbox{Let  $s=min\{q,p+1-w_r\}$\hspace{0.2mm}.}\vspace{0.3mm}

 \noindent a) $w_{r+i}=w_r+i$ for $0 \le i \le s$, and $a_{r+n+i}=1$ for $1 \le i \le s$.

 \noindent b) If $s < q$, then $w_{r+i}=p+1$  for $s \le i \le q$, and $a_{r+n+i}=0$ for $s < i \le q$.

 \noindent c) $w_{r+q}=w_r+s$ and $a_{r+n+1} \cdots a_{r+n+q}=1_s0_{q-s}$.\vspace{0.5mm}

  \noindent {{\textbf{Proof.}}} a) The result is true for  $i=0$\hspace{0.2mm}.
 Suppose $w_{r+i}=w_r+i$ where $0 \le i < s$. Since $i <s \le p+1-w_r$, then $0 \le  w_{r+i}=w_r+i < p+1$\hspace{0.2mm}
 where the first inequality follows from Observation 42.1 a). Since $a_{r+i+1}=0$\hspace{0.2mm},
 then we get from  (42.7)  that

  \hspace{17.2mm}$a_{r+n+i+1}=1$ \mbox{and $w_{r+i+1}=w_{r+i}+1=w_r+(i+1)$\hspace{0.15mm}.}\vspace{0.6mm}

 \noindent b) Suppose $s <q$. Then $s=min\{q,p+1-w_r\}=p+1-w_r$\hspace{0.2mm}, and   a) implies that $w_{r+s}=w_r+s=p+1$\hspace{0.1mm}.
 Suppose $w_{r+i}=p+1$  \mbox{where $s \le i < q$.}
 Since \mbox{$a_{r+i+1}=0$,} then  (42.8) implies that $a_{r+n+i+1}=0$ \mbox{and $w_{r+i+1}=p+1$\hspace{0.2mm}.}\vspace{0.7mm}

 \noindent c) If $s=q$\hspace{0.25mm}, then a) implies  $w_{r+q}=w_{r+s}=w_r+s$. If $s <q$, then we get that $s=p+1-w_r$, and  b) implies
 $w_{r+q}=p+1=w_r+s$\hspace{0.2mm}. The last equality \mbox{follows from a) and b).}\vspace{2mm}

\noindent {\large{\textbf{43. Positive start strings.}}}\vspace{1mm}

 \noindent  Suppose $A=a_1 \cdots a_n$ where $a_1=1$.  Let $P$ be a start string of $A$ satisfying\vspace{0.4mm}

 $\overline{w}(P)=p+1$ \mbox{and $0 < \overline{w}(S) \le p+1$} for each  start string $S \ne \O$ of  $P$.\vspace{0.4mm}

 \noindent Then  we call $P$  a positive start string of $A$ of order $p+1$. We also suppose the \mbox{last bit of $P$ is 1.}
 We will prove that $V(A) \in M_p^+$. We note that $V(A) \in M$.\vspace{0.5mm}

 \noindent Suppose $V(A)=(v_1, \cdots, v_J,v_{J+1})$. Then $A=1_{v_1}0_{v_2} \cdots 1_{v_J}0_{v_{J+1}}$.
 We let

 \hspace{7mm} $P_1=1_{v_1}$, $P_2=1_{v_1}0_{v_2}, P_3=1_{v_1}0_{v_2}1_{v_3},  \cdots, P_J=1_{v_1}0_{v_2} \cdots 1_{v_J}$.\vspace{0.5mm}

 \noindent Let $\rho_0, \rho_1, \cdots, \rho_{J+1}$ be the alternating parameters of $V$. Then $\rho_0=0$,\vspace{0.1mm}

\noindent\hspace{1mm} $ \rho_1=\rho_0+v_1=v_1, \rho_2=\rho_1-v_2=v_1-v_2, \rho_3=\rho_2+v_3=v_1-v_2+v_3,\cdots$.\vspace{0.3mm}

\noindent It is easily seen that
$\overline{w}(P_i)=\rho_i$ for $1 \le i \le J$. Since  $P$ ends with $1$, then $P=1_{p+1}$ where $1 \le p+1 \le v_1$ or there exists an even $j$ such that
$P=1_{v_1}0_{v_2} \cdots 0_{v_j}1_s=P_j1_s$ where $2 \le j < J$  \mbox{and $1 \le s \le v_{j+1}$.}
 In the first case $v_1 \ge p+1$. Hence, Observation 14.1 implies that $V(A) \in M_p^+$.\vspace{0.4mm}

 \noindent In the second case it is sufficient to prove that

\hspace{30mm}$\rho_i>0$ for $1 \le i \le j$, and $\rho_{j+1} \ge p+1$.\vspace{0.3mm}

\noindent If $1 \le i \le j$,  then  $\rho_i=\overline{w}(P_i) >0$ since $P_i \subset P_j \subset P$.
Since $j$ is even, then\vspace{0.5mm}

\hspace{12.7mm}$\rho_{j+1}=\rho_j+v_{j+1} \ge \overline{w}(P_j)+s=\overline{w}(P_j1_s)=\overline{w}(P)=p+1$.\vspace{2mm}

\noindent {\large{\textbf{44. Auxiliary results.}}}\vspace{1mm}

\noindent {\textbf{Observation 44.1.}} Suppose $r \ge 0$.\vspace{0.3mm}

\noindent a) Suppose
$0 < w_{\hspace{0.25mm}r} < p+1$. Then $a_{r+n+1}=a_{r+1}^{\prime}$ and
 $w_{r+1} =w_{r}-\overline{w}(a_{r+1})$. \vspace{0.4mm}

\noindent b)  Suppose $w_{r+1}< w_{r}$. Then  $a_{r+1}=1$ and $w_{r+1}=w_r-1=w_r-\overline{w}(a_{r+1})$.\vspace{0.8mm}

\noindent {\textbf{Proof.}} a) If $a_{r+1}=1$, then (42.5) implies that $w_{r+1}=w_r-1=w_r-\overline{w}(a_{r+1})$. If $a_{r+1}=0$, then (42.7) implies that $w_{r+1}=w_r+1=w_r-\overline{w}(a_{r+1})$.\vspace{0.3mm}

\noindent b)  By Observation 42.1 a) we get that  $0 \le w_r \le p+1$.
If $a_{r+1}=0$, \mbox{then (42.7)} or (42.8) implies that $w_{r+1}\ge w_r$.
Hence, $a_{r+1}=1$. Then we get from  Observation 42.1 a) and (42.5)  that\vspace{0.4mm}

\hspace{10mm}$0 \le w_{r+1}<w_r \le p+1$
and $w_{r+1} =w_{r}-1=w_{r}-\overline{w}(a_{r+1})$.\vspace{2mm}

\noindent {\textbf{Proposition 44.2.}} Suppose $w_r=p+1$, $w_{r+x}=0$ and $0 < w_{r+i} <p+1$ \mbox{for $1 \le i < x$,} where  $0 < x \le n$ and $r \ge 0$. Let $A_r=a_{r+1} \cdots a_{r+n}$.\vspace{0.5mm}

\noindent  a) $w_{r+i+1}=w_{r+i}-\overline{w}(a_{r+i+1})$
for $0 \le i < x$.

\noindent b)  $w_{r+i}=p+1-\overline{w}(a_{r+1} \cdots {a}_{r+i})$  for $1 \le i \le x$.

\noindent c) $\overline{w}(a_{r+1} \cdots a_{r+x})=p+1$ and
 $0<\overline{w}(a_{r+1} \cdots a_{r+i})< p+1$ for $1 \le i < x$.\vspace{0.8mm}

 \noindent d) $w(A_r)=k+p+1$, $A_r$ starts with 1 and $V(A_r) \in M_p^+$.\vspace{0.5mm}

\noindent {\textbf{Proof.}} a) Since $w_{r+1}<w_r$, then we get according to Observation 44.1 b)  that
 $w_{r+1}=w_r-\overline{w}(a_{r+1})$. Suppose $1 \le i < x$. Then $0 < w_{r+i}<p+1$ and the result follows from Observation 44.1 a).\vspace{0.5mm}

\noindent b) By a) we get that $w_{r+1}=w_r-\overline{w}(a_{r+1})=p+1-\overline{w}(a_{r+1})$. Next,   we suppose \mbox{that $w_{r+i}=p+1-\overline{w}(a_{r+1} \cdots {a}_{r+i})$  where $1 \le i < x$. Then a) implies that}

\hspace{14.5mm}$w_{r+i+1}=w_{r+i}-\overline{w}(a_{r+i+1})
=p+1-\overline{w}(a_{r+1} \cdots {a}_{r+i+1})$.\vspace{0.7mm}

\noindent c) follows from b)  since $w_{r+x}=0$ and $0 < w_{r+i} < p+1$ for $1 \le i < x$.\vspace{0.7mm}

\noindent d)  $w(A_r)=w_r+k=k+p+1$. Since $w_{r}>w_{r+1}$ and $w_{r+x-1}>w_{r+x}$, then Observation 44.1 b) implies that $a_{r+1}=a_{r+x}=1$. Hence, $A_r$ starts with 1.
 By  using c) we get that $P=a_{r+1} \cdots a_{r+x}$ is a positive start string \mbox{of $A_r$} of order $p+1$
  that ends with 1. Hence, we get according to Section 43 \mbox{that $V(A_r) \in M_p^+$.}\vspace{2mm}

\hspace{45mm} {\large{\textbf{PART 7.}}\vspace{ 0.5mm}}

\noindent  We  study vector representations of infinite strings generated by symmetric shift registers. \mbox{The main results are  Proposition 45.4 and (46.1).}\vspace{10.5mm}

\noindent {\large{\textbf{45. The infinite vector representation.}}}\vspace{1mm}

\noindent Suppose  $A=a_1 \cdots a_n $ where $a_1=1$\hspace{0.2mm}.   Let $A^{\,\infty}=a_1 \cdots a_na_{n+1} \cdots $ be   generated from $A$ by the symmetric shift register $\theta$ with parameters  $k$, $p$ \mbox{and $n$}\hspace{0.1mm}, and suppose
$w(A)=k+p+1$\hspace{0.15mm}. Since $a_1=1$, then  $a_2+ \cdots+a_n=k+p$. Hence,  $a_{n+1}=a_1^{\,\prime}=0$.  As in Section 6 we decompose\vspace{0.6mm}

\noindent (45.1)  \hspace{18.8mm}  $A^{\,\infty}=1_{q_1}0_{q_2}1_{q_3}  \cdots $
  where $q_i \ge 1$ for $i \ge 1$\hspace{0.15mm}.\vspace{0.6mm}

\noindent The vector representation of $A^{\,\infty}$ is $V(A^{\,\infty})=(q_1,q_2, \cdots )$\hspace{0.2mm}.
 Next, we will  prove the following statement:\vspace{0.5mm}

 \noindent (45.2)  \hspace{2.5mm}  If $a_{r+1}=a_{r+2}= \cdots =a_{r+q}$ where $q \ge 1$ and $r \ge 0$\,, then $q \le n$\hspace{0.25mm}.\vspace{0.5mm}

 \noindent Suppose  $a_{r+1}=a_{r+2}= \cdots =a_{r+q}$ where  $r \ge 0$ and $q>n$\hspace{0.2mm}.
 Then we get \mbox{that $\theta(a_{r+1} \cdots a_{r+n})=a_{r+2} \cdots a_{r+n+1}=a_{r+1} \cdots a_{r+n}$}\hspace{0.1mm}.
 Hence, the period \mbox{of $A^{\,\infty}=a_1a_2 \cdots  $} is one.  This is a contradiction since
 $a_1=1$ and $a_{n+1}=0$\hspace{0.2mm}. Hence, (45.2) is true. By (45.2)  we get that\vspace{0.5mm}

 \noindent (45.3) \hspace{38.5mm}$1 \le q_j \le n$ for $j \ge 1$\hspace{0.2mm}.\vspace{0.5mm}

 \noindent  Since $A=a_1 \cdots a_n$ is a start string of $1_{q_1}0_{q_2}1_{q_3}  \cdots $  succeeded \mbox{by $a_{n+1}=0$\hspace{0.35mm},}
 then there exists an odd integer $J>0$ such that\vspace{0.5mm}

 \noindent (45.4)  \hspace{15mm} $A=1_{q_1} 0_{q_2}1_{q_3}0_{q_4} \cdots 1_{q_{J}}0_{e_{\hspace{0.12mm}0}}$ where $0 \le e_0 <  q_{J+1}$\hspace{0.2mm}.\vspace{0.5mm}

 \noindent As in Section 42 we let $A_r=a_{r+1} \cdots a_{r+n}$ for $r \ge 0$\hspace{0.35mm}. In particular, $A_0=A$\hspace{0.25mm}.\vspace{0.5mm}

 \noindent Suppose $j >0$ is an even vector  period of $V(A^{\,\infty})$. Then\vspace{0.5mm}

\noindent (45.5)  \hspace{11mm}$r=sum(V(A^{\,\infty}),j)=q_1+ \cdots+q_j$ is a period of $A^{\,\infty}$.\vspace{0.5mm}

\noindent In fact, $a_{\hspace{0.2mm}r+1}a_{\hspace{0.2mm}r+2} \cdots =1_{q_{j+1}}0_{q_{j+2}}1_{q_{j+3}}  \cdots=1_{q_1}0_{q_2}1_{q_3}0_{q_4}  \cdots=a_1a_2 \cdots $
where the first equality is true since $a_1 \cdots a_{\hspace{0.2mm}r}=1_{q_1}0_{q_2}1_{q_3} \cdots 0_{q_j}$\hspace{0.25mm}, the second is true since $q_{j+i}=q_i$ for $i \ge 1$\hspace{0.25mm}, and the last follows from (45.1).\vspace{2mm}

 \noindent {\textbf{Observation 45.1.}} Suppose  $r>0$ is a period of $A^{\,\infty}$. Then\vspace{0.3mm}

 \hspace{28mm}$a_r=0$\hspace{0.35mm}, $a_{r+1}=1$ and $w(A_r)=k+p+1$\hspace{0.2mm}.\vspace{0.5mm}

 \noindent {\textbf{Proof.}}  Since  $r$ is a period of $A^{\,\infty}$, then $a_{r+i}=a_i$ for $i \ge 1$\hspace{0.2mm}. Hence,  $A_{\hspace{0.2mm}r}=A$\hspace{0.2mm} and
 $a_{r+1}=a_1=1$\hspace{0.2mm}.  Moreover,
 $w(A_{\hspace{0.2mm}r})=w(A)=k+p+1$\hspace{0.1mm}.
 \mbox{If $a_r=1$,} then Observation 42.1 b) implies  that \mbox{$k \le w(A_r)<k+p+1$\hspace{0.1mm}.}\mbox{ This is a contradiction.}\vspace{2mm}

 \noindent {\textbf{Observation 45.2.}}  Suppose    $r$ is a period of $A^{\,\infty}$. Then there exists  an even vector period $j$ of $V(A^{\,\infty})$ such that $r=q_1+ \cdots +q_j$\hspace{0.3mm}.\vspace{0.7mm}

\noindent {\textbf{Proof.}} By Observation 45.1 we get that $a_{\hspace{0.2mm}r}=0$ and $a_{r+1}=1$\hspace{0.1mm}. By (45.1) there exists an even $j \ge 1$ such that $a_{\hspace{0.2mm}1} \cdots a_{\hspace{0.2mm}r}=1_{q_1}0_{q_2}1_{q_3} \cdots 0_{q_j}$\hspace{0.2mm}. Hence,

 \noindent  $r=q_1 + \cdots +q_j$ and $ 1_{q_{j+1}}0_{q_{j+2}}1_{q_{j+3}} \cdots =a_{r+1}a_{r+2} \cdots =a_1 a_2 \cdots =1_{q_1}0_{q_2}1_{q_3}  \cdots$.\vspace{0.4mm}

  \noindent Then $q_{j+i}=q_i$ for $i \ge 1$\hspace{0.1mm}. That means, $j$ is an even vector period  of $V(A^{\,\infty})$.\vspace{20mm}

  \noindent {\textbf{Observation 45.3.}} $V(A^{\,\infty})$ has an even vector period.\vspace{0.4mm}

  \noindent {\textbf{Proof.}} Since $A^{\,\infty}$ has a period,  this follows from Observation 45.2.\vspace{2mm}

\noindent {\textbf{Proposition 45.4.}} If $j$ is the least even vector period of $V(A^{\,\infty})$, then\vspace{0.4mm}

\hspace{7mm}  $sum(V(A^{\,\infty}),j)=q_1+ \cdots+q_j$ is the minimal period
of  $A^{\,\infty}$.\vspace{0.5mm}

\noindent {\textbf{Proof.}}   Suppose $j$ is the least even vector period of $V(A^{\,\infty})$. By (45.5) we get that  $r=q_1+ \cdots+q_j$ is a period
of  $A^{\,\infty}$. Suppose $r^*$ is a period \mbox{of $A^{\,\infty}$.} By Observation 45.2 there exists an even vector period $i$ of $V(A^{\,\infty})$
such \mbox{that $r^*=q_1+ \cdots+q_{\hspace{0.2mm}i}$\hspace{0.2mm}.} By the minimality of $j$ we get that $i \ge j$ \mbox{and $r^* \ge r$.} Hence,  $r=q_1+ \cdots+q_j$ is the least period
of  $A^{\,\infty}$.\vspace{2mm}

\noindent {\large{\textbf{46. A crucial representation.}}}\vspace{0.5mm}

\noindent Let  $A^{\infty}=A^{\infty}_p(Q)$ where $Q \in M_p$ and $p \ge 0$. Then $A^{\infty}$  is generated as in Section 18 \mbox{from $A=A(Q)$} by the symmetric shift register with respect to the parameters $k=w(A)-(p+1)$, $p$ \mbox{and $n=length(A)$.} Then $w(A)=k+p+1$ \mbox{and $A$} starts with 1.
In the next section we will prove that\vspace{0.5mm}

\noindent (46.1) \hspace{40mm}$C_p^{\,\infty}(Q)=V(A^{\infty})$.\vspace{0.5mm}

\noindent Suppose  $A^{\infty}=a_1 \cdots a_na_{n+1} \cdots$ and  $V(A^{\infty})=(q_1, q_2, q_3,\cdots )$. Then\vspace{0.4mm}

\noindent (46.2)  \hspace{16mm} $q_{j+1}>0$ for $j \ge 0$, and   $A^{\infty}=1_{q_1}0_{q_2}1_{q_3} \cdots $  \vspace{0.5mm}

\noindent
 Let
$r_0=0$
and $r_j=q_1+ \cdots +q_j$ for $j \ge 1 $.  Then\vspace{0.5mm}

\noindent (46.3)\hspace{16.4mm}$r_0=0 < r_1 < r_2 \cdots$ and $r_{j+1}=r_j+q_{j+1}$ for $j \ge 0$.
\vspace{0.5mm}

\noindent Let $H_{j+1}=a_{r_j+n+1} \cdots a_{r_{j+1}+n}=a_{r_j+n+1} \cdots a_{r_{j}+n+q_{j+1}}$ for $j \ge 0$. By (46.3)  we get that $r_0+n=n < r_1+n < r_2+n \cdots$ and \vspace{0.5mm}

\noindent (46.4)\hspace{38mm}   $A^{\infty}=AH_1H_2H_3 \cdots$.\vspace{0.5mm}

\noindent Since $length(A)=n$, then $A=a_1 \cdots a_n$.
Let $A_r=a_{r+1} \cdots a_{r+n}$  for $r \ge 0$. In particular, $A_0=A$. Moreover, let $w_r=w(A_r)-k$ for $r \ge 0$.\vspace{0.5mm}

\noindent We note that $w_{r_0}=w_0=w(A_0)-k=w(A)-k=p+1$.
By (45.4)  there exists an odd integer $J>0$ such that\vspace{0.5mm}

 \noindent (46.5)  \hspace{18mm} $A=1_{q_1} 0_{q_2}1_{q_3}0_{q_4} \cdots 1_{q_{J}}0_{e_{\hspace{0.12mm}0}}$ where $0 \le e_0 <  q_{J+1}$\hspace{0.2mm}.\vspace{0.6mm}

  \noindent Since $A=A(Q)$, then Observation 7.2 imply that $Q=(q_1, \cdots,q_J,e_0)$.\vspace{2mm}

  \noindent {\textbf{Observation 46.1.}} Suppose $j \ge 0$. Then\vspace{0.3mm}

   \hspace{14.5mm} $a_{r_j+1} \cdots a_{r_{j+1}}=a_{r_j+1} \cdots a_{r_{j}+q_{j+1}}=1_{q_{j+1}}$ if $j$ is even,

   \hspace{14.5mm} $a_{r_j+1} \cdots a_{r_{j+1}}=a_{r_j+1} \cdots a_{r_{j}+q_{j+1}}=0_{q_{j+1}}$ if $j$ is odd.\vspace{0.5mm}

   \noindent {\textbf{Proof.}} Since $r_0=0$, $r_1=q_1$ and $1_{q_1}$ starts $A^{\infty}=a_1a_2\cdots $, then\vspace{0.5mm}

   \hspace{34mm}$a_{r_0+1} \cdots a_{r_{1}}=a_{1} \cdots a_{q_{1}}=1_{q_{1}}$,\vspace{0.5mm}

   \noindent and the result is true for $j =0$. Suppose the result is true for $j$
   where $j \ge 0$.\vspace{0.5mm}

   \noindent If $j$ is even, then $a_{r_j+1} \cdots a_{r_{j+1}}=1_{q_{j+1}}$ is succeeded by
   $0_{q_{j+2}}$. Hence,\vspace{0.5mm}

  \hspace{20mm} $a_{r_{j+1}+1} \cdots a_{r_{j+2}}= a_{r_{j+1}+1} \cdots a_{r_{j+1}+q_{j+2}}=0_{q_{j+2}}$.\vspace{0.5mm}

   \noindent If $j$ is odd, then $a_{r_j+1} \cdots a_{r_{j+1}}=0_{q_{j+1}}$ is succeeded by
   $1_{q_{j+2}}$. Hence,\vspace{0.5mm}

   \hspace{20mm} $a_{r_{j+1}+1} \cdots a_{r_{j+2}}= a_{r_{j+1}+1} \cdots a_{r_{j+1}+q_{j+2}}=1_{q_{j+2}}$.\vspace{0.5mm}

   \noindent We conclude that the result is true for $j+1$.\vspace{2mm}

    \noindent {\textbf{Observation 46.2.}} a)   $a_{r_j+i}=1$ for $1 \le i \le q_{j+1}$ if $j\ge 0 $ is even.\vspace{0.5mm}

\noindent b)  $a_{r_j+i}=0$ for $1 \le i \le q_{j+1}$ if $j\ge 0 $ is odd.\vspace{0.5mm}

\noindent c)  $j-1$ is odd and $a_{r_j}=a_{r_{j-1}+q_j}=0$ if $j \ge 1$ is even.\vspace{0.5mm}

\noindent d) $j-1$ is even and $a_{r_j}=a_{r_{j-1}+q_j}=1$ if $j \ge 1$ is odd.\vspace{0.5mm}

 \noindent {\textbf{Proof.}} a) and b) follow from Observation 46.1. Moreover, c) and d) follow from a) and b).\vspace{2mm}

 \noindent {\textbf{Observation 46.3.}} Suppose $j \ge 0$. Then $0 <  w_{r_{j}} \le  p+1$ if $j$ is even, \mbox{and $0 \le  w_{r_{j}} < p+1$} if $j$ is odd.\vspace{0.5mm}

\noindent {\textbf{Proof.}} This is true for $j=0$ since  $w_{r_{0}} = p+1$.  If $j \ge 1$ is even, then   Observation 42.1 c) and 46.2 c) imply that
$a_{r_{j}}=0$ and $0 <  w_{r_{j}} \le  p+1$.\vspace{0.5mm}

\noindent Suppose $j \ge 1$ is odd. Then we get by  Observation 42.1 b) and 46.2 d)  \mbox{that $a_{r_{j}}=1$ and $0 \le   w_{r_{j}} <  p+1$.}\vspace{2mm}

\noindent {\large{\textbf{47. Deductions.}}}\vspace{0.5mm}

\noindent   Let $\lambda_j= w_{r_j}$ for $j \ge 0$. Then $\lambda_0=w_{r_0}=p+1$. Moreover, let\vspace{0.5mm}

\noindent (47.1) \hspace{8.9mm} $s_{j+1}=min\{q_{j+1},\lambda_j\}=min\{q_{j+1},w_{r_j}\}$ if $j \ge 0$ is even,

\noindent (47.2) \hspace{1mm} $s_{j+1}=min\{q_{j+1},p+1-\lambda_j\}=min\{q_{j+1},p+1-w_{r_j}\}$ if $j \ge 0$ is odd,

\noindent (47.3) \hspace{31mm} $e_{j+1}=q_{j+1}-s_{j+1}$ for $j \ge 0$.\vspace{1.5mm}

\noindent {\textbf{Observation 47.1.}} Suppose $j \ge 0$. Then $0 < \lambda_j \le p+1$ if $j$ is even, \mbox{and $0 \le  \lambda_j < p+1$} if $j$ is odd.\vspace{0.5mm}

\noindent {\textbf{Proof.}} Follows from Observation 46.3 since $\lambda_j= w_{r_j}$ for $j \ge 0$.\vspace{2mm}

\noindent {\textbf{Observation 47.2.}} $1 \le s_{j+1} \le q_{j+1}$ and  $e_j \ge 0$ for $j \ge 0$.\vspace{0.6mm}

\noindent {\textbf{Proof.}} By
(46.2)  we get that  $q_{j+1}>0$ for $j \ge 0$. If $j \ge 0$ is even,  then Observation 47.1 and (47.1) imply that $\lambda_j>0$ and $1 \le s_{j+1} \le q_{j+1}$. \mbox{If $j \ge 0$} is odd, then we get according to Observation 47.1 and (47.2)  that
$\lambda_j \le p$ and $1 \le s_{j+1} \le q_{j+1}$.
 Hence, we get by (47.3) that  $e_{j+1}=q_{j+1}-s_{j+1} \ge 0$ \mbox{for $j \ge 0$.} Finally, \mbox{by (46.5)} we get that $e_0 \ge 0$.\vspace{2mm}

\noindent {\textbf{Observation 47.3.}} Suppose $j \ge 0$ is even. Then

\noindent a) $a_{r_j+n+1} \cdots a_{r_j+n+q_{j+1}}=0_{s_{j+1}}1_{q_{j+1}-s_{j+1}}$
and $w_{r_j+q_{j+1}}=w_{r_j}-s_{j+1}$.

\noindent b) $H_{j+1}=a_{r_j+n+1} \cdots a_{r_{j+1}+n}=0_{s_{j+1}}1_{e_{j+1}}$ and $\lambda_{j+1}=\lambda_j-s_{j+1}$.\vspace{10.98mm}

\noindent {\textbf{Proof.}}
a) follows from  Proposition 42.2 c), Observation 46.2 a) and (47.1).\vspace{0.4mm}

\noindent b) We note that   $r_{j+1}=r_j+q_{j+1}$, $\lambda_{j+1}=w_{r_{j+1}}=w_{r_j+q_{j+1}}$,
$\lambda_{j}=w_{r_j}$ \mbox{and
$e_{j+1}=q_{j+1}-s_{j+1}$.} Hence, b) follows from a).\vspace{1.6mm}

\noindent {\textbf{Observation 47.4.}} Suppose $j \ge 0$ is odd. Then

\noindent a) $a_{r_j+n+1} \cdots a_{r_j+n+q_{j+1}}=1_{s_{j+1}}0_{q_{j+1}-s_{j+1}}$
and $w_{r_j+q_{j+1}}=w_{r_j}+s_{j+1}$.

\noindent b) $H_{j+1}=a_{r_j+n+1} \cdots a_{r_{j+1}+n}=1_{s_{j+1}}0_{e_{j+1}}$ and $\lambda_{j+1}=\lambda_j+s_{j+1}$.\vspace{0.6mm}

\noindent {\textbf{Proof.}}
a) follows from  Proposition 42.3 c), Observation 46.2 b) and (47.2).\vspace{0.4mm}

\noindent b) We note that  $r_{j+1}=r_j+q_{j+1}$, $\lambda_{j+1}=w_{r_{j+1}}=w_{r_j+q_{j+1}}$,
$\lambda_{j}=w_{r_j}$ \mbox{and
$e_{j+1}=q_{j+1}-s_{j+1}$.} Hence, b) follows from a).\vspace{1.6mm}

\noindent {\textbf{Observation 47.5.}} $q_{J+j+1}=e_j+s_{j+1}$ for $j \ge 0$.\vspace{0.6mm}

\noindent {\textbf{Proof.}}
 Observation 47.3 b) and 47.4 b) imply that\vspace{0.5mm}

 \hspace{13mm}$H_1=0_{s_1}1_{e_1}$, $H_2=1_{s_2}0_{e_2}$,
$H_3=0_{s_3}1_{e_3}, \,H_4=1_{s_4}0_{e_4}, \cdots$.\vspace{0.5mm}

\noindent Hence, we get from (46.4), (46.5) and Observation 47.2   that\vspace{0.4mm}

\noindent \hspace{4mm}$A^{\infty}=AH_1H_2H_3H_4 \cdots =1_{q_1}0_{q_2} \cdots 0_{q_{J-1}}1_{q_{J}}0_{e_0}0_{s_1}1_{e_1}1_{s_2}
 0_{e_2}0_{s_3}1_{e_3}1_{s_4}0_{e_4} \cdots $

\noindent \hspace{0.5mm}$=1_{q_1}0_{q_2} \cdots 0_{q_{J-1}}1_{q_{J}}0_{e_0+s_1}1_{e_1+s_2}0_{e_2+s_3} 1_{e_3+s_4}\cdots $ where $e_i+s_{i+1}>0$ for $i \ge 0$.\vspace{0.6mm}

 \noindent  According to  (46.2) we get that $A^{\infty}=1_{q_1}0_{q_2} \cdots 0_{q_{J-1}}1_{q_{J}}0_{q_{J+1}}1_{q_{J+2}}\cdots$. Hence,\vspace{0.5mm}

 \hspace{17mm}  $q_{J+1}=e_0+s_1,\, q_{J+2}=e_1+s_2,\, q_{J+3}=e_2+s_3,\, \cdots$.\vspace{1.4mm}

 \noindent {\textbf{Proof of (46.1).}} We get that\vspace{0.5mm}

\noindent (47.4)\hspace{7mm} $s_{j+1}=min\{q_{j+1},\lambda_j\}$ and $\lambda_{j+1}=\lambda_j-s_{j+1}$ if
$j \ge 0$ is even,

\noindent (47.5)\hspace{2mm} $s_{j+1}=min\{q_{j+1},p+1-\lambda_j\}$ and $\lambda_{j+1}=\lambda_j+s_{j+1}$ if $j \ge 0$ is odd,

\noindent (47.6)\hspace{31mm} $e_{j+1}=q_{j+1}-s_{j+1}$ for $j \ge 0$,\vspace{0.2mm}

\noindent (47.7) \hspace{29.5mm} $q_{J+j+1}=e_{j}+ s_{j+1}$ for $j \ge 0$,\vspace{0.5mm}

\noindent  where (47.4) follows from (47.1) and Observation 47.3 b), (47.5) follows \mbox{from (47.2)} and Observation 47.4 b), (47.6) from (47.3) and (47.7) from Observation 47.5.
Since $\lambda_0=p+1$ and $Q=(q_1, \cdots,q_J,e_0)$, \mbox{then we get} \mbox{from (47.4), $\cdots$, (47.7)}  that $C_p^{\,\infty}(Q)=(q_1,q_2,\cdots )$. Hence, $C_p^{\,\infty}(Q)=V(A^{\infty})$.\vspace{2.2mm}

\hspace{45mm} {\large{\textbf{PART 8.}}  \vspace{ 1 mm}}

\noindent {\large{\textbf{48. Uniqueness properties.}}}\vspace{1mm}

\noindent Suppose $Q^{\,\infty}=(q_1,q_2, \cdots)$ is complete. Let $\pi(Q^{\,\infty})=(q_1^*,q_2^*,  \cdots)$ be the contraction vector and  $D^{\,\infty}=(d_1,d_2, \cdots)$  the distance vector of $Q^{\,\infty}$. \mbox{Let $\tau$} be the distance function of $Q^{\,\infty}$.
  Let $r_0, r_1, \cdots$  and $t_1, t_2, \cdots$ be the $r$-indexes and $t$-indexes of $Q^{\,\infty}$. We also let  $\beta_i=\tau(r_i)$ and $y_i=y(r_i)$ for $i \ge 0$.
Moreover, let $t_{max}(r)$ and $next(r)$ be as in Section 33. Suppose  $j \ge 1$\hspace{0.1mm}. We will prove in Section 49 and 50 that\vspace{0.5mm}

 \noindent (48.1)\hspace{6mm} $q_{r_j+r}=q_r$ for $r \ge 1$ $\Leftrightarrow$
 $ d_{y_j+i}=d_{\hspace{0.1mm}i}+\beta_j$ and $q_{j+i}^*=q_i^*$ for $i \ge 1$\hspace{0.1mm}.\vspace{20.6mm}

 \noindent  By (37.1) we get that $y_j=r(D^{\,\infty},\beta_j)$. Hence, according to (35.7) we get that (48.1) is equivalent to\vspace{0.5mm}

 \noindent (48.2) \hspace{5.5mm}$r_j$ is a vector period of $Q^{\,\infty}$ $\Leftrightarrow$    $j$ is a vector period of $\pi(Q^{\,\infty})$

\hspace{9.8mm} and $\beta_j$ is a progression parameter of $D^{\,\infty}$.\vspace{2.5mm}

\noindent {\large{\textbf{49.  Deductions.}}}\vspace{0.7mm}

\noindent  Suppose $j \ge 1$. Moreover,  suppose $q_{r_j+r}=q_r$ for $r \ge 1$. We will   prove that
$ d_{y_j+i}=d_{\hspace{0.1mm}i}+\beta_j$ and $q_{j+i}^*=q_i^*$ for $i \ge 1$\hspace{0.1mm}. We refer to \mbox{Observation 49.5 and 49.6.}\vspace{2mm}

\noindent {\textbf{Observation 49.1.}} $\tau(r_j+i)=\beta_j+\tau(i)$  for $i \ge 0$\hspace{0.2mm}.

\noindent {\textbf{Proof.}} If $i=0$, this is trivial. Otherwise, if $i >0$, then  (32.8)  implies that

$\tau(r_j+i)=\tau(r_j)+\delta(q_{r_j+1}, \cdots, q_{r_j+i})=\beta_j+\delta(q_{1}, \cdots, q_{i})=\beta_j+\tau(i)$\hspace{0.1mm}.\vspace{2mm}

 \noindent {\textbf{Observation 49.2.}}
 $r_{j+i}=r_j+r_i$ for $i \ge 0$, and  $t_{j+i}=t_{i}$ for $i \ge 1$.\vspace{0.5mm}

 \noindent {\textbf{Proof.}}   Since $r_0=0$\hspace{0.2mm}, then $r_{j+0}=r_j+r_0$. Suppose $r_{j+i}=r_j+r_i$ \mbox{where $i \ge 0$.}
 By  (33.3)  we get that

     \hspace{17mm} $q_{r_{j+i}+2x}=q_{r_j+r_i+2x}=q_{r_i+2x}=1$ for $1 \le x \le t_{i+1}$,

    \hspace{17mm} $q_{r_{j+i}+2t_{i+1}+2}=q_{r_j+r_i+2t_{i+1}+2}=q_{r_i+2t_{i+1}+2}>1$.\vspace{0.7mm}

 \noindent Hence,   $t_{i+1}=t_{max}(r_{j+i})=t_{j+i+1}$ where the last equality follows from the definition of $t_{j+i+1}$ in Section 33. Then we get by (33.2)  that\vspace{0.5mm}

 \hspace{10mm} $r_{j+i+1}=r_{j+i}+2t_{i+1}+1=r_j+r_i+2t_{i+1}+1=r_j+r_{i+1}$.\vspace{2mm}

\noindent {\textbf{Observation 49.3.}} $y_{j+i}=y_j+y_i$ for $i \ge 0$.\vspace{0.5mm}

\noindent {\textbf{Proof.}} As in Section 37 we get that $y_0=0$. Hence, the result is true \mbox{for $i=0$.}
Suppose $y_{j+i}=y_j+y_i$ where $i \ge 0$. Then Observation 37.1 b) \mbox{and 49.2} imply that
$y_{j+i+1}=y_{j+i}+t_{j+i+1}=y_j+y_i+t_{i+1}=y_j+y_{i+1}$.\vspace{2mm}

\noindent {\textbf{Observation 49.4.}} If  $1 \le m \le t_{i+1}$ where $i \ge 0$, then
$d_{y_j+y_i+m}=d_{y_i+m}+\beta_j$.\vspace{0.5mm}

\noindent {\textbf{Proof.}}  Suppose $1 \le m \le t_{i+1}$ where $i \ge 0$.
Hence, we get by using Observation 37.2, 49.1, 49.2 and 49.3   that $1 \le m \le t_{j+i+1}$ and $d_{y_j+y_i+m}=d_{y_{j+i}+m}=\tau(r_{j+i}+2m)=\tau(r_j+r_i+2m)
=\tau(r_{i}+2m)+\beta_j=d_{y_i+m}+\beta_j$.\vspace{2mm}

\noindent {\textbf{Observation 49.5.}} $q_{j+i}^*=q_i^*$   for $i \ge 1$\hspace{0.1mm}.\vspace{0.5mm}

\noindent {\textbf{Proof.}}  If $i \ge 1$, then   (33.6),   Observation 49.1 and 49.2   imply that\vspace{0.5mm}

\noindent  \ $q_{j+i}^*=\tau(r_{j+i})-\tau(r_{j+i-1})=\tau(r_j+r_{i})-\tau(r_j+r_{i-1})
=\tau(r_{i})-\tau(r_{i-1})=q_{i}^*$.\vspace{2mm}

\noindent {\textbf{Observation 49.6.}}  $ d_{y_j+x}=d_{\hspace{0.1mm}x}+\beta_j$  for $x \ge 1$\hspace{0.1mm}.\vspace{0.5mm}

\noindent {\textbf{Proof.}}    Let $x \ge 1$. According to Observation  37.4 a) there exists $i \ge 0$ such
\mbox{that $y_i < x \le y_{i+1}$.} By Observation 37.1 b) we get that $y_{i+1}=y_i+t_{i+1}$.
Hence,  $x=y_i+m$  \mbox{where $i \ge 0$} and $1 \le m \le t_{i+1}$.
Then Observation 49.4 implies that $d_{y_j+x}=d_{y_j+y_i+m}=d_{y_i+m}+\beta_j=d_x+\beta_j$.\vspace{2mm}

\noindent {\large{\textbf{50.  Additional deductions.}}}\vspace{0.7mm}

 \noindent In this section we suppose $j \ge 1$ and\vspace{0.3mm}

 \noindent (50.1) \hspace{22mm} $ d_{y_j+i}=d_{\hspace{0.1mm}i}+\beta_j$ and $q_{j+i}^*=q_i^*$ for $i \ge 1$.\vspace{0.5mm}

 \noindent We will prove in the end of this section that $q_{r_j+r}=q_r$ for $r \ge 1$.
  We refer to Observation 50.7 b). By Observation 37.5 a) we get that  $\beta_j>0$.
  \mbox{Since
   $d_{y_j+1}=d_1+\beta_j>d_1$, then $y_j \ge 1$.} By (37.1) we get that $y_i=r(D^{\,\infty},\beta_i)$ \mbox{for $i \ge 0$.} Suppose $i \ge 0$. Then we get by (35.5) and (35.6) that\vspace{0.5mm}

   \noindent (50.2) \hspace{12.6mm}$d_{y_i} \le \beta_i< d_{y_i+1}$ if $y_i>0$,\, and $\beta_i < d_{y_i+1}$ if $y_i=0$.\vspace{0.4mm}

   \noindent In particular, since $y_j>0$, then $d_{y_j} \le \beta_j < d_{y_j+1}$.\vspace{2mm}

\noindent {\textbf{Observation 50.1.}} $\beta_{j+i}=\beta_j+\beta_i$ for $i \ge 0$.\vspace{0.4mm}

 \noindent {\textbf{Proof.}} By Section 37 we get that $\beta_0=0$. Then $\beta_{j+0}=\beta_j+\beta_0$. Next, we suppose  $\beta_{j+i}=\beta_j+\beta_i$ where $i \ge 0$.  Then  (50.1) and Observation 37.5 b) imply that
 $\beta_{j+i+1}=\beta_{j+i}+q_{j+i+1}^*=\beta_j+\beta_i+q_{i+1}^*=\beta_j+\beta_{i+1}$.\vspace{2mm}

  \noindent {\textbf{Observation 50.2.}} Suppose $i \ge 0$. Then $y_{j}+y_i=y_{j+i}$.\vspace{0.3mm}

  \noindent {\textbf{Proof.}} By Observation 37.4 a) and 37.5 a)  we get that  $0 <y_j \le y_{j+i}$
  \mbox{and $\beta_j \le \beta_{j+i}$.} By
  (50.1), (50.2) and Observation 50.1  we get that

 \hspace{24mm} $\beta_{j+i}=\beta_j+\beta_i < \beta_j+d_{y_i+1}=d_{y_j+y_i+1}$,\vspace{0.5mm}

 \hspace{20.5mm} $d_{y_j+y_i}=d_{y_i}+\beta_j \le \beta_i+\beta_j=\beta_{j+i}$   if $y_i >0$,\vspace{0.5mm}

  \hspace{28.5mm} $d_{y_j+y_i}=d_{y_j} \le \beta_j \le \beta_{j+i}$   if $y_i =0$.\vspace{0.5mm}

  \noindent We  conclude that $d_{y_j+y_i} \le \beta_{j+i} < d_{y_j+y_i+1}$. Hence, we get from  (35.5)  \mbox{that
  $r(D^{\,\infty},\beta_{j+i})=y_j+y_i$.} Since $y_{j+i}=r(D^{\,\infty},\beta_{j+i})$, then $y_{j+i}=y_j+y_i$.\vspace{2mm}

\noindent {\textbf{Observation 50.3.}}   $t_{j+i+1}=t_{i+1}$  and $r_{j+i}=r_j+r_i$ for $i \ge 0$.\vspace{0.3mm}

 \noindent {\textbf{Proof.}} Let $i \ge 0$. By Observation 37.1 b) and  50.2    we get that

\hspace{6mm}$t_{j+i+1}=y_{j+i+1}-y_{j+i}=(y_j+y_{i+1})-(y_j+y_i)=y_{i+1}-y_i=t_{i+1}$.\vspace{0.5mm}

\noindent Moreover, Observation 37.3 b) and 50.2  imply that

\hspace{12mm} $r_{j+i}=2y_{j+i}+(j+i)=2y_j+2y_i+j+i=r_j+r_i$.\vspace{2mm}

\noindent {\textbf{Observation 50.4.}} If  $i \ge 0$, then
 $\tau(r_{j+i})=\beta_{j+i}=\beta_j+\beta_i=\tau(r_i)+\beta_j$.\vspace{0.3mm}

 \noindent {\textbf{Proof.}}  The result follows from Observation 50.1.\vspace{1.4mm}

 \noindent {\textbf{Observation 50.5.}}  Suppose $1 \le x \le t_{i+1}$ and $i \ge 0$. Then

 $\tau(r_{j+i}+2x)=d_{y_{j+i}+x}=d_{y_j+y_i+x}=d_{y_i+x}+\beta_j=\tau(r_i+2x)+\beta_j$,

 $\tau(r_{j+i}+2x-1)=\tau(r_{j+i}+2x)=\tau(r_{i}+2x)+\beta_j=\tau(r_{i}+2x-1)+\beta_j$.\vspace{0.5mm}

 \noindent {\textbf{Proof.}} By Observation 50.3 we get that $t_{j+i+1}=t_{i+1}$. Hence, the results follow from (50.1), Observation 37.2 and 50.2.\vspace{2mm}

 \noindent {\textbf{Observation 50.6.}}  $\tau(r_{j}+r_i+m)=\tau(r_{j+i}+m)=\tau(r_i+m)+\beta_j$ \mbox{for $0 \le m < 2t_{i+1} +1$ and $i \ge 0$.}\vspace{0.3mm}

  \noindent {\textbf{Proof.}} The first equality follows from Observation 50.3, and the last
   equality follows from  Observation 50.4 and 50.5.\vspace{2mm}

  \noindent {\textbf{Observation 50.7.}} a) $\tau(r_{j}+r)=\tau(r)+\beta_j$ for $r \ge 0$.

  \noindent b) $q_{r_j+r}=q_{r}$ for $r \ge 1$.\vspace{0.5mm}

  \noindent {\textbf{Proof.}} a) By (33.1)  there exists $i \ge 0$ such that $r_i \le r < r_{i+1}$. By (33.2) we also get that  $r=r_i+m$ where $0 \le m < 2t_{i+1} +1$. Then Observation 50.6 implies that
  $\tau(r_{j}+r)=\tau(r_{j}+r_i+m)=\tau(r_i+m)+\beta_j=\tau(r)+\beta_j$.\vspace{0.5mm}

  \noindent b) Let $r \ge 1$. It is sufficient to prove that $q_{\,r_j+r}^-=q_{r}^-$. By a) we get that\vspace{0.5mm}

  \hspace{2mm} $q_{\,r_j+r}^-=\tau(r_j+r)-\tau(r_j+r-1)=
  (\tau(r)+\beta_j)-(\tau(r-1)+\beta_j)$

  \hspace{13.5mm}$=\tau(r)-\tau(r-1)=q_{\,r}^-$.\vspace{2mm}

 \hspace{45mm} {\large{\textbf{PART 9.}}  \vspace{ 0.5 mm}}

 \noindent  We will  derive a lot of auxiliary results about  shift symmetric vectors.\vspace{1.4mm}

\noindent {\large{\textbf{51. Assumptions and basic properties.}}}\vspace{1mm}

\noindent  In this part  we suppose $p \ge 0$, $Q=(q_1, \cdots, q_J,e_0) \in M$,  and   $Q^{\,\infty}=C_p^{\infty}(Q)$ is  the   shift  symmetric  vector  generated by $Q$ with respect to  $p$\hspace{0.2mm}.
We also suppose    $Q^{\,\infty}=(q_1,q_2, \cdots )$. Let $(s_1, s_2, \cdots )$, $(e_0,e_1, \cdots )$ and $(\lambda_0, \lambda_1, \cdots )$   be  the associated sequences.
Then $\lambda_0=p+1$. By Observation 28.4 we get \mbox{that $q_{j+1}>0$}  for $j \ge 0$.
Since $Q \in M$, then $e_0 \ge 0$. Let $\alpha=\delta(Q)+1$.\vspace{2mm}

\noindent {\textbf{Observation 51.1.}} a) $1 \le s_{j+1} \le q_{j+1}$ for $j \ge 0$.

\noindent b) If $q_r=1$ where $r \ge 1$, then $s_r=1$ and $e_r=0$.\vspace{0.5mm}

\noindent {\textbf{Proof.}} a) follows from (28.5) and Observation 28.4.

 \noindent b) If $q_r=1$ where $r \ge 1$, then a) and (28.3) imply that $s_r=1$ and $e_r=0$.\vspace{1mm}

\noindent {\textbf{Observation 51.2.}} Suppose $s_{j+1}<q_{j+1}$ where $j \ge 0$.

\noindent a) If $j$ is even, then $s_{j+1}=\lambda_j$ and $\lambda_{j+1}=\lambda_j-s_{j+1}=0$.

\noindent b) If $j$ is odd, then $s_{j+1}=p+1-\lambda_j$ and $\lambda_{j+1}=\lambda_j+s_{j+1}=p+1$.\vspace{0.6mm}

\noindent {\textbf{Proof.}} a) and b) follow from (28.1) and (28.2) respectively.\vspace{1mm}

 \noindent {{\textbf{Observation  51.3.}}}
 If $e_r>0$ or $s_{r+1} \ge 2$ where $r \ge 0$\hspace{0.3mm}, then $q_{J+r+1}\ge 2$\hspace{0.2mm}.\vspace{0.5mm}

 \noindent {{\textbf{Proof.}}}
 By (28.7) and Observation 28.4 we get that $e_r \ge 0$ and $s_{r+1} \ge 1$. Hence, we get  from (28.4)
 that $q_{J+r+1}=e_r+s_{r+1}\ge 2$ if $e_r>0$ or $s_{r+1} \ge 2$ .\vspace{2mm}

\noindent {\textbf{Observation 51.4.}} $s_1>1$ and $q_{J+1}>1$ if $q_1>1$ and $p>0$.\vspace{0.5mm}

 \noindent {{\textbf{Proof.}}} Suppose $p>0$. Then $\lambda_0=p+1>1$. If $q_1>1$, then (28.1) and Observation 51.3 imply that $s_{1}=min\{q_{1},\lambda_0\} >1$
 and $q_{J+1}>1$.\vspace{10.5mm}

\noindent {\textbf{Observation 51.5.}} Suppose $s_{r+2}=1$ and $\lambda_{r+1}=0$ where $r \ge 0$ is even. Then $\lambda_{r+2}=1, s_{r+3}=1$ and
$\lambda_{r+3}=0$.\vspace{0.4mm}

\noindent {\textbf{Proof.}} Since $r+1$ is odd, then (28.2) implies  that  \mbox{$\lambda_{r+2}=\lambda_{r+1}+s_{r+2}=1$.} Since $r+2$ is even and $q_{r+3}>0$, then we get by  (28.1)  that

\hspace{12mm}$s_{r+3}=min\{q_{r+3}, \lambda_{r+2}\}=1$ and $\lambda_{r+3}=\lambda_{r+2}-s_{r+3}=0$.\vspace{1mm}

\noindent {\textbf{Observation 51.6.}} Suppose $s_{r+2}=1$ and $\lambda_{r+1}=p+1$ where $r \ge 0$ is odd. Then $\lambda_{r+2}=p, s_{r+3}=1$ and
$\lambda_{r+3}=p+1$.

\noindent {\textbf{Proof.}} Since $r+1$ is even,  $r+2$ is odd and $q_{r+3}>0$, then (28.1) and (28.2) imply
 that $\lambda_{r+2}=\lambda_{r+1}-s_{r+2}=p$,
$s_{r+3}=min\{q_{r+3}, p+1-\lambda_{r+2}\}=1$ \mbox{and $\lambda_{r+3}=\lambda_{r+2}+s_{r+3}=p+1$.}\vspace{2mm}

\noindent {\large{\textbf{52. Properties of the distance function.}}}\vspace{1mm}

\noindent In this section we let  $\tau$ be the distance function of $Q^{\,\infty}$ defined by $\tau(0)=0$ and $\tau(r)=q_1^-+ \cdots + q_r^-$ for $r \ge 1$.
 Since $Q=(q_1, \cdots, q_J,e_{\hspace{0.2mm}0})$\hspace{0.2mm} \mbox{where $q_i>0$} for $1 \le i \le J$ and $e_0 \ge 0$, then\vspace{0.7mm}

\noindent (52.1)  \hspace{5mm} $\alpha=\delta(Q)+1=q_1^-+ \cdots+q_J^-+e_0^-+1=\tau(J)+e_{\hspace{0.2mm}0} \ge \tau(J)$.\vspace{0.5mm}

\noindent By (28.8) and (28.9) we get that\vspace{0.4mm}

\noindent (52.2)\hspace{5.9mm}$q^-_{J+i+1}-q^-_{i+1}=e_i+s_{i+1}^--(s_{i+1}^-+e_{i+1})=e_i-e_{i+1}$ for $i \ge 0$\,.\vspace{2mm}

\noindent {{\textbf{Observation  52.1.}}}  $\tau(J+i)-\tau(i)=\alpha-e_i$ for $i \ge 0$.\vspace{0.5mm}

\noindent {{\textbf{Proof.}}}  Since $\tau(0)=0$, then (52.1) implies that $\tau(J)- \tau(0)=\tau(J)=\alpha-e_0$. Next, we suppose $\tau(J+i)- \tau(i)=\alpha-e_i$ where $i \ge 0$.\vspace{0.5mm}

\noindent By (52.2)
we get that
$\tau(J+i+1)-\tau(i+1)=\tau(J+i)+q^-_{J+i+1}-(\tau(i)+q^-_{i+1})$
is equal to $\tau(J+i)-\tau(i)+e_{i}-e_{i+1}=\alpha-e_i +e_i-e_{i+1}=\alpha-e_{i+1}$.\vspace{2mm}

\noindent {{\textbf{Observation  52.2.}}}  If $q_i=1$ where $i \ge 1$, then $\tau(J+i)=\tau(i)+\alpha$.\vspace{0.3mm}

\noindent {{\textbf{Proof.}}}  Suppose  $q_i=1$ where $i \ge 1$. By  Observation 51.1 b) and 52.1 we get that $e_i=0$ \mbox{and $\tau(J+i)=\tau(i)+\alpha$.}\vspace{1.5mm}

 \noindent {\large{\textbf{53. Properties of the $\lambda$\,-\,parameters.}}}\vspace{1mm}

 \noindent In this section we suppose $r \ge 0$, $t \ge 0$ and $q_{r+2i}=1$ for $1 \le i \le t$.
 By Observation 51.1 b) we get that $s_{r+2i}=1$ for $1 \le i \le t$.\vspace{2mm}

  \noindent {\textbf{Observation 53.1.}} Suppose $0 \le i \le 2t$ where $i$ is odd.
  Then  $s_{r+i+1}=1$.\vspace{0.4mm}

   \noindent {\textbf{Proof.}} Since $i$ is odd and $i+1$ is even, then $1 \le i < 2t$
  and $i+1=2i^*$
  where  $1 \le i^* \le t$. Hence, $s_{r+i+1}=s_{r+2i^*}=1$.\vspace{2mm}

  \noindent {\textbf{Observation 53.2.}} Suppose $r$ is even. Let $1 \le i \le 2t+1$. Then

 \hspace{19mm} $\lambda_{r+i}=\lambda_r-(s_{r+1}^-+ \cdots +s_{r+i}^-)$ \,if $i$ is even,\vspace{0.2mm}

 \hspace{19mm}  $\lambda_{r+i}=\lambda_r-(s_{r+1}^-+ \cdots +s_{r+i}^-)-1$ if $i$ is odd.\vspace{0.6mm}

   \noindent {\textbf{Proof.}} By (28.1) we get that
   $\lambda_{r+1}=\lambda_r-s_{r+1}=\lambda_r-s_{r+1}^--1$. Suppose the result is true for $i$
   where $1 \le i < 2t+1$. We note that $r+i$ is even if and only if $i$ is even. If $i$ is even, then (28.1)  implies that\vspace{0.4mm}

   \hspace{6mm}$\lambda_{r+i+1}=\lambda_{r+i}-s_{r+i+1}=\lambda_r-(s_{r+1}^-+ \cdots +s_{r+i}^-)-s_{r+i+1}^--1$.\vspace{0.6mm}

  \noindent If $i$ is odd, then Observation 53.1 and (28.2) imply that $s_{r+i+1}=1$ and\vspace{0.4mm}

   \hspace{10.2mm}$\lambda_{r+i+1}=\lambda_{r+i}+s_{r+i+1}=\lambda_r-(s_{r+1}^-+ \cdots +s_{r+i}^-)-1+1$

  \hspace{5mm} $=\lambda_r-(s_{r+1}^-+ \cdots +s_{r+i}^-)=\lambda_r-(s_{r+1}^-+ \cdots +s_{r+i}^-)-s_{r+i+1}^-$.\vspace{0.4mm}

  \noindent Hence, the results are true for $i+1$.\vspace{2mm}

 \noindent {\textbf{Observation 53.3.}} Suppose $r$ is odd. Let $1 \le i \le 2t+1$. Then

 \hspace{18mm} $\lambda_{r+i}=\lambda_r+(s_{r+1}^-+ \cdots +s_{r+i}^-)$ \,if $i$ is even,

 \hspace{18mm}  $\lambda_{r+i}=\lambda_r+(s_{r+1}^-+ \cdots +s_{r+i}^-)+1$ if $i$ is odd.\vspace{0.4mm}

   \noindent {\textbf{Proof.}} By (28.2) we get that
   $\lambda_{r+1}=\lambda_r+s_{r+1}=\lambda_r+s_{r+1}^-+1$. Suppose the result is true for $i$
   where $1 \le i < 2t+1$. If $i$ is even, then  (28.2)  implies that\vspace{0.3mm}

  \hspace{6mm} $\lambda_{r+i+1}=\lambda_{r+i}+s_{r+i+1}=\lambda_r+(s_{r+1}^-+ \cdots +s_{r+i}^-)+s_{r+i+1}^-+1$\vspace{0.6mm}

  \noindent where we have used that $r+i$ is odd. If $i$ is odd, then $r+i$ is even. Hence, Observation 53.1 and (28.1) imply that $s_{r+i+1}=1$ and\vspace{0.4mm}

  \hspace{9.9mm} $\lambda_{r+i+1}=\lambda_{r+i}-s_{r+i+1}=\lambda_r+(s_{r+1}^-+ \cdots +s_{r+i}^-)+1-1$

  \hspace{5mm} $=\lambda_r+(s_{r+1}^-+ \cdots +s_{r+i}^-)=\lambda_r+(s_{r+1}^-+ \cdots +s_{r+i}^-)+s_{r+i+1}^-$.\vspace{0.3mm}

  \noindent Hence, the results are true for $i+1$.\vspace{2mm}

\noindent {\large{\textbf{54. Periodic properties.}}}\vspace{1mm}

 \noindent We let $f_{2i}=q_{2i+1}+ q_{2i+3}+ q_{2i+5}+\cdots +q_{2i+J}$   for $i \ge 0$.
 If $i \ge 0$, then\vspace{0.5mm}

 \noindent (54.1) \hspace{3mm}$f_{2(i+1)}=q_{2i+3}+ q_{2i+5}+ \cdots +q_{2(i+1)+J}=
 f_{2i}+q_{2(i+1)+J}-q_{2i+1}$.\vspace{1.5mm}

 \noindent {{\textbf{Observation  54.1.}}} Suppose $i \ge 0$.

 \noindent a) $\lambda_{2(i+1)}=\lambda_{2i+2}=\lambda_{2i+1}+s_{2i+2} =\lambda_{2i}-s_{2i+1} +s_{2i+2}$.

\noindent b)  $f_{2(i+1)}=f_{2i}+q_{2(i+1)+J}-q_{2i+1}=f_{2i}-s_{2i+1}+s_{2i+2}$.

\noindent c)  $f_{2(i+1)}-\lambda_{2(i+1)}=f_{2i}-\lambda_{2i}$.\vspace{0.7mm}

 \noindent {{\textbf{Proof.}}} a) follows from (28.1) and (28.2). The first equality in b) follows \mbox{from (54.1).} The second equality is true since (28.3) and (28.4) imply that

 \noindent\hspace{0.1mm} $q_{2(i+1)+J}-q_{2i+1}=q_{2i+2+J}-q_{2i+1}=e_{2i+1}+s_{2i+2}-(s_{2i+1}+e_{2i+1})
 $.\vspace{0.5mm}

\noindent  Moreover, c) follows from a) and b).\vspace{1.5mm}

  \noindent {{\textbf{Observation  54.2.}}} a) $f_{2i}-\lambda_{2i}=f_{0}-\lambda_{0}$ for $i \ge 0$.\vspace{0.2mm}

  \noindent b) Suppose $i>0$ and $q_{2i+m}=q_m$  for $m \ge 1$. Then  $\lambda_{2i}=p+1$.\vspace{0.5mm}

  \noindent {{\textbf{Proof.}}} a) follows from Observation 54.1 c) by induction.\vspace{0.2mm}

  \noindent b) Since $q_{2i+m}=q_m$  for $m \ge 1$, then $f_{2i}=q_{2i+1}+ q_{2i+3}+\cdots +q_{2i+J}$ is equal to $q_{1}+ q_{3}+\cdots +q_{J}=f_0$. Hence, a) implies that  $\lambda_{2i}=\lambda_0=p+1$.\vspace{20mm}

 \hspace{45mm} {\large{\textbf{PART 10.}}  \vspace{ 0.5 mm}}

  \noindent In this part we define and describe properties of progression parameters.
   The main results are (57.2) and (57.3) in Section 57.\vspace{1.5mm}

  \noindent {\large{\textbf{55. Assumptions and definitions.}}}\vspace{0.5mm}

\noindent Suppose  $D^{\,\infty}=(d_1,d_{\hspace{0.2mm}2}, \cdots )$ where\vspace{0.4mm}

\noindent  (55.1)\hspace{16.4mm} $ 0<d_1 \le d_{\hspace{0.2mm}2} \le d_{\hspace{0.2mm}3} \le \cdots$ and $d_j$ $\rightarrow$ $\infty$ if $j \rightarrow \infty$\hspace{0.2mm}.\vspace{0.4mm}

\noindent  Let $r(D^{\infty},\beta)=\#\{i \ge 1: d_i \le \beta\}$ for $\beta > 0$.\vspace{0.6mm}
 \mbox{If $y \ge 1$, then  (55.1) implies that}\vspace{0.4mm}

\noindent (55.2) \hspace{28mm} $r(D^{\infty},\beta)=y$ $\Leftrightarrow$   $d_y \le \beta < d_{y+1}$.\vspace{0.6mm}

\noindent We call $\beta>0$ a progression parameter of $D^{\,\infty}$ if\vspace{0.5mm}

\noindent (55.3) \hspace{17mm} $d_{y+i}=d_i+\beta$ for $i \ge 1$, where $y =r(D^{\infty},\beta)$.\vspace{0.5mm}

\noindent Then $y>0$, since $d_{y+1}=d_1+\beta>d_1$.\vspace{2mm}

\noindent {\large{\textbf{56. Basic properties.}}}\vspace{1mm}

\noindent Suppose $\beta$ and $\beta^*$ are progression parameters of $D^{\infty}$. Let $y =r(D^{\infty},\beta)$ \mbox{and $y^* =r(D^{\infty},\beta^*)$.} Then\vspace{0.5mm}

 \noindent (56.1)  \hspace{18mm} $d_{y+i}=d_i+\beta$ and $d_{y^*+i}=d_i+\beta^*$  for $i \ge 1$.\vspace{0.5mm}

 \noindent We note that $y \ge 1$ and $y^* \ge 1$. Moreover,  by (55.2) we get that\vspace{0.5mm}

 \noindent (56.2)  \hspace{24mm}$d_y \le \beta < d_{y+1}$  and $d_{y^*} \le \beta^* < d_{y^*+1}$.\vspace{2mm}

 \noindent {\textbf{Observation 56.1.}}  $\beta+\beta^*$ is a progression parameter of $D^{\infty}$ and

\hspace{38.7mm} $r(D^{\infty},\beta+\beta^*)=y+y^*$.\vspace{0.5mm}

 \noindent {\textbf{Proof.}} By  (56.1) and (56.2) we get that\vspace{0.5mm}

 \noindent (56.3) \hspace{16mm}$d_{y+y^*+i}=d_{y+i}+\beta^*=d_i+\beta+\beta^*$ for $i \ge 1$,\vspace{0.2mm}

 \noindent (56.4)  \hspace{4mm}$d_{y+y^*}=d_y+\beta^* \le \beta+\beta^*$  and $d_{y+y^*+1}=d_1+\beta+ \beta^*>\beta+ \beta^*$.\vspace{0.5mm}

 \noindent By (55.2) and (56.4) we get that $r(D^{\infty},\beta+\beta^*)=y+y^*$. Hence, (56.3) implies that
$\beta+\beta^*$ is a progression parameter of $D^{\infty}$.\vspace{2mm}

\noindent {\textbf{Observation 56.2.}}   Suppose $\beta^*>\beta$. Then $\beta^*-\beta$ is a progression parameter \mbox{of $D^{\infty}$}
 and $r(D^{\infty},\beta^*-\beta)=y^*-y$.\vspace{0.5mm}

\noindent {\textbf{Proof.}} We note that $y^*= r(D^{\infty},\beta^*)\ge r(D^{\infty},\beta)= y$. Hence, $y^*-\,y+i \ge 1$ for  $i \ge 1$. By (56.1) we get that\vspace{0.2mm}

\noindent (56.5) \hspace{8mm} $d_i+\beta^*=d_{y^*+i}=d_{y+(y^*-y+i)}=d_{y^*-y+i}+\beta$ for $i \ge 1$,\vspace{0.6mm}

 \noindent (56.6) \hspace{28mm}$d_{y^*-\,y+i}=d_i+\beta^*-\beta$ for $i \ge 1$,\vspace{0.6mm}

\noindent (56.7)  \hspace{28mm}$d_{y^*-\,y+1}=d_1+\beta^*-\beta>\beta^*-\beta$.\vspace{0.6mm}

\noindent Since $\beta^*-\beta>0$, then we get by (56.7) that $y^*-y \ge 1$. Then we get  according to (56.1) that $d_{y^*}=d_{y+(y^*-y)}=d_{y^*-y}+\beta$. Hence,\vspace{0.4mm}

\noindent (56.8)  \hspace{32mm} $d_{y^*-\,y}=d_{y^*}-\beta \le \beta^*-\beta$\vspace{0.4mm}

\noindent where the last inequality follows from (56.2). By (55.2), (56.7) and (56.8)
we get that $r(D^{\infty},\beta^*-\beta)=y^*-y$. Hence, (56.6) implies that
$\beta^*-\beta$ is a progression parameter of $D^{\infty}$.\vspace{2mm}

\noindent {\large{\textbf{57. Least progression parameters.}}}\vspace{0.5mm}

\noindent In this section we suppose $\alpha >0$ is a progression parameter of  $D^{\,\infty}$. Then\vspace{0.3mm}

 \noindent (57.1) \hspace{17mm}$d_{\gamma+i}=d_i+\alpha$ for $i \ge 1$, where $\gamma=r(D^{\,\infty},\alpha)$.\vspace{0.3mm}

\noindent We note that $\gamma >0$ and $ 0<d_1 \le d_{\hspace{0.2mm}2} \le \cdots  \le d_{\hspace{0.2mm}\gamma} \le \alpha < d_{\gamma+1} \le \cdots$. Let $D=(d_1, \cdots, d_{\gamma})$. Moreover, let
 $\alpha^*$ and  $\gamma^* $ be the least progression parameters of $D$ with respect to $\alpha$. In the end of this section we will prove:\vspace{0.4mm}

 \noindent (57.2)  \hspace{30.2mm} $r(D^{\,\infty},x\alpha^*)=x\gamma^*$ for $x >0$\hspace{0.2mm},\vspace{0.2mm}

 \noindent   (57.3) \hspace{3.7mm}$\beta$ is a progression parameter of $D^{\,\infty}$  $\Leftrightarrow$ $\beta=x\alpha^*$ where
$x >0$\hspace{0.15mm}.\vspace{0.5mm}

\noindent  Let
 $\beta^*$  be the least    progression parameter  of $D^{\,\infty} $ and $r^*=r(D^{\,\infty},\beta^*)$. Then

  \noindent (57.4) \hspace{35.2mm}$d_{r^*+i}=d_i+\beta^*$ for $i \ge 1$. \vspace{0.5mm}

 \noindent Let $F=\{m >0: m$ factor of $gcd(\alpha,\gamma)\}$ and\vspace{0.5mm}

\noindent (57.5) \hspace{4mm}$E=(e_1, \cdots, e_{2\gamma})$ where $e_{i}=d_i$ and $e_{\gamma+i}=d_i+\alpha$ for $1 \le i \le \gamma$.\vspace{0.5mm}

\noindent By  (57.1) it is easily seen by an inductive argument that\vspace{0.5mm}

 \noindent (57.6) \hspace{23.5mm}$d_{m\gamma+i}=d_i+m\alpha$ for $i \ge 1$ and $m \ge 0$.\vspace{0.5mm}

 \noindent Suppose  $d_{r+i}=d_i+\beta$  for $1 \le i \le \gamma$, where $r>0$ and $\beta>0$\,. Then\vspace{0.8mm}

  \noindent (57.7) \hspace{2.5mm} $d_{\hspace{0.2mm}r+i}=d_i+\beta$  for $i \ge 1$\hspace{0.1mm},
   and $d_{\hspace{0.2mm}xr+i}=d_i+x\beta$  for $i \ge 1$\hspace{0.1mm} and $x \ge 1$.\vspace{0.5mm}

   \noindent  These results are proved as follows. Let $i \ge 1$. Decompose
$i=m\gamma+ j$ such that $m \ge 0$ and $1 \le j\le \gamma$.
   By (57.6) we get that

\hspace{2mm}$d_{\hspace{0.2mm}r+i}=d_{\hspace{0.2mm}r+m\gamma+j}=d_{\hspace{0.2mm}r+j}+m\alpha
=d_j+\beta+m\alpha=d_{m\gamma+j}+\beta=d_i+\beta.$\vspace{0.3mm}

\noindent  Suppose $d_{\hspace{0.2mm}xr+i}=d_i+x\beta$  for $i \ge 1$,\hspace{0.1mm} where $x \ge 1$. Then we get  that\vspace{0.4mm}

\noindent \hspace{2mm} $d_{(x+1)r+i}=d_{r+(xr+i)}=d_{xr+i}+\beta=d_i+x\beta+\beta=d_i+(x+1)\beta$ for $i \ge 1$.\vspace{2mm}

\noindent {\textbf{Observation 57.1.}} a) $\beta=x\beta^*$ is a progression parameter of $D^{\,\infty}$  for $x \ge 1$.

 \noindent b) $\beta$ is a progression parameter of $D^{\,\infty}$  $\Leftrightarrow$ $\beta=x\beta^*$  where $x >0$.

 \noindent c) $r(D^{\,\infty},x\beta^*)=xr^*$ for $x  \ge 1$.\vspace{0.4mm}

 \noindent {\textbf{Proof.}} a) This is true for $x=1$. Suppose $x\beta^*$ is a progression parameter \mbox{of $D^{\,\infty}$}  where $x \ge 1$. By Observation 56.1 we get that $(x+1)\beta^*=x\beta^*+\beta^*$
 is a progression parameter of $D^{\,\infty}$.\vspace{0.2mm}

 \noindent b) Suppose $\beta$ is a progression parameter of $D^{\,\infty}$. By the minimality of $\beta^*$ there exists $x \ge 1$ such that $x\beta^* \le \beta < (x+1)\beta^*$.  Suppose
 $x\beta^* < \beta$. By a) and Observation 56.2 we get that $\beta-x\beta^*$ a progression parameter of $D^{\,\infty}$. Since $\beta-x\beta^* < \beta^*$, this is a contradiction. Hence, $\beta=x\beta^*$. The reverse implication follows from a).\vspace{0.4mm}

 \noindent c) This is true for $x=1$. Suppose $r(D^{\,\infty},x\beta^*)=xr^*$ where $x >0$.
 By a) we get that $\beta=x\beta^*$ is a progression parameter of $D^{\,\infty}$. Since $\beta^*$  is a    progression parameters  of $D^{\,\infty} $ and $r^*=r(D^{\,\infty},\beta^*)$, then Observation 56.1 implies that\vspace{0.3mm}

\hspace{4mm} $r(D^{\,\infty},(x+1)\beta^*)=r(D^{\,\infty},x\beta^*+\beta^*)=xr^*+r^*=(x+1)r^*$.
\vspace{1.6mm}

\noindent {\textbf{Observation 57.2.}}  There exists $m^* \ge 1$ such that $\alpha=m^*\beta^*$ and $\gamma=m^*r^*$.\vspace{0.4mm}

\noindent {\textbf{Proof.}}  Since $\alpha >0$ is a progression parameter of  $D^{\,\infty}$, then
according to  Observation 57.1 b) there exists $m^* \ge 1$ such that $\alpha=m^*\beta^*$. Then we also get from
Observation 57.1 c) that $\gamma=r(D^{\,\infty},\alpha)=r(D^{\,\infty},m^*\beta^*)=m^*r^*$.\vspace{2mm}

\noindent {\textbf{Observation 57.3.}} Suppose $m\beta=\alpha$ and $mr=\gamma$ where $m>0$.

\noindent If $d_{r+i}=d_i+\beta$ for $1 \le i\le \gamma$, then $\beta$ is a progression parameter of $D^{\,\infty}$.\vspace{0.5mm}

\noindent {\textbf{Proof.}} Suppose  $d_{r+i}=d_i+\beta$ for $1 \le i\le \gamma$.  Then (57.7)  implies that\vspace{0.5mm}

\hspace{8mm}$d_{r}+(m-1)\beta =d_{(m-1)r+r}=d_{mr}=d_{\gamma} \le \alpha=m\beta$,\vspace{0.5mm}

\hspace{8mm}$d_{r+1}+(m-1)\beta =d_{(m-1)r+r+1}=d_{mr+1}=d_{\gamma+1} > \alpha=m\beta$.\vspace{0.5mm}

\noindent Hence, $d_r \le \beta < d_{r+1}$. Then (55.2) implies that $r(D^{\infty},\beta)=r$. Moreover, we get by (57.7)  that $d_{r+i}=d_i+\beta$ for $i \ge 1$. Hence, $\beta$ is a progression parameter of $D^{\,\infty}$.\vspace{2mm}

\noindent {\textbf{Observation 57.4.}} $e_i =d_i$ for $1 \le i \le 2\gamma$.\vspace{0.4mm}

\noindent {\textbf{Proof.}} If $1 \le i \le \gamma$, this is trivial. Let $\gamma < i \le 2\gamma$.
Then $i =\gamma+i^*$ \mbox{where $1 \le i^* \le \gamma$,} and (57.1) implies that $e_i=e_{\gamma+i^*}=d_{i^*}+\alpha=d_{\gamma+i^*}=d_i$.\vspace{2mm}

\noindent Next, we let $m^*$ be as in Observation 57.2. Then we get that\vspace{0.5mm}

\noindent (57.8) \hspace{24mm} $\beta^*=\frac{\alpha}{m^*}$ and $r^*=\frac{\gamma}{m^*}$ where $m^*\in F$.\vspace{2mm}

\noindent {\textbf{Observation 57.5.}} $m^*$ is the maximal progression coefficient of $D$ with respect to $\alpha$.\vspace{0.5mm}

\noindent {\textbf{Proof.}}  Observation 57.4 and (57.4) imply that $e_{r^*+i}=d_{r^*+i}=d_i+\beta^*$  \mbox{for $1 \le i \le \gamma$.}   By (57.8) we get that $\beta^*=\frac{\alpha}{m^*}$ and $r^*=\frac{\gamma}{m^*}$ where $m^*\in F$. Hence, $m^*$ is a progression coefficient of $D$ with respect to $\alpha$.\vspace{0.5mm}

 \noindent Suppose $m$ is a progression coefficient of $D$ with respect to $\alpha$. \mbox{Then $m \in F$.}
Let $\beta=\frac{\alpha}{m}$ and $r=\frac{\gamma}{m}$. Then $e_{r+i}=d_i+\beta$ for $1 \le i\le \gamma$.   \mbox{Observation 57.4} implies that $d_{r+i}=d_i+\beta$ for $1 \le i\le \gamma$. Then we get by Observation 57.3  that $\beta$ is a progression parameter of $D^{\,\infty}$. By the minimality property of $\beta^*$ we conclude  that $\frac{\alpha}{m}=\beta \ge \beta^*=\frac{\alpha}{m^*}$.
Hence, $m \le m^*$.\vspace{2mm}

\noindent {\textbf{Proof of (57.2) and (57.3).}} By Observation 57.5 and (57.8) we get \mbox{that
$\beta^*=\frac{\alpha}{m^*}$ and $r^*=\frac{\gamma}{m^*}$} where
$m^*$ is the maximal progression coefficient of $D$ with respect to $\alpha$. As in Section 16 we therefore get that $\beta^*$ and $r^*$ are  the least progression parameters of $D$ with respect to $\alpha$. Hence, $\alpha^*=\beta^*$ \mbox{and $\gamma^*=r^*$.}
Then (57.2) \mbox{and (57.3)} follow from Observation 57.1 b) and c).\vspace{2mm}

\hspace{45mm} {\large{\textbf{PART 11.}}}\vspace{ 1 mm}

\noindent We will derive results about $r$-indexes, $t$-indexes, the contraction vector and the distance vector of  complete shift symmetric vectors. The main result about the indexes and the distance vector are  (65.1), Proposition 66.2 \mbox{and 67.1.}  Moreover, (71.13) is the key result about the contraction vector.\vspace{2mm}

\noindent {\large{\textbf{58.  Assumptions and notation.}}}\vspace{0.5mm}

\noindent  We suppose $p>0$ and $Q^{\,\infty}=C_p^{\infty}(Q)$ is the shift symmetric vector generated by $Q=(q_1, \cdots, q_J,e_0) \in M_p^+$ with respect to $p$. We also suppose $D(Q) \ne \O$.

\noindent Let $Q^{\,\infty}=(q_1,q_2, \cdots )$, and \mbox{let   $( s_1, s_2, \cdots)$,}  $(e_0,e_1,  \cdots)$ and $(\lambda_0, \lambda_1, \cdots)$ be the associated sequences. In particular, $\lambda_0=p+1$.
 By Observation 14.3 \mbox{and 14.4} we get that $Q \in M_p$ and $Q \in M^*$\hspace{0.3mm}.  Moreover, by Observation 28.4 we get that $q_i \ge 1$ for $i \ge 1$. Since $Q \in M^*$, then\vspace{0.5mm}

\noindent (58.1) \hspace{14mm}$J \ge 1$ is odd, $e_0 \ge 0$,  $q_1>1$ and $q_i \ge 1$ for $i \ge 2$\hspace{0.15mm}.\vspace{0.5mm}

\noindent
 Let  $\tau(r)$ be the distance function of $Q^{\,\infty}$.
By  (46.1) we get \mbox{that $Q^{\,\infty}=V(A^{\,\infty})$}  where $A^{\,\infty}=A^{\,\infty}_p(Q)$. Moreover, by (18.1) we  get that $V(A^{\,\infty})$ has an even vector period. Hence,\vspace{0.5mm}

\noindent (58.2) \hspace{27.8mm} $Q^{\,\infty}$ has an even vector period.\vspace{0.2mm}

\noindent Since $p>0$ \mbox{and $q_1>1$\hspace{0.1mm},} then\vspace{0.2mm}

\noindent (58.3) \hspace{23.5mm} $\lambda_0=p+1>1$\hspace{0.1mm}, $s_1>1$ and $q_{J+1}>1$\hspace{0.1mm}\vspace{0.3mm}

\noindent where the last inequalities follow from Observation 51.4. Since $D(Q) \ne \O$,\vspace{0.4mm}

\noindent (58.4)   \hspace{18mm}there exist $i$ such that $q_i=1$
and $1 < i \le J$.

\noindent  Since $q_1>1$  \mbox{and $q_{J+1}>1$} \mbox{where $J+1$} is even, then (58.2) implies that\vspace{0.5mm}

\noindent (58.5) \ $\#\{ i \ge 1: q_i>1$ and  $i$ even\} = $\#\{ i \ge 1: q_i>1$ and  $i$ odd\} = $\infty$\hspace{0.3mm}.\vspace{0.5mm}

\noindent  By (58.2) and (58.4) we get that

\noindent (58.6)  \hspace{34mm} $\#\{ i \ge 1: q_i=1$\} = $\infty$\hspace{0.3mm}.\vspace{0.2mm}

\noindent Hence, by (58.1), (58.5) and (58.6)  we get that $Q^{\,\infty}$ is a complete vector.\vspace{0.5mm}

\noindent Suppose $r \ge 0$.   Let $t_{max}(r)=t$ and $next(r)=r+2t+1$ where $t \ge 0$ is maximal such that   $q_{r+2i}=1$ for $1 \le i \le t$. Moreover, let $\alpha=\delta(Q)+1$.\vspace{2mm}

\noindent {\large{\textbf{59.  The $r$- and $t$-indexes.}}}\vspace{0.5mm}

\noindent We suppose $r_0, r_1, r_2, \cdots$ and  $t_1, t_2, \cdots$ are the $r$- and $t$- indexes of $Q^{\,\infty}$. We choose $I$ maximal such that $r_I \le J$.  By (33.1) we get that

\noindent (59.1) \hspace{23mm} $r_0=0 < r_1 < \cdots < r_I \le J < r_{I+1} < \cdots$

\noindent where $r_{j+1}=next(r_j)$ for $j \ge 0$. Moreover, $t_{j+1}=t_{max}(r_j)$ for $j \ge 0$.\vspace{0.7mm}

\noindent We also let $z_0 \ge 0$ be maximal such \mbox{that $r_I+2z_0 \le J$.} Then $r_I+2z_0=J-1$ or $r_I+2z_0=J$. We note that $r_I+2z_0+1 \le J+1$.\vspace{1.7mm}

\noindent {\textbf{Observation 59.1.}} $r_I+2z_{\hspace{0.3mm}0}=J$, $q_{\hspace{0.1mm}r_I+2i}=1$  for $1 \le i \le z_{\hspace{0.3mm}0}$\hspace{0.2mm}, and $I$ is odd.\vspace{0.7mm}

\noindent {\textbf{Proof.}}
 By (33.2) we get that $r_I+2t_{I+1}+1=r_{I+1}>J$. Hence, $0 \le z_0 \le t_{I+1}$,  and (33.3) implies that $q_{\hspace{0.1mm}r_I+2i}=1$  for $1 \le i \le z_{\hspace{0.3mm}0}$. If  $r_I+2z_{\hspace{0.3mm}0}=J-1$, then we get from  (58.3)   that $q_{r_I+2z_0+2}=q_{J+1}>1$.
 Hence, $t_{I+1}=t_{max}(r_I)=z_0$.
Then  (33.2) implies that   $r_{I+1}=r_I+2t_{I+1}+1=r_I+2z_0+1 \le J$. This is a contradiction since $r_{I+1}> J$. Hence, $r_I+2z_{\hspace{0.3mm}0}=J$.\vspace{0.5mm}

 \noindent Since $r_I+2z_{\hspace{0.3mm}0}=J$ where $J$ is odd, then $r_I$ is odd,  and  (33.4) implies  that $I$ is odd.\vspace{2mm}

\noindent  {\large{\textbf{60. Properties of the $\lambda$\,-\,parameters.}}\vspace{ 1 mm}}

 \noindent If $j \ge 0$, then  (33.3) and Observation 51.1 b) imply that\vspace{0.4mm}

 \noindent (60.1)\hspace{17mm} $q_{r_j+2i}=s_{r_j+2i}=1$ and $e_{r_j+2i}=0$ for $1 \le i \le t_{j+1}$.\vspace{0.5mm}

\noindent  If $s_{r_j+1}>1$ where $j \ge 0$, then Observation 28.4 implies that\vspace{0.4mm}

 \noindent (60.2)\hspace{10.4mm} $s_{r_j+i} \ge 1$ for $i \ge 2$, and  ${s}^-_{r_j+1}+ \cdots +{s}^-_{r_j+i}>0$ for $i \ge 1$.\vspace{0.3mm}

 \noindent Suppose $j \ge 0$. By (33.4) and Observation 28.3 b) we get that\vspace{0.3mm}

 \noindent (60.3)\hspace{16mm} $r_j$ is even if $j$ is even, and $r_j$ is odd if $j$ is odd,\vspace{0.5mm}

  \noindent (60.4)\hspace{8mm} $1 \le \lambda_{r_j} \le p+1$  if $j$ is even, and $0 \le \lambda_{r_j} < p+1$ if $j$ is odd.\vspace{2mm}

\noindent{{\textbf{Observation 60.1.}}} Suppose  $1 \le i \le 2t_{j+1}+1$ where $j \ge 0$ is even.

 \noindent a)  $\lambda_{r_j+i}=\lambda_{r_j}- ({s}^-_{r_j+1}+ \cdots +{s}^-_{r_j+i})$ if $i$ is even.\vspace{0.2mm}

\noindent b)  $\lambda_{r_j+i}=\lambda_{r_j}- ({s}^-_{r_j+1}+ \cdots +{s}^-_{r_j+i})-1$ if $i$ is odd.\vspace{0.8mm}

\noindent{{\textbf{Proof.}}}  a) and b)  follow from (60.1) and Observation 53.2 since $r_j$ is even.\vspace{2mm}

\noindent{{\textbf{Observation 60.2.}}} Suppose  $1 \le i \le 2t_{j+1}+1$ where $j \ge 0$ is odd.

 \noindent a)  $\lambda_{r_j+i}=\lambda_{r_j}+ ({s}^-_{r_j+1}+ \cdots +{s}^-_{r_j+i})$ if $i$ is even.\vspace{0.1mm}

\noindent b)  $\lambda_{r_j+i}=\lambda_{r_j}+({s}^-_{r_j+1}+ \cdots +{s}^-_{r_j+i})+1$ if $i$ is odd.\vspace{0.8mm}

\noindent{{\textbf{Proof.}}}   a) and b) follow from (60.1) and Observation 53.3 since $r_j$ is odd.\vspace{2mm}

\noindent {{\textbf{Observation 60.3.}}} Suppose $j \ge 0$ is even.

\noindent a) $\lambda_{r_{j+1}}=\lambda_{r_j+2t_{j+1}+1}
=\lambda_{r_j}- ({s}^-_{r_j+1}+ \cdots +{s}^-_{r_j+2t_{j+1}+1})-1$.\vspace{0.5mm}

\noindent b) Suppose $\lambda_{r_j}>1$. Then $s_{r_j+1}>1$ and $\lambda_{r_{j+1}} <p$.\vspace{0.8mm}

\noindent{{\textbf{Proof.}}} a) follows from (33.2) and Observation 60.1 b).\vspace{0.5mm}

\noindent b)  By (33.5) and (60.3) we get that $q_{r_j+1}>1$ and $r_j$ is even. Since $\lambda_{r_j}>1$, then   (28.1) implies that
\mbox{ $s_{r_j+1}=min\{q_{r_j+1},\lambda_{r_j}\}>1$.} Then we get according to (60.2) and (60.4) that ${s}^-_{r_j+1}+ \cdots +{s}^-_{r_j+2t_{j+1}+1}>0$ and  $\lambda_{r_j} \le p+1$. Hence, by a) we conclude   that  $\lambda_{r_{j+1}}=\lambda_{r_j}- ({s}^-_{r_j+1}+ \cdots +{s}^-_{r_j+2t_{j+1}+1})-1 <p$.\vspace{200mm}

\noindent {{\textbf{Observation 60.4.}}} Suppose $j \ge 0$ is odd.

\noindent a) $\lambda_{r_{j+1}}=\lambda_{r_j+2t_{j+1}+1}
=\lambda_{r_j}+ ({s}^-_{r_j+1}+ \cdots +{s}^-_{r_j+2t_{j+1}+1})+1$.\vspace{0.3mm}

\noindent b) Suppose $\lambda_{r_j}<p$. Then $s_{r_j+1}>1$ and $\lambda_{r_{j+1}} >1$.\vspace{1mm}

\noindent{{\textbf{Proof.}}} a) follows from (33.2) and Observation 60.2 b).\vspace{0.5mm}

\noindent b)  By (33.5) and (60.3) we get that $q_{r_j+1}>1$ and $r_j$ is odd. Since $\lambda_{r_j}<p$, then  (28.2) implies that
\mbox{ $s_{r_j+1}=min\{q_{r_j+1},p+1-\lambda_{r_j}\}>1$.} Then we get from (60.2) and (60.4) that ${s}^-_{r_j+1}+ \cdots +{s}^-_{r_j+2t_{j+1}+1}>0$ and  $\lambda_{r_j} \ge 0$. Hence, by a) we conclude  that  $\lambda_{r_{j+1}}=\lambda_{r_j}+ ({s}^-_{r_j+1}+ \cdots +{s}^-_{r_j+2t_{j+1}+1})+1 >1$.\vspace{1mm}

\noindent {{\textbf{Observation 60.5.}}} Suppose $j \ge 0$. Then\vspace{0.3mm}

\hspace{24mm} $\lambda_{r_j}>1$ if $j$ is even, and  $\lambda_{r_j} < p$ \mbox{if $j$} is odd.\vspace{0.8mm}

\noindent {{\textbf{Proof.}}}  By (58.3) we get that $\lambda_0>1$.
Since $r_0=0$,
the result is true \mbox{for $j=0$.} Next, suppose the result is true for $j \ge 0$. If $j$ is even, then $\lambda_{r_j}>1$. Hence, by Observation 60.3 b) we get that
  $\lambda_{r_{j+1}}<p$.
 If $j$ is odd, then $\lambda_{r_j}<p$. Hence, by Observation 60.4 b) we get that $\lambda_{r_{j+1}}>1$. That means, the result is true for $j+1$.\vspace{1mm}

 \noindent {{\textbf{Observation 60.6.}}} $s_{r_j+1}>1$ for $j \ge 0$.\vspace{0.5mm}

 \noindent {{\textbf{Proof.}}} If $j \ge 0$ is even,  then  Observation  60.5 and 60.3 b) imply that $\lambda_{r_j}>1$ and $s_{r_j+1}>1$.  If $j \ge 0$ is odd,  then  we get from Observation  60.5 \mbox{and 60.4 b)}  that $\lambda_{r_j}<p$ and $s_{r_j+1}>1$.\vspace{2mm}

\noindent {{\textbf{Observation 60.7.}}} Suppose $j \ge 0$. Then
$q_{r_j+J+2t_{j+1}+2}=q_{r_{j+1}+J+1}>1$.\vspace{1mm}

\noindent {{\textbf{Proof.}}}  By Observation 60.6  we get that $s_{r_{j+1}+1}>1$. Hence,
(33.2) and Observation 51.3
imply that $q_{r_j+J+2t_{j+1}+2}=q_{r_{j+1}+J+1}>1$.\vspace{1mm}

\noindent {{\textbf{Observation 60.8.}}} Suppose $j \ge 0$ is even and $1 \le i \le 2t_{j+1}$ where $i$ is even. Then ${s}^-_{r_j+1}+ \cdots +{s}^-_{r_j+i}>0$ and
 $\lambda_{r_j+i}=\lambda_{r_j}-({s}^-_{r_j+1}+ \cdots +{s}^-_{r_j+i})\le p$.\vspace{1mm}

\noindent {{\textbf{Proof.}}} By  (60.4) and Observation 60.6  we get that $\lambda_{r_j} \le p+1$ and $s_{r_j+1}>1$. Hence, (60.2) and Observation 60.1 a) imply that the results are true.\vspace{1.5mm}

\noindent {\large{\textbf{61. Additional parameters and results.}}}\vspace{0.6mm}

\noindent   For $j \ge 0$ we let  $m_{j+1}$ be the maximal integer such that\vspace{0.5mm}

\noindent (61.1) \hspace{16.8mm}$0 \le m_{j+1} \le t_{j+1}$ and $e_{r_j+i}=0$ for $1 \le i \le 2m_{j+1}$\hspace{0.2mm}.\vspace{0.4mm}

\noindent We also let $z_{j+1}=t_{j+1}-m_{j+1}$  for $j \ge 0$.  If $j \ge 0$, then\vspace{0.4mm}

\noindent (61.2) \hspace{14mm} $t_{j+1}=m_{j+1}+z_{j+1}$ where $m_{j+1} \ge 0$ and $z_{j+1} \ge 0$,\vspace{0.3mm}

\noindent (61.3) \hspace{17mm}$r_j+2m_{j+1}+2z_{j+1}+1=r_j+2t_{j+1}+1=r_{j+1}$,

\noindent (61.4) \hspace{30.7mm}$q_{r_j+i}=s_{r_j+i}$ for $1 \le i \le 2m_{j+1}$,\vspace{0.2mm}

\noindent where   (61.2) is trivial, (61.3) follows from (33.2) and (61.2), and (61.4) follows from (28.3) and (61.1).\vspace{100mm}

\noindent {\textbf{Observation 61.1.}}  If $e_{r_j+2m_{j+1}+1}=0$ where $j \ge 0$, then $m_{j+1}=t_{j+1}$.\vspace{0.6mm}

\noindent {\textbf{Proof.}} Suppose  $e_{r_j+2m_{j+1}+1}=0$.
 If $m_{j+1} <t_{j+1}$, then (60.1) implies that \mbox{$e_{r_j+2m_{j+1}+2}=e_{r_j+2(m_{j+1}+1)}=0$.}
  Hence,  $e_{r_j+i}=0$ for $1 \le i \le 2(m_{j+1}+1)$. \mbox{By (61.1)} this is a contradiction
  and we get that $m_{j+1}=t_{j+1}$.\vspace{1.5mm}

\noindent {\large{\textbf{62. The even case.}}}\vspace{0.6mm}

\noindent {{\textbf{Observation  62.1.}}} Suppose
 $e_{r_j+2m_{j+1}+1}>0$ where $j \ge 0$ is even. Then\vspace{0.4mm}

 \noindent \hspace{4mm}  $\lambda_{r_j+2i+1}=0$  for $m_{j+1} \le i \le t_{j+1}$, \mbox{and $s_{r_j+2i+1}=1$  for $m_{j+1} < i \le t_{j+1}$.}\vspace{0.7mm}

\noindent {{\textbf{Proof.}}} By (60.3)  we get that $r_j$ and $r_j+2m_{j+1}$ are  even. Hence, (28.3) and Observation 51.2 a)  imply that
 $s_{r_j+2m_{j+1}+1}< q_{r_j+2m_{j+1}+1}$ \mbox{and  $\lambda_{r_j+2m_{j+1}+1}=0$.}
 Next, we suppose $\lambda_{r_j+2i+1} =0$ where $m_{j+1} \le i <t_{j+1}$.
By (60.1) we get \mbox{that $s_{r_j+2i+2}=s_{r_j+2(i+1)}=1$.} Since $r_j+2i$ is even, then Observation 51.5 implies that $s_{r_j+2(i+1)+1}=s_{r_j+2i+3}=1$ and $\lambda_{r_j+2(i+1)+1}=\lambda_{r_j+2i+3}=0$.\vspace{2mm}

\noindent {\large{\textbf{63. The odd case.}}}\vspace{0.6mm}

\noindent {{\textbf{Observation  63.1.}}} Suppose
 $e_{r_j+2m_{j+1}+1}>0$ where $j \ge 0$ is odd. Then

 \noindent  \hspace{2mm} $\lambda_{r_j+2i+1}=p+1$  for $m_{j+1} \le i \le t_{j+1}$, \mbox{and $s_{r_j+2i+1}=1$  for $m_{j+1} < i \le t_{j+1}$.}\vspace{0.7mm}

\noindent {{\textbf{Proof.}}} By (60.3)  we get that $r_j$ and $r_j+2m_{j+1}$ are odd. Then (28.3) and Observation 51.2 b)  imply that
 $s_{r_j+2m_{j+1}+1}< q_{r_j+2m_{j+1}+1}$ and  $\lambda_{r_j+2m_{j+1}+1}=p+1$. Suppose $\lambda_{r_j+2i+1} =p+1$ where $m_{j+1} \le i <t_{j+1}$.
\mbox{By (60.1)}  we get that $s_{r_j+2i+2}=s_{r_j+2(i+1)}=1$. Since $r_j+2i$ is odd, then  Observation 51.6 implies that
$s_{r_j+2(i+1)+1}=s_{r_j+2i+3}=1$ and $\lambda_{r_j+2(i+1)+1}=\lambda_{r_j+2i+3}=p+1$.\vspace{2mm}

\noindent {\large{\textbf{64. Properties of the parameters.}}}\vspace{1mm}

\noindent {{\textbf{Observation 64.1.}}} Suppose  $m_{j+1}<t_{j+1}$ where $j \ge 0$. Then   $s_{r_j+2i+1}=1$ for $m_{j+1} < i \le t_{j+1}$, and
$s_{r_j+i}=1$ for $2m_{j+1}+1 < i \le 2t_{j+1}+1$.\vspace{0.6mm}

\noindent {{\textbf{Proof.}}}  Observation 61.1 implies that $e_{r_j+2m_{j+1}+1}>0$.
Hence, the results  follow from (60.1), Observation 62.1 and 63.1.\vspace{2mm}

\noindent {{\textbf{Proposition 64.2.}}} Suppose $j \ge 0$.

\noindent a) $q_{r_j+J+2i}=1$ for $1 \le i \le m_{j+1}$\hspace{0.2mm}.

\noindent b)  $q_{\hspace{0.2mm}r_j+J+2i+1}=1$ for $m_{j+1} < i \le t_{j+1}$\hspace{0.2mm}.

\noindent c)  $q_{r_j+J+2m_{j+1}+2}>1$.\vspace{0.7mm}

\noindent {{\textbf{Proof.}}}   a) By using (60.1) and (61.1) we get that  $e_{\hspace{0.2mm}r_j+2i-1}=0$ and $s_{\hspace{0.2mm}r_j+2i}=1$ for $1 \le i \le m_{j+1}$\hspace{0.2mm}. Then (28.4) implies that

\hspace{20mm}\mbox{$q_{\hspace{0.2mm}r_j+J+2i}=e_{\hspace{0.2mm}r_j+2i-1}+s_{r_j+2i}=1$} \mbox{for $1 \le i \le m_{j+1}$\hspace{0.2mm}.}\vspace{0.6mm}

\noindent b) If $m_{j+1}=t_{j+1}$\hspace{0.2mm}, this is trivial. Otherwise,   (60.1) and   Observation 64.1   imply that $e_{r_j+2i}=0$ and $s_{r_j+2i+1}=1$ for $m_{j+1}<i \le t_{j+1}$\hspace{0.2mm}. Hence, we get from (28.4)   that
 $q_{\hspace{0.2mm}r_{j}+J+2i+1}=e_{\hspace{0.2mm}r_j+2i}+s_{\hspace{0.2mm}r_j+2i+1}=1$ for $m_{j+1}<i \le t_{j+1}$.\vspace{0.6mm}

\noindent c)   First, we suppose $m_{j+1}= t_{j+1}$\hspace{0.2mm}. Then we get by using Observation 60.7 that $q_{r_j+J+2m_{j+1}+2}=q_{r_j+J+2t_{j+1}+2}>1$. Next, we suppose $m_{j+1}< t_{j+1}$. Then  Observation 61.1
implies that $e_{r_j+2m_{j+1}+1}>0$. Hence, we get from Observation 51.3 \mbox{that  $q_{\hspace{0.2mm}r_j+J+2m_{j+1}+2}>1$.}\vspace{2mm}

 \noindent {\large{\textbf{65. Crucial results.}}}\vspace{0.6mm}

\noindent We  prove in the end of this section that\vspace{0.4mm}

\noindent (65.1) \hspace{8mm}$r_{I+j}=r_{j-1}+J+2m_{j}+1$ and $t_{I+j}=z_{j-1}+m_{j}$ for  $j \ge 1$.\vspace{1.4mm}

\noindent {{\textbf{Observation 65.1.}}} Suppose $j \ge 0$, $t \ge 0$, $q_{r_j+2i}=1$ for $1 \le i \le t$, \mbox{and $q_{r_j+2t+2}>1$.} Then  $t_{j+1}=t$ and $r_{j+1}=r+2t+1$.\vspace{0.5mm}

\noindent {{\textbf{Proof.}}} We note that $t=t_{max}(r_j)=t_{j+1}$.   Hence, by  (33.2) we also get that\vspace{0.3mm}

\hspace{30mm}$r_{j+1}=r_j+2t_{j+1}+1=r+2t+1$.\vspace{2mm}

\noindent {{\textbf{Observation 65.2.}}} a) $q_{r_I+2i}=1$ for $1 \le i \le z_0+m_1$.

\noindent b)
$q_{r_I+2(z_0+m_1)+2}=q_{r_0+J+2m_1+2}>1$.

\noindent c) $t_{I+1}=z_0+m_1$   and $r_{I+1}=r_I+2(z_0+m_1)+1=r_0+J+2m_1+1$.\vspace{0.7mm}

\noindent {{\textbf{Proof.}}} Since $r_0=0$, then by Observation 59.1 we get
that $r_I+2z_0=r_0+J$.\vspace{0.4mm}

\noindent a) According to Observation 59.1 we get that a) is true for $1 \le i \le z_0$.  Suppose $z_0 < i \le z_0+m_1$. Then $1 \le i-z_0 \le m_1$ and  Proposition 64.2 a) implies that
$q_{r_I+2i}=q_{r_I+2z_0+2(i-z_0)}=q_{r_0+J+2(i-z_0)}=1$.\vspace{0.6mm}

\noindent b) follows from Proposition 64.2 c) and the equality $r_I+2z_0=r_0+J$.\vspace{0.5mm}

\noindent    c) follows
from a), b), Observation 65.1 and the equality $r_I+2z_0=r_0+J$.\vspace{1mm}

\noindent {{\textbf{Observation 65.3.}}} Suppose $r_{I+j}=r_{j-1}+J+2m_j+1$ where $j \ge 1$.

\noindent a) $r_{I+j}+2z_j=r_{j-1}+J+2m_j+1+2z_j=r_j+J$.

\noindent b) $q_{r_{I+j}+2i}=1$ for $1 \le i \le z_j+m_{j+1}$.

\noindent c)
$q_{r_{I+j}+2(z_j+m_{j+1})+2}=q_{r_{j}+J+2m_{j+1}+2} >1$.\vspace{0.5mm}

\noindent d) $t_{I+j+1}=z_j+m_{j+1}$  and $r_{I+j+1}=r_{I+j}+2(z_j+m_{j+1})+1=r_{j}+J+2m_{j+1}+1$.\vspace{0.7mm}

\noindent {{\textbf{Proof.}}} a) follows from (61.3).

\noindent b) Suppose $1  \le i \le z_j$. By (61.2) and Proposition 64.2 b) we get that\vspace{0.1mm}

\hspace{4mm}$m_j < m_j+i \le m_j+z_j=t_{j}$ and $q_{r_{I+j}+2i}=q_{r_{j-1}+J+2(m_j+i)+1}=1$.\vspace{0.5mm}

\noindent Suppose $z_j  < i \le z_j+m_{j+1}$. Then $1 \le i-z_j \le m_{j+1}$. By using a) and Proposition 64.2 a)
 we get that $q_{r_{I+j}+2i}=q_{r_{I+j}+2z_j+2(i-z_j)}=q_{r_j+J+2(i-z_j)}=1$.\vspace{0.5mm}

\noindent c) follows from a) and Proposition 64.2 c).\vspace{0.4mm}

\noindent d) follows from a), b), c) and Observation 65.1.\vspace{1mm}

\noindent {{\textbf{Proof of (65.1).}}} By Observation 65.2 c) we get that (65.1) is true for $j=1$.
Suppose $r_{I+j}=r_{j-1}+J+2m_j+1$ where $j \ge 1$. Then we get according to Observation 65.3 d) that (65.1) is true for $j+1$.\vspace{2mm}

\noindent {\large{\textbf{66. Properties of the even vector periods.}}}\vspace{1mm}

\noindent {\textbf{Observation 66.1.}} Suppose $r \ge  0$ is even and $r \ne r_j$ for each $j \ge 0$\hspace{0.35mm}. Then\vspace{0.1mm}

\hspace{41mm}$q_{r+1}=1$ or $\lambda_{\hspace{0.25mm}r}\le p$\hspace{0.1mm}.\vspace{0.7mm}

\noindent {\textbf{Proof.}}  By (33.1) there exists $j \ge 0$  such that $r_{j}<r <r_{j+1}$\hspace{0.1mm}. By (33.2) we get that $r_{j+1}=r_j+2t_{j+1}+1$\hspace{0.1mm}.   \mbox{There exists $i$} such that $r=r_j+2i-1$ \mbox{or $r=r_j+2i$} where $1 \le i \le t_{j+1}$\hspace{0.15mm}. If $r=r_j+2i-1$\hspace{0.1mm}, then we get according \mbox{to (33.3)}  that $q_{\hspace{0.2mm}r+1}=q_{r_j+2i}=1$\hspace{0.1mm}. If $r=r_j+2\hspace{0.15mm}i$\hspace{0.1mm}, then $r_j$ is even since $r$ is even. Hence, (60.3) and Observation 60.8 imply  that $j$ is even and
 $\lambda_{\hspace{0.25mm}r}=\lambda_{r_j+2i} \le p$.\vspace{1mm}

\noindent {\textbf{Proposition 66.2.}}  Suppose $r$ is an even vector period of
$Q^{\,\infty}$. Then there exists  $j >0$ such that $r=r_j$\hspace{0.24mm}.\vspace{0.5mm}

\noindent {\textbf{Proof.}} By (58.1) we note that $q_{r+1}=q_1>1$. We choose $i \ge 0$ such
\mbox{that $r=2i$.} Then $q_{2i+m}=q_m$ for $m \ge 1$, and we get from \mbox{Observation 54.2 b)}   that $\lambda_r=\lambda_{2i}=p+1$. Hence, the result follows from Observation 66.1.\vspace{2mm}

\noindent {\large{\textbf{67.  Properties of the distance vector.}}}\vspace{0.5mm}

\noindent     Let   $c_{\hspace{0.45mm}0}=0 <c_{\hspace{0.25mm}1}< c_{\hspace{0.25mm}2}<\cdots $  be the $c$\,-\,indexes of $Q^{\,\infty}=(q_1, q_2, \cdots )$, and
let $\tau$ be the distance function of $Q^{\,\infty}$.
 The distance vector of $Q^{\,\infty}$ is given by\vspace{0.4mm}

\noindent (67.1) \hspace{16.5mm} $D^{\,\infty}=(d_1,d_2,  \cdots)$ where $d_{\hspace{0.35mm}i}=\tau(c_{\hspace{0.45mm}i})$ for $i \ge 1$\hspace{0.1mm}.\vspace{0.5mm}

\noindent Since $c_1$ is the least index $>c_0+1=1$ such that $q_{c_1}=1$, then (58.4) implies that $c_1 \le J$. Let $\alpha=\delta(Q)+1$ \mbox{and $\gamma=r(D^{\,\infty},\alpha)$.} By (32.10) and (52.1) we get that $d_1=\tau(c_1) \le \tau(J) \le \alpha$. Hence, $\gamma=r(D^{\,\infty},\alpha) \ge 1$. Moreover,\vspace{0.5mm}

\noindent (67.2)\hspace{12mm} $D(Q)=(d_1, \cdots, d_{\gamma})$ and
 $d_{\gamma+i}=d_i+\alpha$ for $i \ge 1$,\vspace{0.5mm}

 \noindent where (67.2) follow from  Proposition 69.1 and 70.6. By (35.4) we get that\vspace{0.5mm}

 \noindent (67.3) \hspace{20.4mm}  $0 <d_1 \le d_2 \le \cdots$ and  $d_i \rightarrow \infty$ if $\rightarrow \infty$\,.\vspace{2mm}

\noindent {\textbf{Proposition 67.1.}} Suppose $\alpha^*$  and $\gamma^*$ are  the least progression parameters  of $D(Q)$ with respect to $\alpha$. Then $r(D^{\,\infty},x\alpha^*)=x\gamma^*$ for $x >0$\hspace{0.2mm}, and \vspace{0.05mm}

\noindent \hspace{13mm} $\beta$ is a progression parameter of $D^{\,\infty}$ $\Leftrightarrow$ $\beta=x\alpha^*$ where $x >0$\hspace{0.35mm}.\vspace{0.7mm}

\noindent {\textbf{Proof.}} By (67.2) we get that $d_{\gamma+i}=d_{i}+\alpha$  for $ i \ge 1$. Since $\gamma=r(D^{\,\infty},\alpha)$, then $\alpha >0$ is a progression parameter of  $D^{\,\infty}$. Let  $D=(d_1, \cdots, d_{\gamma})$.
\mbox{By (67.2)} we get that
$D(Q)=D$. Hence,   $\alpha^*$  and $\gamma^*$ are  the least progression parameters  of $D$ with respect to $\alpha$. The results follow from (57.2) and (57.3).\vspace{2mm}

\noindent {\large{\textbf{68. Notation and auxiliary results.}}}

\noindent   As in Section 37 we suppose $ y_j =y(r_j)$ for $j \ge 0$. According to  (35.8), and Observation 37.1 b) we get that\vspace{0.4mm}

 \noindent (68.1)  \hspace{16.7mm}$y_{j+1}=y_j+t_{j+1}$ and $c_{y_j} \le r_j < c_{y_j+1}$   for $j \ge 0$.\vspace{10mm}

 \noindent According to   Observation 37.1 a), 37.2 and 37.4 a) we get that\vspace{0.4mm}

 \noindent (68.2) \hspace{2mm}$c_{\hspace{0.35mm}y_j+\hspace{0.35mm}i}=r_j+2i$ and $d_{\hspace{0.35mm}y_j+\hspace{0.35mm}i}=\tau(r_j+2i)$ for $1 \le i \le t_{j+1}$ and $j \ge 0$,\vspace{0.4mm}

\noindent (68.3) \hspace{16.7mm}$y_0=0 \le y_1 \le y_2 \le \cdots$  and $y_j \rightarrow \infty$ if $j \rightarrow \infty$.\vspace{0.5mm}

\noindent  By (28.9) and (58.3)  we get that $s_1>1$ and
 $e_0 < e_0+s_1^-=q_{J+1}^-$. Hence, we get by  (52.1)  that\vspace{0.5mm}

\noindent (68.4) \hspace{15.4mm}$\tau(J)\le \alpha=\tau(J)+e_0 <\tau(J)+q_{J+1}^-=\tau(J+1)$.\vspace{2mm}

 \noindent {\textbf{Observation  68.1.}}  $c_{y_I+z_0} \le J < c_{y_I+z_0+1}$ where $y_I+z_0 >0$.\vspace{0.4mm}

  \noindent {\textbf{Proof.}}  By Observation 59.1  we get that $J=r_I+2z_0$. If $z_0=0$, \mbox{then $r_I=J$.} By (68.1) we get that  $c_{y_I} \le r_I < c_{y_I+1}$. Hence, $c_{y_I+z_0} \le J < c_{y_I+z_0+1}$. Suppose $z_0>0$.   Then  (65.1) and (68.2) imply that $1 \le z_0 \le z_0+m_1=t_{I+1}$ and
 $c_{y_I+z_0}=r_I+2z_0=J$. Hence, $c_{y_I+z_0} \le J < c_{y_I+z_0+1}$ also in this case.\vspace{0.5mm}

\noindent  By Section 67 we get that $c_1 \le J$. Since $ J<c_{y_I+z_0+1}$,
  then $y_I+z_0 >0$.\vspace{2mm}

  \noindent {\textbf{Observation  68.2.}}    $\gamma=y_I+z_0$ and $c_{\gamma} \le J < c_{\gamma+1}$.\vspace{0.4mm}

\noindent {\textbf{Proof.}}  By Observation 68.1 we get that
 $c_{y_I+z_0} \le J < c_{y_I+z_0+1}$ \mbox{and $y_I+z_0 >0$.} Then (32.10), (67.1) and (68.4)
  imply that\vspace{0.4mm}

 \hspace{2mm}$d_{y_I+z_0}=\tau(c_{y_I+z_0}) \le \tau(J) \le \alpha < \tau(J+1) \le \tau(c_{y_I+z_0+1})=d_{y_I+z_0+1}$.\vspace{0.6mm}

\noindent Hence, we get from (35.5) that  $\gamma=r(D^{\,\infty},\alpha)=y_I+z_0$
and  $c_{\gamma} \le J < c_{\gamma+1}$.\vspace{2mm}

\noindent {\large{\textbf{69. The distance vector of $Q$.}}}\vspace{0.8mm}

\noindent Let $\tau_Q$ be the distance vector of  $Q=(q_1, \cdots, q_J,e_0)$.
 Since $Q^{\,\infty}=(q_1, q_2, \cdots )$, then  $\tau_Q(c)=\tau(c)$  for $0 \le c \le J$. By Observation 68.2 we get that

 \noindent (69.1)    \hspace{10mm} $c_{\hspace{0.45mm}0}=0 < c_1 < c_2 < \cdots <c_{\hspace{0.35mm}\gamma}\le J < c_{\hspace{0.35mm}\gamma+1} < c_{\hspace{0.35mm}\gamma+2} < \cdots$\,.\vspace{2mm}

\noindent {{\textbf{Proposition  69.1.}}}
 $D(Q)=(d_1, \cdots d_{\hspace{0.25mm}\gamma})$\hspace{0.15mm}.\vspace{0.5mm}

 \noindent {{\textbf{Proof.}}} a) By (35.2) and (69.1) we get that
 $c_{\hspace{0.45mm}0}=0 < c_1 < c_2 < \cdots <c_{\hspace{0.25mm}\gamma}\le J$,

 \noindent  \hspace{9mm} $c_{i+1}$ is the least index $>c_i+1$ such that $q_{c_{i+1}}=1$
 for $0 \le i<\gamma$,\vspace{0.2mm}

\hspace{37.7mm} $q_c>1$ if $c_{\hspace{0.35mm}\gamma}+1 < c \le J$.\vspace{0.3mm}

\noindent Then we get by (15.1), (15.2) and (15.3)  that $c_{\hspace{0.45mm}0}\hspace{0.35mm}, c_1 , c_2, \cdots ,c_{\hspace{0.35mm}\gamma}$ are the $c$\,-\,indexes \mbox{of $Q=(q_1, \cdots, q_J,e_0)$\hspace{0.15mm}.} Hence, we get from (15.5) and (67.1) that\vspace{0.5mm}

\hspace{3.6mm}$D(Q)=(\tau_Q(c_1), \cdots, \tau_Q(c_{\gamma}))=(\tau(c_1), \cdots, \tau(c_{\gamma}))=(d_1, \cdots, d_{\gamma})$.\vspace{2mm}

\noindent {\large{\textbf{70.  Arithmetical  properties.}}}\vspace{0.5mm}

\noindent {{\textbf{Observation 70.1.}}}  $\gamma+y_j=y_{\hspace{0.15mm}I+j}+z_j$ for $j \ge 0$\hspace{0.35mm}.\vspace{0.5mm}

\noindent {\textbf{Proof.}}
Since $y_0=0$\hspace{0.15mm}, we get from Observation 68.2 that $\gamma+y_0=y_{\hspace{0.15mm}I+0}+z_0$\hspace{0.15mm}.\vspace{0.4mm}

  \noindent  If $\gamma+y_j=y_{\hspace{0.15mm}I+j}+z_j$ where $j \ge 0$\hspace{0.2mm},  then (61.2), (65.1) and (68.1)  imply that\vspace{0.4mm}

 \hspace{15.9mm}$\gamma+y_{j+1}=\gamma+y_j+t_{j+1}=y_{I+j}+z_j+m_{j+1}+z_{j+1}$

 \hspace{32mm}$=y_{I+j}+t_{I+j+1}+z_{j+1}=y_{I+j+1}+z_{j+1}$\hspace{0.2mm}.\vspace{200mm}

 \noindent {{\textbf{Observation 70.2.}}} Suppose $j \ge 0$.

 \noindent a) $\tau(r_j+J+2i)= \tau(r_j+2i)+\alpha$ for $1 \le i \le t_{j+1}$.

 \noindent b) $\tau(r_j+J+2i+1)= \tau(r_j+J+2i)=\tau(r_j+2i)+\alpha$ for $m_{j+1} < i \le t_{j+1}$.\vspace{0.5mm}

 \noindent {{\textbf{Proof.}}} a) Suppose $1 \le i \le t_{j+1}$.
 By (33.3) we get that $q_{r_j+2i}=1$. Hence, a) follows from Observation 52.2.\vspace{0.2mm}

\noindent b) Suppose $m_{j+1} < i \le t_{j+1}$.  We get by Proposition 64.2 b) \mbox{that $q_{r_j+J+2i+1}=1$.} Hence, the first equality is true. The second equality follows from a).\vspace{2mm}

\noindent {{\textbf{Observation 70.3.}}} a) $r_{I+j}+2z_j=r_j+J$  for $j \ge 0$.

\noindent b) $r_{I+j+1}-2m_{j+1}=r_j+J+2m_{j+1}+1-2m_{j+1}=r_{j}+J+1$ for $j \ge 0$.

\noindent c) $y_{I+j}+z_j+m_{j+1}=y_{I+j}+t_{I+j+1}=y_{I+j+1}$
for $j \ge 0$.\vspace{0.5mm}

\noindent {\textbf{Proof.}} a) Since $r_0=0$, then  Observation 59.1 implies that $r_{I}+2z_0=r_0+J$.
If $j \ge 1$, then we get by (65.1) that $r_{I+j}=r_{j-1}+J+2m_j+1$. Hence, Observation 65.3 a) implies that  $r_{I+j}+2z_j=r_j+J$.\vspace{0.5mm}

\noindent b) and c) follow from  from (65.1) and (68.1).\vspace{2mm}

\noindent {{\textbf{Observation 70.4.}}} Let $i=y_j+i^*$ where $1 \le i^* \le m_{j+1}$ and $j \ge 0$.
Then\vspace{0.5mm}

\hspace{1.5mm}$d_{\gamma+i}=d_{\gamma+y_j+i^*}=d_{y_{I+j}+z_j+i^*}=
\tau(r_{I+j}+2(z_j+i^*))=\tau(r_j+J+2i^*)$\vspace{0.5mm}

\hspace{10mm}$=\tau(r_j+2i^*)+\alpha=d_{y_j+i^*}+\alpha=d_i+\alpha$.\vspace{0.6mm}

\noindent {\textbf{Proof.}}  The first equality is trivial. The second equality follows from Observation 70.1.
By  (65.1) we get that $1 \le  z_j+i^* \le z_j+m_{j+1}=t_{I+j+1}$. Hence,  the third  follows from (68.2). By (61.1) we get that $1 \le i^* \le m_{j+1}\le t_{j+1}$. The next equalities follow from Observation 70.3 a), Observation 70.2 a) \mbox{and (68.2).}\vspace{2mm}

 \noindent {{\textbf{Observation 70.5.}}} Let $i=y_j+i^*$ where $m_{j+1} < i^* \le t_{j+1}$ and $j \ge 0$.
Then

\noindent\hspace{3mm}$d_{\gamma+i}=d_{\gamma+y_j+i^*}=d_{y_{I+j}+z_j+i^*}=d_{y_{I+j+1}+i^*-m_{j+1}}=
\tau(r_{I+j+1}+2(i^*-m_{j+1}))$\vspace{0.3mm}

\noindent \hspace{11.2mm}$=\tau(r_j+J+2i^*+1)=\tau(r_j+2i^*)+\alpha=d_{y_j+i^*}+\alpha=d_i+\alpha$.\vspace{0.6mm}

\noindent {\textbf{Proof.}} The first equality is trivial, the second follows from Observation 70.1 and
the third  follows from Observation 70.3 c). By (61.2) and (65.1) we get that $1 \le i^*-m_{j+1} \le t_{j+1}-m_{j+1}=z_{j+1}  \le z_{j+1}+m_{j+2}=t_{I+j+2}$, and the forth  follows from (68.2).
We note that $0 \le m_{j+1} < i^* \le t_{j+1}$. The next equalities follow from Observation 70.3 b), Observation 70.2 b) and (68.2).\vspace{2mm}

\noindent {{\textbf{Proposition 70.6.}}} $d_{\hspace{0.25mm}\gamma+\hspace{0.4mm}i}=d_{\hspace{0.25mm}i}+\alpha$ for $i \ge 1$\hspace{0.15mm}.\vspace{0.5mm}

\noindent {\textbf{Proof.}} Let $i \ge 1$. By (68.1) and (68.3) there exists $j \ge 0$  such that $i=y_j+i^*$ where $1 \le i^* \le t_{j+1}$. By (61.1) we get that $0 \le m_{j+1} \le t_{j+1}$.\vspace{0.5mm}

\noindent If $1 \le i^* \le m_{j+1}$, then the result follows by Observation 70.4.\vspace{0.5mm}

\noindent If $m_{j+1} < i^* \le t_{j+1}$, then the result follows from Observation 70.5.\vspace{2mm}

\newpage

\noindent {\large{\textbf{71.  Properties of the contraction vectors.}}}\vspace{1mm}

 \noindent Let $Q^*=\pi(Q)$ and $Q^{\,\infty}_*=\pi(Q^{\,\infty})$ be the contraction vectors of $Q$ and $Q^{\,\infty}$. By (33.6)  we get that\vspace{0.4mm}

\noindent (71.1) \hspace{10mm}$Q^{\,\infty}_*=(q_1^*,q_2^*, \cdots )$ where $q_{j+1}^*=\tau(r_{j+1})-\tau(r_j)$  for $j \ge 0$.\vspace{0.4mm}

 \noindent By (33.1) we get that
 $r_0=0 < r_1 < r_2 <\cdots $.  We decompose\vspace{0.6mm}

\noindent (71.2)\hspace{5.5mm} $Q^{\,\infty}=(G_1,G_2, \cdots )$ where $G_{j+1}=(q_{r_j+1}, \cdots , q_{r_{j+1}})$ for $j \ge 0$.\vspace{0.6mm}

\noindent Then we get by (32.8) and (71.1) that\vspace{0.5mm}

\noindent (71.3) \hspace{2.5mm}$\delta(G_{j+1})=q_{r_j+1}^- + \cdots + q_{r_{j+1}}^-=\tau(r_{j+1})-\tau(r_j)=q_{j+1}^*$ for $j \ge 0$,\vspace{0.7mm}

\noindent (71.4)\hspace{15mm}  $\delta(G_1, \cdots, G_j)=\delta(q_1, \cdots, q_{r_j})=\tau(r_j)$  for $j \ge 1$.\vspace{0.3mm}

\noindent  Next, we let

\noindent  (71.5)\hspace{12.5mm} $p^*=p-1$ and $e_0^*=\alpha-\tau(r_I)$ where $\alpha=\delta(Q)+1$,\vspace{0.5mm}

\noindent  (71.6)\hspace{7mm}
 $s_{j+1}^*=s_{r_j+1}^- + \cdots + s_{r_{j+1}}^-$ and
 $e_{j+1}^*=q_{j+1}^*-s_{j+1}^*$ for $j \ge 0$,\vspace{0.8mm}

 \noindent  (71.7)\hspace{7mm}  $\lambda_j^*=\lambda_{r_j}-1$ if $j \ge 0$ is even, and $\lambda_j^*=\lambda_{r_j}$ if $j\ge 0$ is odd.\vspace{0.5mm}

\noindent    Suppose $j \ge 0$. We will prove  that\vspace{0.5mm}

  \noindent (71.8)  \hspace{9.5mm} $s^*_{j+1}=min\{q_{j+1}^*,\lambda_j^*\}$ and
$\lambda^*_{j+1}=
\lambda_j^*-s_{j+1}^*$ if $j$ is even,\vspace{0.5mm}

\noindent (71.9) \hspace{3mm} $s^*_{j+1}=min\{q_{j+1}^*,p^*+1-\lambda_j^*\}$ and
 $\lambda^*_{j+1}=\lambda_j^*-s_{j+1}^*$
if $j$ is odd,\vspace{0.5mm}

\noindent (71.10) \hspace{37mm} $e_{j+1}^*=q_{j+1}^*-s_{j+1}^*$,

\noindent (71.11) \hspace{37mm}  $q_{J+j+1}^*=e_{j}^*+s_{j+1}^*$.\vspace{0.9mm}

\noindent In fact, (71.10) is trivial. Moreover, (71.8), (71.9) and (71.11) follow from Observation 74.2, 74.3 and 74.4 in Section 74.
 Since $\lambda_0=p+1$, we note that\vspace{0.3mm}

 \noindent (71.12) \hspace{21mm} $\lambda_0^*=\lambda_{r_0}-1=\lambda_{0}-1=p=p^*+1$.

 \noindent By Observation 72.2 we get that  $Q^*=(q_1^*, \cdots, q_I^*,e_0^*)$.
   Hence, we get according to (71.8),  $\cdots$  ,(71.12) \mbox{that
 $Q^{\,\infty}_*$} is   the   shift  symmetric  vector  generated by $Q^*$ with respect to  $p^*$\hspace{0.2mm}.  That means,
 $Q^{\,\infty}_*=C_{p^*}^{\infty}(Q^*)$. Alternatively,\vspace{0.5mm}

\noindent (71.13) \hspace{32.5mm}  $\pi(Q^{\,\infty})=C_{p-1}^{\infty}(\pi(Q))$.\vspace{2mm}

\noindent{{\textbf{Observation 71.1.}}} Let $j \ge 0$. Then $G_{j+1}$ is an  odd  component
 succeeded by a coordinate larger than $1$ in $Q^{\infty}$.\vspace{0.7mm}

\noindent{{\textbf{Proof.}}}  By (33.2) and (33.3)  we get that $G_{j+1}=(q_{r_{j}+1}, \cdots ,q_{r_j+2t_{j+1}+1})$ where $q_{r_j+2i}=1$ for $1 \le i \le t_{j+1}$, and $q_{r_j+2t_{j+1}+2}>1$.\vspace{2mm}

\noindent  {\large{\textbf{72. The contraction vector of $Q$.}}}\vspace{1mm}

\noindent As in Section 59  let $I$ be maximal such that $r_I \le J$.  By (33.1) we get that

\noindent (72.1) \hspace{20mm} $r_0=0 < r_1 < \cdots < r_I \le J < r_{I+1} < \cdots$\,.\vspace{0.5mm}

\noindent  Let $Q^*=\pi(Q)$. By (72.1) we can decompose $Q=(q_1, \cdots,q_J,e_0)$
as\vspace{0.7mm}

\noindent (72.2) $Q =(G_1, \cdots, G_I,F_0)$  where $G_{j+1}=(q_{r_j+1}, \cdots, q_{r_{j+1}})$  for $0 \le j < I$.\vspace{0.4mm}

\noindent By Observation 59.1 we get that $I$ is odd. In particular $I>0$, and the decomposition in (72.2) is well-defined. If $r_I=J$, then $F_0=(e_0)$.
  If $r_I<J$, then Observation 59.1 implies that\vspace{0.5mm}

\noindent (72.3) \hspace{15mm}  $F_0=(q_{r_I+1}, \cdots, q_J,e_0)=(q_{r_I+1}, \cdots,q_{r_I+2z_0},e_0)$

 \hspace{23mm}where  $q_{r_I+2i}=1$ for $1 \le i \le z_0$.\vspace{0.5mm}

 \noindent  Since  $\alpha=\delta(Q)+1$, then (71.4), (71.5) and (72.2) imply that\vspace{0.6mm}

\noindent (72.4) \hspace{7mm}  $\delta(F_0)+1=\delta(Q)-\delta(G_1, \cdots, G_I)+1=\alpha-\tau(r_I)=e_0^*$.\vspace{2mm}

\noindent{{\textbf{Observation 72.1.}}} a)  $G_{j+1}$ is a proper odd  component in $Q$ for  $0 \le j < I$.\vspace{0.3mm}

\noindent b)  $F_0$ is a proper odd  component in $Q$.\vspace{0.7mm}

\noindent{{\textbf{Proof.}}} a) Let $0 \le j < I$. By Observation 71.1 we get that $G_{j+1}$ is an  odd  component
 succeeded by a coordinate $>1$ in $Q^{\infty}$. If $G_{j+1}$ is succeeded by more than one coordinate in $Q$, then $G_{j+1}$ is
 succeeded by a coordinate $>1$ in $Q$. Hence, we get that $G_{j+1}$ is a proper odd component  in $Q$.\vspace{0.4mm}

\noindent b) Since $F_0$ ends $Q$, it is sufficient to prove that $F_0$ is an odd component. \mbox{If $r_I=J$,} then $F_0=(e_0)$ is an odd component. If $r_I <J$, then (72.3) implies that
$F_0$ is an odd component.\vspace{2mm}

\noindent{{\textbf{Observation 72.2.}}}  $Q^*=(\delta(G_1), \cdots, \delta(G_I),\delta(F_0)+1)=(q_1^*, \cdots, q_I^*,e_0^*)$.\vspace{1mm}

\noindent{{\textbf{Proof.}}}   By (72.2) and Observation 72.1   we get that
 $(G_1, \cdots, G_I,F_0)$ is the component decomposition of $Q$.\vspace{0.5mm}

 \noindent By (71.3) and (72.4) we get
 $\delta(G_{j+1})=q_{j+1}^*$ for $1 \le j < I$, and $\delta(F_0)+1=e_0^*$.
 \mbox{Hence,  Observation 13.1} implies that\vspace{0.4mm}

\hspace{6mm} $Q^*=\pi(Q)=(\delta(G_1), \cdots, \delta(G_I),\delta(F_0)+1)=(q_1^*, \cdots, q_I^*,e_0^*)$.\vspace{2.1mm}

 \noindent {\large{\textbf{73. Crucial relations.}}}\vspace{0.9mm}

 \noindent Suppose $j \ge 0$. Then \vspace{0.5mm}

 \noindent (73.1) \hspace{1mm}$q_{j+1}^*=q_{r_j+1}^- + \cdots + q_{r_{j}+2t_{j+1}+1}^-$
 and $s_{j+1}^*=s_{r_j+1}^- + \cdots + s_{r_{j+2t_{j+1}+1}}^-$,\vspace{0.2mm}

 \noindent  (73.2) \hspace{21.3mm} $s_{r_j+i} \le q_{r_j+i}$ for $i \ge 1$, and $s_{j+1}^* \le  q_{j+1}^*$, \vspace{0.7mm}

\noindent where (73.1) follows from (33.2), (71.3) and (71.6), and  (73.2)  \mbox{from  (28.5).}\vspace{2mm}

\noindent {{\textbf{Observation 73.1.}}} a) $s_{j+1}^*=s_{r_j+1}^-+ \cdots + s_{r_j+2m_{j+1}+1}^-$ for $j \ge 0$.

\noindent b) $s_{j+1}^*=q_{r_j+1}^-+ \cdots + q_{r_j+2m_{j+1}+1}^- -e_{r_j+2m_{j+1}+1}$ for $j \ge 0$.\vspace{0.9mm}

\noindent c) $s_{j+1}^*=\tau(r_j+2m_{j+1}+1)-\tau(r_j) -e_{r_j+2m_{j+1}+1}$ for $j \ge 0$.\vspace{1mm}

\noindent {{\textbf{Proof.}}} a)  follows from Observation 64.1 and (73.1).\vspace{0.7mm}

\noindent  Let $j \ge 0$. By  (61.4) we get that $s_{r_j+i}^-=q_{r_j+i}^-$ for $1 \le i \le 2m_{j+1}$.  Moreover, by (28.8) we get that
$s_{r_j+2m_{j+1}+1}^-= q_{r_j+2m_{j+1}+1}^- -e_{r_j+2m_{j+1}+1}$. Hence, b) follows from a). Moreover, c) follows from b) and (32.8).\vspace{20mm}

\noindent {{\textbf{Observation 73.2.}}} $\tau(r_{I+j+1})=\tau(r_j)+s_{j+1}^*+\alpha$ for
 $j \ge 0$.\vspace{0.5mm}

 \noindent {{\textbf{Proof.}}} Let $j \ge 0$. By (65.1), Observation 52.1 and 73.1 c) we get that\vspace{0.5mm}

\noindent\hspace{1mm} $\tau(r_{I+j+1})=\tau(r_j+J+2m_{j+1}+1)=\tau(r_j+2m_{j+1}+1)-e_{r_j+2m_{j+1}+1}+\alpha$

\noindent\hspace{1mm} $=\tau(r_j)+(\tau(r_j+2m_{j+1}+1)-\tau(r_j)-e_{r_j+2m_{j+1}+1})+\alpha=\tau(r_j)+s_{j+1}^*
+\alpha$.\vspace{2mm}

\noindent  {\large{\textbf{74. Deductions.}}\vspace{1mm}}

\noindent {\textbf{Observation 74.1.}} If $e_{r_j+2m_{j+1}+1}=0$ where $j \ge 0$, then
 $s_{j+1}^* = q_{j+1}^*$.\vspace{0.6mm}

\noindent {\textbf{Proof.}} Suppose $e_{r_j+2m_{j+1}+1}=0$ where $j \ge 0$. By
Observation 61.1 we get that $m_{j+1}=t_{j+1}$. Moreover, we get according to (28.3) and (61.4)  \mbox{that
$s_{r_j+2m_{j+1}+1}=q_{r_j+2m_{j+1}+1}$} and $s_{r_j+i}=q_{r_j+i}$ for $1 \le i \le 2m_{j+1}$.\vspace{0.5mm}

\noindent Since $m_{j+1}=t_{j+1}$, we conclude that $s_{r_j+i}=q_{r_j+i}$ for $1 \le i \le 2t_{j+1}+1$. Hence, we get according
\mbox{to (73.1)} that $s_{j+1}^* = q_{j+1}^*$.\vspace{1.7mm}

\noindent{{\textbf{Observation 74.2.}}} Suppose $j \ge 0$ is even. Then

 \hspace{23.7mm} $\lambda^*_{j+1}=
\lambda_j^*-s_{j+1}^*$  and $s^*_{j+1}=min\{q_{j+1}^*,\lambda_j^*\}$.\vspace{0.6mm}

\noindent {\textbf{Proof.}}  Since $j$ is even and $j+1$ is odd, then $\lambda^*_{j+1}=\lambda_{r_{j+1}}$ and $\lambda_{r_j}-1=\lambda_j^*$. Hence,
we get by Observation 60.3 a) and (73.1) that\vspace{0.01mm}

\noindent (74.1) \hspace{6mm} $\lambda^*_{j+1}=\lambda_{r_{j+1}}=\lambda_{r_j+2t_{j+1}+1}
=\lambda_{r_j}- s_{j+1}^*-1=\lambda_j^*-s_{j+1}^*$.\vspace{0.6mm}

\noindent By (60.4) and (74.1) we get that $\lambda_j^*-s_{j+1}^*=\lambda_{r_{j+1}} \ge 0$. Hence, we get \mbox{that $s_{j+1}^* \le \lambda_j^*$.} Then (73.2) implies that $s^*_{j+1}\le min\{q_{j+1}^*,\lambda_j^*\}$. It is therefore sufficient to prove that $s^*_{j+1}=q_{j+1}^*$ or $s^*_{j+1}=\lambda_j^*$.\vspace{0.5mm}

 \noindent If $e_{r_j+2m_{j+1}+1}=0$, then Observation 74.1  implies that
  $s_{j+1}^* = q_{j+1}^*$.\vspace{0.6mm}

  \noindent Next, suppose $e_{r_j+2m_{j+1}+1}>0$. Then we get from Observation 62.1 \mbox{and (74.1)}  \mbox{that
$\lambda_j^*-s_{j+1}^*=\lambda_{r_j+2t_{j+1}+1}=0$.} Hence, $s_{j+1}^*=\lambda_j^*$.\vspace{2mm}

\noindent{{\textbf{Observation 74.3.}}} Suppose $j \ge 0$ is odd. Then\vspace{0.3mm}

 \hspace{16.7mm} $\lambda^*_{j+1}=
\lambda_j^*+s_{j+1}^*$  and $s^*_{j+1}=min\{q_{j+1}^*,p^*+1-\lambda_j^*\}$.\vspace{0.7mm}

\noindent {\textbf{Proof.}} Since $j$ is odd and $j+1$ is even, then $\lambda^*_{j+1}=\lambda_{r_{j+1}}-1$ and $\lambda_{r_j}=\lambda_j^*$. Hence,
we get by Observation 60.4 a) and (73.1) that\vspace{0.01mm}

\noindent (74.2) \hspace{1mm} $\lambda^*_{j+1}=\lambda_{r_{j+1}}-1=\lambda_{r_j+2t_{j+1}+1}-1
=\lambda_{r_j}+ s_{j+1}^*+1-1=\lambda_j^*+s_{j+1}^*$.\vspace{0.8mm}

\noindent By (60.4) and (74.2) we get that
$\lambda_j^*+s_{j+1}^*=\lambda_{r_{j+1}}-1 \le p=p^*+1$.

\noindent Hence, $s_{j+1}^* \le p^*+1-\lambda_j^*$. Then (73.2) implies that\vspace{0.5mm}

\hspace{37mm}$s^*_{j+1}\le min\{q_{j+1}^*,p^*+1-\lambda_j^*\}$.\vspace{0.5mm}

 \noindent Hence, it is  sufficient to prove \mbox{that $s^*_{j+1}=q_{j+1}^*$ or $s^*_{j+1}=p^*+1-\lambda_j^*$.}\vspace{0.5mm}

 \noindent If $e_{r_j+2m_{j+1}+1}=0$, then Observation 74.1 implies that
  $s_{j+1}^* = q_{j+1}^*$.\vspace{0.6mm}

  \noindent Suppose $e_{r_j+2m_{j+1}+1}>0$. Then we get by Observation 63.1 and (74.2)  \mbox{that
$\lambda_j^*+s_{j+1}^*=\lambda_{r_j+2t_{j+1}+1}-1=p=p^*+1$.} Hence, $s_{j+1}^*=p^*+1-\lambda_j^*$.\vspace{200mm}

\noindent {{\textbf{Observation 74.4.}}} $q_{I+j+1}^*= e_{j}^*+s_{j+1}^*$ for $j \ge 0$.\vspace{0.5mm}

\noindent {{\textbf{Proof.}}}    Since $\tau(r_0)=\tau(0)=0$, then  we get according to Observation 73.2  that $\tau(r_{I+1})=s_1^*+\alpha$. Since $e_0^*=\alpha-\tau(r_I)$, then we get from (71.3) that\vspace{0.5mm}

\noindent
\hspace{20mm}$ q_{I+1}^*=\tau(r_{I+1})-\tau(r_I)=(s_1^*+\alpha)- (\alpha-e_0^*)=e_0^*+s_1^*$.\vspace{0.5mm}

\noindent  Let $j \ge 1$.  By  (71.3), Observation 73.2 and (71.10) we get that\vspace{0.5mm}

\noindent       \hspace{2.9mm}     $q_{I+j+1}^*=\tau(r_{I+j+1})-\tau(r_{I+j})=(\tau(r_j)+s_{j+1}^*+\alpha)- (\tau(r_{j-1})+s_{j}^*+\alpha)$.\vspace{0.2mm}

\noindent \hspace{16.7mm}$=(\tau(r_j) -\tau(r_{j-1}))-s_j^*+s_{j+1}^*=q_j^*-s_j^*+s_{j+1}^*=e_j^*+s_{j+1}^*$.\vspace{2.8mm}

\hspace{47mm} {\large{\textbf{PART 12.}}}  \vspace{ 0.5 mm}

\noindent We suppose $A^{\,\infty}=a_1a_2 \cdots  $ is   generated from $A=a_1 \cdots a_n$ by the symmetric shift register $\theta$ with parameters $k$, $p$ and $n$ where $k \le w(A) \le k+p+1$. Let\vspace{0.4mm}

\noindent\hspace{1mm} $A_r=a_{r+1} \cdots, a_{r+n}$ and  $w_r=w(a_{r+1} \cdots a_{r+n})-k=w(A_r)-k$ for $r \ge 0$.\vspace{0.6mm}

\noindent We note that $A_r=\theta^{\,r}(A)$ for $r \ge 0$. We also let  $b_r=w_r+w_{r+n}$ for $r \ge 0$.  The goal of this part is to determine appropriate start strings. In fact, in Section 80 we will determine $r \ge 0$ such that $A_r=\theta^{\,r}(A)$ satisfies (19.1).\vspace{2mm}

\noindent {\large{\textbf{75. Properties of the weight parameters.}}}\vspace{0.8mm}

\noindent In this section we suppose $r \ge 0$. By Observation 42.1 a) we get that

\noindent (75.1) \hspace{23mm} $0 \le w_r \le p+1$ and $0 \le w_{r+n} \le p+1$.\vspace{0.5mm}

\noindent  The next observations follow from (42.5), $\cdots$, (42.8) and (75.1).\vspace{1mm}

\noindent  {{\textbf{Observation 75.1.}}} Suppose $a_{r+1}=1$, $w_r>0$ and $w_{r+n}>0$.

\noindent a) If  $w_{r+n}< p+1$, then $w_{r+1}=w_r-1$, $a_{r+n+1}=0$ and $w_{r+n+1}=w_{r+n}+1$.

\noindent b) If  $w_{r+n}= p+1$, then $w_{r+1}=w_r-1$, $a_{r+n+1}=0$ and $w_{r+n+1}=w_{r+n}$.\vspace{2mm}

\noindent  {{\textbf{Observation 75.2.}}} Suppose $a_{r+1}=0$, $w_r>0$ and $w_{r+n}>0$.

\noindent a) If  $w_r <p+1$, then $w_{r+1}=w_r+1$, $a_{r+n+1}=1$
and $w_{r+n+1}=w_{r+n}-1$.

\noindent b) If  $w_r =p+1$, then $w_{r+1}=w_r$, $a_{r+n+1}=0$ and $w_{r+n+1} \ge w_{r+n}$.\vspace{2mm}

\noindent  {{\textbf{Observation 75.3.}}} Suppose $a_{r+1}=0$, $w_r<p+1$ and $w_{r+n}<p+1$.

\noindent a) If  $w_{r+n}>0$, then $w_{r+1}=w_r+1$, $a_{r+n+1}=1$ and $w_{r+n+1}=w_{r+n}-1$.

\noindent b) If  $w_{r+n}= 0$, then $w_{r+1}=w_r+1$, $a_{r+n+1}=1$ and $w_{r+n+1}=w_{r+n}$.\vspace{2mm}

\noindent  {{\textbf{Observation 75.4.}}} Suppose $a_{r+1}=1$, $w_r<p+1$ and $w_{r+n}<p+1$.

\noindent a) If  $w_r >0$,  then $w_{r+1}=w_r-1$, $a_{r+n+1}=0$
and $w_{r+n+1}=w_{r+n}+1$.

\noindent b) If  $w_r =0$, then $w_{r+1}=w_r$, $a_{r+n+1}=1$ and $w_{r+n+1} \le w_{r+n}$.\vspace{2mm}

\noindent {\large{\textbf{76. Auxiliary results.}}}\vspace{0.9mm}

\noindent {{\textbf{Observation 76.1.}}} Suppose $0 < m < p+1$, $b_r=w_r+w_{r+n}$ and $r \ge 0$.

\noindent a) If $b_r \ge m+p+1$\,, then $m \le w_r \le p+1$ and $m \le w_{r+n} \le p+1$.\vspace{0.4mm}

\noindent b) If $b_r \le m$\,, then $0 \le w_r \le m$ and $0 \le w_{r+n} \le m$.\vspace{0.4mm}

\noindent {{\textbf{Proof.}}} Since $b_r=w_r+w_{r+n}$, the results follow from (75.1).\vspace{2mm}

\noindent {{\textbf{Observation 76.2.}}}
  a) If $ A_{r+z}=A_r$ where $z>0$,  then $z$ is a period of $A^{\,\infty}$.

\noindent  b) If $z$ is period of $A^{\,\infty}$, then $A_{i+z}=A_i$ and $w_{i+z}=w_i$ for $i \ge 0$.\vspace{0.5mm}

\noindent  c) Suppose $0<   w_{r+i} < p+1$ for $0 \le  i  < 2n$. Then $2n$ is a period of $A^{\,\infty}$.\vspace{0.6mm}

\noindent {{\textbf{Proof.}}} a) and b) are trivial.\vspace{0.1mm}

\noindent  c) Since  $0 < w_{r+i} < p+1$ and  $0 < w_{r+n+i} < p+1$ \mbox{for $0 \le  i  < n$,}
then Observation 44.1 a) implies that
 $a_{r+n+i}=a_{r+i}^{\,\prime}$  and $a_{r+2n+i}=a_{r+n+i}^{\,\prime}=a_{r+i}$ \mbox{for $1 \le i \le n$.} Hence, $A_{r+2n}=A_r$ and the result follows from  a).\vspace{2mm}

\noindent {{\textbf{Observation 76.3.}}} Suppose $r \ge 0$ and $m \ge 0$. If $w_i \le m$ for $i >r$, \mbox{then $w_i \le m$ for $i \ge 0$.}
 If $w_i \ge m$ for $i >r$, then $w_i \ge m$ for $i \ge 0$.\vspace{0.4mm}

\noindent {{\textbf{Proof.}}}  Suppose  $w_i \le m$ for $i >r$. Let  $i \ge 0$. We choose a period $z >r$   \mbox{of $A^{\,\infty}$.}
Since $i+z >r$, then  Observation 76.2 b) implies that $w_i=w_{i+z} \le m$. The last statement is proved in the same way.\vspace{2mm}

\noindent {{\textbf{Observation 76.4.}}}  Suppose  $z$ is a period of $A^{\,\infty}$ and $r \ge 0$.

\noindent a) If $w_{r+i} \le m$  for $0 \le i < z$,
then $w_i \le m$ for $i \ge 0$.

\noindent b) If
$ w_{r+i}\ge m$  for $0 \le i < z$,
then $ w_i\ge m$ for $i \ge 0$.\vspace{0.5mm}

\noindent {{\textbf{Proof.}}} a) According to Observation 76.3 it is sufficient to prove that $w_i \le m$
\mbox{for $i \ge r+z$.} Let $i \ge r+z$. Then $i=r+xz+i^*$ where $x \ge 1$ and $0 \le i^* < z$. Since $xz$ is a period of
$A^{\,\infty}$, then Observation 76.2 b) implies that

\hspace{38mm}$w_i=w_{r+xz+i^*}=w_{r+i^*} \le m$.\vspace{0.1mm}

\noindent b) is proved in the same way.\vspace{2mm}

\noindent {\large{\textbf{77. Basic results.}}}\vspace{1mm}

\noindent We suppose in this section that   $0 < m < p+1$. We define\vspace{0.5mm}

\hspace{9.5mm}$B_m^+=\{r \ge 0: b_r \ge m+p+1$ and $w_{r+i} \ge m$ for $1 \le i \le n\}$,

\hspace{9.5mm}$B_m^-=\{r \ge 0: b_r \le m$ and $w_{r+i} \le m$ for $1 \le i \le n\}$.\vspace{1.5mm}

\noindent  {{\textbf{Observation 77.1.}}} Suppose $r \in B_m^+$ and $b_{r+1} \ge m+p+1$.
Then $r+1 \in B_m^+$.

\noindent {\textbf{Proof.}} It is sufficient to prove that $w_{r+1+i} \ge m$ for $1 \le i \le n$.

\noindent Since $b_{r+1} \ge m+p+1$, then Observation 76.1 a) implies that $w_{r+1+n} \ge m$.
Moreover, since $r \in B_m^+$, then $w_{r+1+i}=w_{r+(i+1)} \ge m$ for $1 \le i < n$.\vspace{2mm}

\noindent  {{\textbf{Observation 77.2.}}} Suppose $r \in B_m^+$.\vspace{0.3mm}

\noindent a) If  $a_{r+1} =1$ and $w_{r+n}< p+1$, then $b_{r+1}=w_{r+n+1}+w_{r+1} = w_{r+n}+w_{r}$.

\noindent b) If  $a_{r+1} =1$ and $w_{r+n}= p+1$, then
$b_{r+1}=w_{r+n+1}+w_{r+1} = w_{r+n}+w_{r}-1$.

\noindent c) If  $a_{r+1} =0$, then   $b_{r+1}=w_{r+n+1}+w_{r+1} \ge w_{r+n}+w_{r}$.

\noindent d) $b_{r+1} \ge m+p+1$ and $r+1 \in B_m^+$.\vspace{0.5mm}

\noindent {\textbf{Proof.}} By assumption   $w_{r+n}+w_{r}=b_r \ge m+p+1$.
  Observation 76.1 a) implies  that $0 <m \le w_r \le p+1$ and  $0 <m \le w_{r+n} \le p+1$.\vspace{0.4mm}

 \noindent a), b) and c) follow from Observation 75.1 a), 75.1 b) and 75.2 respectively.\vspace{0.5mm}

 \noindent d) By Observation 77.1 it is sufficient to prove that $b_{r+1} \ge m+p+1$.

\noindent If  $a_{r+1} =1$ and $w_{r+n}< p+1$, this follows from a). If  $a_{r+1} =0$, this follows from c).
  Otherwise, $a_{r+1}=1$ and $w_{r+n}=p+1$. If $w_{r}=m$, \mbox{then  (42.5)} implies  that $w_{r+1}=m-1$. This is a contradiction since $r \in B_m^+$.
Hence,  $w_{r} > m$ \mbox{and $b_r=w_r+w_{r+n}>m+p+1$.} Then b) implies that $b_{r+1} \ge m+p+1$.\vspace{2mm}

\noindent  {{\textbf{Observation 77.3.}}} Suppose there exists  $r \in B_m^+$. \mbox{Then $w_i \ge m$ for $i \ge 0$.}\vspace{0.6mm}

\noindent {\textbf{Proof.}}   Observation 77.2 d) implies that $r+i \in B_m^+$ for $i \ge 0$. Hence, $w_{r+i} \ge m$ for $i \ge 1$. By  Observation 76.3 we  get that $w_i \ge m$ for $i \ge 0$.\vspace{2mm}

\noindent  {{\textbf{Observation 77.4.}}} Suppose $r \in B_m^-$ and $b_{r+1} \le m$.
Then $r+1 \in B_m^-$.

\noindent {\textbf{Proof.}} It is sufficient to prove that $w_{r+1+i} \le m$ for $1 \le i \le n$.

\noindent Since $b_{r+1} \le m$, then Observation 76.1 b) implies that $w_{r+1+n} \le m$.
Moreover, since $r \in B_m^-$, then $w_{r+1+i}=w_{r+(i+1)} \le m$ for $1 \le i < n$.\vspace{2mm}

\noindent  {{\textbf{Observation 77.5.}}} Suppose $r \in B_m^-$.\vspace{0.3mm}

\noindent a) If  $a_{r+1} =0$ and $w_{r+n}>0$, then $b_{r+1}=w_{r+n+1}+w_{r+1} = w_{r+n}+w_{r}$.

\noindent b) If  $a_{r+1} =0$ and $w_{r+n}= 0$, then
$b_{r+1}=w_{r+n+1}+w_{r+1} = w_{r+n}+w_{r}+1$.

\noindent c) If  $a_{r+1} =1$, then   $b_{r+1}=w_{r+n+1}+w_{r+1} \le w_{r+n}+w_{r}$.

\noindent d) $b_{r+1} \le m$ and $r+1 \in B_m^-$.\vspace{0.5mm}

\noindent {\textbf{Proof.}} By assumption   $w_{r+n}+w_{r}=b_r \le m$.
  Observation 76.1 b) implies  that $0 \le w_r \le   m< p+1 $ and  $0 \le  w_{r+n} \le m < p+1$.\vspace{0.4mm}

 \noindent a), b) and c) follow from Observation 75.3 a), 75.3 b) and 75.4 respectively.\vspace{0.7mm}

 \noindent d) By Observation 77.4 it is sufficient to prove that $b_{r+1} \le m$.
 \mbox{If  $a_{r+1} =0$} \mbox{and $w_{r+n}> 0$,} this follows from a). If  $a_{r+1} =1$, this follows from c).
  Otherwise, $a_{r+1}=0$ and $w_{r+n}=0$. If $w_{r}=m$, then we get according \mbox{to  (42.7)}
  \mbox{that $w_{r+1}=m+1$.} This is a contradiction since $r \in B_m^-$.
We conclude that  $w_{r} < m$ \mbox{and $b_r=w_r+w_{r+n}<m$ in this case.} Then b) implies  $b_{r+1} \le m$.\vspace{2mm}

\noindent  {{\textbf{Observation 77.6.}}} Suppose there exists $r \in B_m^-$. Then $w_i \le m$ for $i \ge 0$.\vspace{0.5mm}

\noindent {\textbf{Proof.}}   Observation 77.5 d) implies that $r+i \in B_m^-$ for $i \ge 0$. Hence, $w_{r+i} \le m$ for $i \ge 1$. By  Observation 76.3 we  get that $w_i \le m$ for $i \ge 0$.\vspace{2mm}

\noindent {\large{\textbf{78. The determination of upper and lower bounds.}}}\vspace{1mm}

\noindent Let $x=min\{w_i : 0 \le i \le 2n\}$ and $y=max\{w_i : 0 \le i \le 2n\}$.\vspace{1.5mm}

\noindent  {{\textbf{Proposition 78.1.}}} $min\{w_i :  i \ge 0\}=x$ and $max\{w_i : i \ge 0\}=y$.\vspace{0.5mm}

\noindent  {{\textbf{Proof.}}}  It is sufficient to prove \mbox{that
 $x \le w_i \le y$} for $i \ge 0$.  By (75.1) we note that $0 \le x < y \le p+1$.
 The proof is divided into 4 subcases.\vspace{1mm}

\noindent {\textbf{1)}  If $x =0$ and $y =p+1$, then (75.1) implies that  $x \le w_i \le y$ for $i \ge 0$.\vspace{0.5mm}

\noindent {\textbf{2)}  Suppose $0 <x < y < p+1$. Then $0 < w_i < p+1$ for $0 \le i \le 2n$, and Observation 76.2 c) implies that $2n$ is a period of $A^{\,\infty}$.
Since $x \le w_{i} \le y$ \mbox{for $0 \le i \le 2n$,} then Observation 76.4 implies that $x \le w_i \le y$ for $i \ge 0$.\vspace{0.5mm}

\noindent {\textbf{3)}  Suppose $0 <x < y = p+1$. Choose $j$ such that $w_{j}=p+1$ \mbox{and $0 \le j \le 2n$}. If $j \le n$, let $r =j$. If $j >n$, let $r=j-n$.
Then $w_r=p+1$ or $w_{r+n}=p+1$ where $0 \le r \le n$. Moreover, $0 \le r+i \le 2n$ and
 $w_{r+i} \ge x$ for $0 \le i \le n$. Hence, $b_r=w_r+w_{r+n} \ge p+1+x$ and $r \in B_x^+$.\vspace{0.3mm}

\noindent Then
 $x \le w_i \le p+1=y$ for $i \ge 0$, where the first inequality follows from Observation 77.3 and the  last inequality follows from (75.1).\vspace{0.65mm}

\noindent {\textbf{4)}} Suppose $0 =x < y <p+1$. Choose $j$ such that $w_{j}=0$ and $0 \le j \le 2n$. \mbox{If $j \le n$,} we let $r =j$. If $j >n$, we let $r=j-n$.
Then $w_r=0$ or $w_{r+n}=0$ where $0 \le r \le n$. Moreover, $0 \le r+i \le 2n$ and $ w_{r+i} \le y$ for $0 \le i \le n$.
Hence, $b_r=w_r+w_{r+n} \le y$ and  $r \in B_y^-$.\vspace{0.3mm}

\noindent Then
$x = 0 \le w_i \le y$ for $i \ge 0$, where the  first inequality follows from (75.1)
and the last inequality follows from Observation 77.6.\vspace{1.5mm}

\noindent {\large{\textbf{79. Adjustment of the parameters.}}}\vspace{1mm}

\noindent We will determine parameters $k^* \ge 0$ and $p^* \ge 0$ such that\vspace{0.5mm}

\noindent (79.1) \hspace{2mm}$min\{w(A_i): i \ge 0\}=k^*$ and $max\{w(A_i): i \ge 0\}=k^*+p^*+1$\hspace{0.15mm},\vspace{0.5mm}

\noindent and $A^{\,\infty}=a_1a_2 \cdots  $ is   generated from $A=a_1 \cdots a_n$  by the symmetric shift register  with parameters $k^*$, $p^*$ and $n$.\vspace{0.4mm}

\noindent Let $x=min\{w_i : 0 \le i \le 2n\}$ and $y=max\{w_i : 0 \le i \le 2n\}$. According to Proposition 78.1 we get that\vspace{0.5mm}

\noindent (79.2) \hspace{17mm}$min\{w_i : i \ge 0\}=x$ and $max\{w_i : i \ge 0\}=y$.\vspace{0.5mm}

\noindent Since $w(A_i)=k+w_i$ for $i \ge 0$, then  we also  get that\vspace{0.6mm}

\noindent (79.3)\hspace{6mm}$min\{w(A_i): i \ge 0\}=k+x$
 and $max\{w(A_i): i\ge 0\}=k+y$.\vspace{0.6mm}

\noindent Let  $k^*=k+x$ and $p^*=y-x-1$. Since $k+x=k^*$ and $k+y=k^*+p^*+1$, then (79.1) follows from (79.3).\vspace{0.6mm}

\noindent We will prove that $A^{\,\infty}=a_1a_2 \cdots  $ is   generated from $A=a_1 \cdots a_n$ by the symmetric shift register  with parameters $k^*$, $p^*$ and $n$.
 If $x=0$ and $y=p+1$, then $k^*=k$ and $p^*=p$ and this is trivial. Otherwise this follows from Proposition 79.1. By (75.1) we get that $0 \le x < y \le p+1$. Hence,\vspace{1mm}

\noindent (79.4)  \hspace{1mm}$k^*=k+x \ge k$ and  $k^*+p^*=k+x+(y-x-1)=k+y-1 \le k+p$.\vspace{1.5mm}

 \noindent {\textbf{Proposition 79.1.}}  Let $r \ge 0$\hspace{0.25mm}. Then\vspace{0.2mm}

 \hspace{6mm} $a_{r+n+1}=a_{r+1}^{\,\prime}$ if and only if $k^* \le a_{r+2}+ \cdots + a_{r+n} \le k^*+p^*$.\vspace{0.75mm}

 \noindent {\textbf{Proof.}} First, we suppose $k^* \le a_{r+2}+ \cdots + a_{r+n} \le k^*+p^*$.
 By (79.4) we get that $k \le a_{r+2}+ \cdots + a_{r+n} \le k+p$\hspace{0.3mm}. Hence, $a_{r+n+1}=a_{r+1}^{\,\prime}$\hspace{0.2mm}.\vspace{0.5mm}

 \noindent Suppose $a_{r+2}+ \cdots + a_{r+n}<k^*$. Then we get that $w(A_r)<k^*$ if $a_{r+1}=0$. Moreover,
 $w(A_{r+n+1}) < k^*$ if  $a_{r+n+1}=0$. By (79.1) we get that $w(A_r)\ge k^*$ and $w(A_{r+n+1})\ge k^*$
  Hence,  $a_{r+n+1}=a_{r+1}=1$.\vspace{0.5mm}

 \noindent Suppose
  $a_{r+2}+ \cdots + a_{r+n}>k^*+p^*$. Then  $w(A_r)>k^*+p^*+1$ if $a_{r+1}=1$. Moreover,
   $w(A_{r+n+1}) > k^*+p^*+1$ if  $a_{r+n+1}=1$.
  By (79.1) we get that $w(A_r)\le k^*+p^*+1$ and $w(A_{r+n+1})\le k^*+p^*+1$.
  Hence,  $a_{r+n+1}=a_{r+1}=0$.\vspace{0.5mm}\vspace{2mm}

  \noindent {\large{\textbf{80. A crucial algorithm.}}}\vspace{0.5mm}

\noindent Suppose $min\{w(A_i): i \ge 0\}=k$ and $max\{w(A_i): i\ge 0\}=k+p+1$.
By Section 79 we can always transfer the determination of minimal periods to this case. We get that  $min\{w_i: i \ge 0\}=0$ and $max\{w_i: i\ge 0\}=p+1$.\vspace{1.5mm}

 \noindent {\textbf{Algorithm 80.1.}} Let $s$ be the least integer  $\ge 0$ \mbox{such that $w_s=p+1$.}
 \mbox{Let $t>s$} be minimal such that $w_{t}=0$. Let  $r$ be maximal such that $s \le r< t$ and $w_{r}=p+1$.\vspace{0.5mm}

 \noindent Since $w_y>0$ for $s < y < t$, and  $w_y<p+1$ for $r < y < t$, then\vspace{0.5mm}

 \noindent (80.1)  \hspace{30.8mm} $0 < w_y < p+1$ if $r < y < t$.\vspace{0.5mm}

 \noindent Let  $x=t-r$.  We will prove that\vspace{0.5mm}

 \noindent (80.2)  \hspace{24.5mm}  $s \le 2n$,  $t \le s+n \le 3n$ and  $r <3n$,\vspace{0.3mm}

 \noindent (80.3)\hspace{7mm}  $w_{r}=p+1$, $w_{r+x}=0$ and $0 < w_{r+i} <p+1$ for $1 \le i < x$,\vspace{0.5mm}

  \noindent (80.4) \hspace{10mm}$w(A_r)=k+p+1$, $A_r$ starts with 1 and $V(A_r) \in M_p^+$.\vspace{0.6mm}

  \noindent By (80.4) we get that $A_r$ satisfies (19.1).  Then we can determine the minimal period of $A_r$ as  in Part 2. This will also be the minimal period of $A$.\vspace{2mm}

\noindent {\textbf{Proof of (80.2).}} Proposition 78.1 implies
 \mbox{that $max\{w_i: 0 \le i \le 2n\}=p+1$.} Hence, $s \le 2n$. Suppose $t>s+n$.   Then we get\vspace{0.5mm}

\hspace{2mm}$m=min\{w_{s+i}: 1 \le i \le n\}>0$ and $b_s=w_s+w_{s+n} \ge m+p+1$.\vspace{0.5mm}

\noindent Hence, $s \in B_m^+$ and Observation 77.3  implies that  $w_i\ge m>0$ for $i \ge 0$.  This is a contradiction. Hence, $t \le s+n \le 3n$ and $r < t \le 3n$.\vspace{2mm}

\noindent {\textbf{Proof of (80.3).}}
 We note that    $w_{r}=p+1$  and $w_{r+x}=w_{t}=0$. \mbox{Let $1 \le i < x$.} \mbox{Since $r < r+i < r+x=t$,} then (80.1) implies that $0 < w_{r+i} < p+1$.\vspace{2mm}

\noindent {\textbf{Proof of (80.4).}}
By (80.2) we get that $0<t-r \le t-s \le n$. Hence, $0 < x \le n$. Then we get by  (80.3) and Proposition 44.2 d)  \mbox{that (80.4) is true.}\vspace{200mm}

\hspace{45mm} {\large{\textbf{PART 13.}}  \vspace{ 1mm}}

 \noindent If $V \in M^*$, we will prove in Section 81 that $\pi(V) \in M$. Moreover, if $V \in M_p^+$ where $p>0$, we will prove in Section 84 that
 $\pi(V) \in M^+_{p-1}$.\vspace{1mm}

 \noindent {\large{\textbf{81. Assumptions and notation.}}}\vspace{1mm}

 \noindent Let $V=(v_1, \cdots , v_{J+1}) \in M^*$, and let $\tau$ be the distance function of $V$. Then\vspace{0.4mm}

  \noindent (81.1)  \hspace{16mm} $J$ is odd, $v_1>1$, $v_i \ge 1$ for $2 \le i \le J$, and $v_{J+1} \ge 0$,\vspace{0.3mm}

  \noindent (81.2)  \hspace{26mm} $v_i^- \ge 0$ for $1 \le i \le J$, and $v_{J+1}^-+1 \ge 0$,\vspace{0.5mm}

  \noindent (81.3)  \hspace{35mm} $0<\tau(1) \le \tau(2) \le   \cdots \le \tau(J)$.\vspace{0.5mm}

\noindent If $0 \le r \le J$, let  $t_{max}(r)=t$ and $next(r)=r+2t+1=r+2t_{max}(r)+1$ where $t \ge 0$ is maximal
such that  $r+2t \le J$ and $v_{r+2i}=1$ for $1 \le i \le t$.

\noindent Let $r_0, \cdots, r_{I+1}$ be the $r$\,-\,indexes of $V$. \mbox{Then $r_{j+1}=next(r_j)$ for $0 \le j \le I$, and}\vspace{0.4mm}

\noindent (81.4)\hspace{22mm}$r_0=0 < r_1 < \cdots < r_{I-1} <r_I < r_{I+1}= J+1$.\vspace{0.4mm}

\noindent Let $t_{j+1}=t_{max}(r_j)$ for $0 \le j \le I$. Suppose $0 \le j \le I$. Then\vspace{0.5mm}

\noindent (81.5) \hspace{13mm} $r_{j+1}=r_j+2t_{j+1}+1$ and $v_{r_j+2i}=1$ for $1 \le i \le t_{j+1}$, \vspace{0.3mm}

 \noindent (81.6) \hspace{14mm} $v_{r_j+2t_{j+1}+2}=v_{r_j+2(t_{j+1}+1)}>1$  if  $r_j+2t_{j+1} \le J-2$.\vspace{0.7mm}

 \noindent We note that  $r_{j+1}-r_j$ is odd for $0 \le j \le I$, and $r_0$ is even. Then\vspace{0.5mm}

\noindent (81.7) \hspace{22mm}  $r_j$ is even $\Leftrightarrow$ $j$ is even for  $0 \le j \le I+1$.\vspace{0.4mm}

\noindent In particular, since $r_{I+1}=J+1$ is even, then $r_I$ is odd. Hence, (81.7) implies that $I$ is odd.
   By (13.1)
   the contraction vector  of $V$ is\vspace{0.4mm}

\noindent (81.8) \hspace{3.5mm}$\pi(V)=(v_1^*, \cdots, v_{I+1}^*)$ where  $v_{j+1}^*=\tau(r_{j+1})-\tau(r_{j})$  for $0 \le j < I$,\vspace{0.2mm}

\hspace{9.5mm}and $v_{I+1}^*=\tau(r_{I+1})-\tau(r_{I})+1=\tau(J+1)-\tau(r_{I})+1$.\vspace{0.6mm}

\noindent By (10.2) and (81.8)  we get that\vspace{0.5mm}

\noindent (81.9) \hspace{9.6mm}$v_{j+1}^*= \tau(r_{j+1})-\tau(r_{j})=v_{r_j+1}^-+ \cdots + v_{r_{j+1}}^-$ for $0 \le j < I$.\vspace{0.5mm}

\noindent Moreover, (10.2), (81.2) and  (81.8)  imply that\vspace{0.4mm}

\noindent (81.10) \hspace{6.5mm} $v_{I+1}^*= \tau(J+1)-\tau(r_{I})+1=v_{r_{I}+1}^- + \cdots+ v_{J+1}^-+1 \ge 0$.\vspace{0.5mm}

\noindent Let $\rho_0,  \cdots, \rho_{J+1}$ be the alternating parameters of $V$ as in Section 14.

\noindent We note that $\rho_0=0$. If $0 \le i \le J$, then by (14.1) we get that\vspace{0.2mm}

\noindent (81.11) \hspace{4mm} $\rho_{i+1}=\rho_{i}+v_{i+1}$ if $i$ is even, and
$\rho_{i+1}=\rho_{i}-v_{i+1}$ if $i$ is odd.\vspace{1.5mm}

\noindent {\textbf{Observation 81.1.}} a) $v_{r_j+1}>1$ for $0 \le j <I$.\vspace{0.5mm}

\noindent b) If $0 \le j <I$, then $r_{j+1}\le r_I \le J$ and $v_{j+1}^*=v_{r_j+1}^-+ \cdots + v_{r_{j+1}}^->0$.\vspace{0.5mm}

 \noindent {\textbf{Proof.}} a) By (81.1) we get that $v_{r_0+1}=v_1>1$. Suppose  $1 \le j < I$.
 Then we get from (81.4) and (81.5)  that $r_j<r_I \le J$ and $r_{j-1}+2t_j=r_j-1 \le J-2$.
 Hence, (81.5) and (81.6) imply that $v_{r_j+1}=v_{r_{j-1}+2t_j+2}>1$.\vspace{0.5mm}

 \noindent b) Let  $0 \le j <I$. By (81.4) we get that $r_{j+1}\le r_I \le J$. Hence,  a), (81.2) and (81.9) imply that $v_{r_j+1}^->0$ and $v_{j+1}^*=v_{r_j+1}^-+ \cdots + v_{r_{j+1}}^->0$.\vspace{2mm}

\noindent {\textbf{Proposition 81.2.}}  $v_{j+1}^*>0$ for $0 \le j < I$,  $v_{I+1}^* \ge 0$
 and $\pi(V) \in M$.\vspace{0.7mm}

 \noindent {\textbf{Proof.}}   Observation 81.1 b) and (81.10) imply that
 $v_{j+1}^*>0$ \mbox{for $0 \le j < I$,} and $v_{I+1}^* \ge 0$. Since $I$ is odd, then  $\pi(V)=(v_1^*, \cdots, v_{I+1}^*) \in M$.\vspace{2mm}

   \noindent {\large{\textbf{82. Auxiliary results.}}}\vspace{0.5mm}

\noindent {\textbf{Observation 82.1.}} Suppose $0 \le j \le I$, $0 \le i \le 2t_{j+1}$ and $i$ is odd.
  Then\vspace{0.5mm}

  \noindent $v_{r_j+i+1}=1$, \,$v^-_{r_j+i+1}=0$ \,and \,$\tau(r_j+i+1)=\tau(r_j+i)+v^-_{r_j+i+1}=\tau(r_j+i)$.\vspace{0.7mm}

   \noindent {\textbf{Proof.}} Since $i$ is odd and $i+1$ is even, then $1 \le i < 2t_{j+1}$
  and $i+1=2i^*$
  where  $1 \le i^* \le t_{j+1}$. Hence, (81.5) implies that $v_{r_j+i+1}=v_{r_j+2i^*}=1$.\vspace{1.4mm}

\noindent {\textbf{Observation 82.2.}} Suppose $0 \le i \le 2t_{j+1}+1$ where $0 \le j \le I$ and $j$ is even.  Then $\rho_{r_j+i}=\rho_{r_j}+\tau(r_j+i)-\tau(r_j)$ if $i$ is even, and

  \hspace{24mm} $\rho_{r_j+i}=\rho_{r_j}+\tau(r_j+i)-\tau(r_j)+1$ if $i$ is odd.\vspace{0.5mm}

  \noindent {\textbf{Proof.}}  By (81.7) we get that $r_j$ is even. The result is true \mbox{for $i=0$.} Suppose the result is true for $i$ where $0 \le i \le 2t_{j+1}$.  We note that $r_j+i$ is even if and only if $i$ is even. Then we get from Observation 82.1 and (81.11)  that

  \hspace{6mm}$\rho_{r_j+i+1}=\rho_{r_j+i}+v_{r_j+i+1}=\rho_{r_j}+\tau(r_j+i)-\tau(r_j)+v_{r_j+i+1}^-+1$

 \hspace{19.2mm}$=\rho_{r_j}+\tau(r_j+i+1)-\tau(r_j)+1$ if $i$ is even,\vspace{0.5mm}

 \hspace{6mm}$\rho_{r_j+i+1}=\rho_{r_j+i}-v_{r_j+i+1}=\rho_{r_j}+\tau(r_j+i)-\tau(r_j)+1-1$

 \hspace{6mm}$=\rho_{r_j}+\tau(r_j+i)-\tau(r_j)=\rho_{r_j}+\tau(r_j+i+1)-\tau(r_j)$ if $i$ is odd.\vspace{0.5mm}

  \noindent Hence, the result is true for $i+1$.\vspace{2mm}

  \noindent {\textbf{Observation 82.3.}} Suppose $0 \le i \le 2t_{j+1}+1$ where $0 \le j \le I$ and $j$ is odd. Then $\rho_{r_j+i}=\rho_{r_j}-\tau(r_j+i)+\tau(r_j)$ if $i$ is even, and

  \hspace{24mm} $\rho_{r_j+i}=\rho_{r_j}-\tau(r_j+i)+\tau(r_j)-1$ if $i$ is odd.\vspace{0.5mm}

  \noindent {\textbf{Proof.}}  By (81.7) we get that $r_j$ is odd. The result is true \mbox{for $i=0$.} Suppose the result is true for $i$ where $0 \le i \le 2t_{j+1}$.  We note that $r_j+i$ is odd if and only if $i$ is even. Then we get from Observation 82.1 and (81.11)  that

  \hspace{6mm}$\rho_{r_j+i+1}=\rho_{r_j+i}-v_{r_j+i+1}=\rho_{r_j}-\tau(r_j+i)+\tau(r_j)-v_{r_j+i+1}^--1$

 \hspace{19.2mm}$=\rho_{r_j}-\tau(r_j+i+1)+\tau(r_j)-1$ if $i$ is even,\vspace{0.5mm}

 \hspace{6mm}$\rho_{r_j+i+1}=\rho_{r_j+i}+v_{r_j+i+1}=\rho_{r_j}-\tau(r_j+i)+\tau(r_j)-1+1$

 \hspace{6mm}$=\rho_{r_j}-\tau(r_j+i)+\tau(r_j)=\rho_{r_j}-\tau(r_j+i+1)+\tau(r_j)$ if $i$ is odd.\vspace{0.5mm}

  \noindent Hence, the result is true for $i+1$.\vspace{2mm}

  \noindent {\textbf{Observation 82.4.}} Suppose $0 \le j < I$. Then

   \hspace{19mm} $\rho_{r_{j+1}}=\rho_{r_j+2t_{j+1}+1}=\rho_{r_j}+v_{j+1}^*+1$ if $j$ is even,

  \hspace{19mm} $\rho_{r_{j+1}}=\rho_{r_j+2t_{j+1}+1}=\rho_{r_j}-v_{j+1}^*-1$ if $j$ is odd.\vspace{0.8mm}

  \noindent {\textbf{Proof.}} By  (81.5) and (81.8) we get that\vspace{0.6mm}

 \noindent (82.1) \hspace{12mm}$\tau(r_j+2t_{j+1}+1)-\tau(r_j)=\tau(r_{j+1})-\tau(r_j)=v_{j+1}^*$.\vspace{12mm}

    \noindent If $j$ is even, then (81.5),  Observation 82.2 and (82.1) imply that\vspace{0.4mm}

 \noindent\hspace{1mm} $\rho_{r_{j+1}}=\rho_{r_j+2t_{j+1}+1}=\rho_{r_j}+\tau(r_j+2t_{j+1}+1)-\tau(r_j)+1
 =\rho_{r_j}+v_{j+1}^*+1$.\vspace{0.5mm}

  \noindent  If $j$ is odd, then  (81.5), Observation 82.3 and (82.1) imply that\vspace{0.4mm}

 \noindent \hspace{1mm} $\rho_{r_{j+1}}=\rho_{r_j+2t_{j+1}+1}=\rho_{r_j}-\tau(r_j+2t_{j+1}+1)+\tau(r_j)-1=\rho_{r_j}-v_{j+1}^*-1$.
  \vspace{2mm}

\noindent {\large{\textbf{83. The alternating parameters of $\pi(V)$.}}}\vspace{0.5mm}

   \noindent Let $\rho_0^*,  \cdots, \rho_{I+1}^*$ be the alternating parameters of $\pi(V)=(v_1^*, \cdots, v_{I+1}^*)$.

\noindent Then $\rho_0^*=0$ and $\rho_1^*=\rho_0^*+v_1^*=v_1^*$. If $0 \le i \le I$,  then by (14.1) we get that\vspace{0.7mm}

\noindent (83.1) \hspace{5mm}$\rho^*_{i+1}=\rho^*_{i}+v_{j+1}^*$ if $i$ is even,
and $\rho^*_{i+1}=\rho^*_{i}-v_{i+1}^*$ if $i$ is odd.\vspace{2mm}

\noindent {\textbf{Observation 83.1.}} Suppose $0 \le j \le I$. Then $\rho_{r_j}=\rho^*_j$
if $j$ is even, \mbox{and $\rho_{r_j}=\rho^*_j+1$ if $j$ is odd.}\vspace{0.5mm}

 \noindent {\textbf{Proof.}} Since $r_0=0$, then $\rho^*_0=0=\rho_0=\rho_{r_0}$.  Suppose
 this is true for $j$
 \mbox{where $0 \le j < I$.}
 Then (83.1) and Observation 82.4  imply that\vspace{0.3mm}

\hspace{5.6mm} $\rho_{r_{j+1}}=\rho_{r_j}+v_{j+1}^*+1=\rho^*_j+v_{j+1}^*+1=\rho^*_{j+1}+1$ \,if $j$ is even,

\hspace{5.6mm} $\rho_{r_{j+1}}=\rho_{r_j}-v_{j+1}^*-1=\rho^*_j+1-v_{j+1}^*-1=\rho^*_{j+1}$ \,if $j$ is odd.\vspace{0.5mm}

\noindent Hence, the result is true for $j+1$.\vspace{1.5mm}

\noindent {\textbf{Observation 83.2.}}  Suppose $0 \le j \le I$ and  $\rho_i>0$ for $1 \le i \le r_j$. Then\vspace{0.5mm}

\hspace{38mm}$\rho^*_m >0$ for $1 \le m \le j$.\vspace{0.5mm}

 \noindent {\textbf{Proof.}} Suppose  $1 \le m \le j$ where $m$ is even. By (81.4) we get that $1 \le r_m \le r_j$. Then Observation 83.1  implies that
 $\rho^*_m=\rho_{r_m} \ge 1$. By Proposition 81.2 we get \mbox{that $\rho_1^*=v_1^*>0$.} Suppose $1 < m \le j$
 where $m$ is odd. Since $m-1$ is even, then $\rho_{m-1}^* \ge 1$. Hence we get by using (83.1) and Proposition 81.2 that $v_{m}^*>0$ and
 $\rho_{m}^*=\rho_{m-1}^*+v_{m}^*>0$.\vspace{1.5mm}

  \noindent {\textbf{Observation 83.3.}}  Suppose $0 \le j \le I$ where $j$ is odd, and  $0 \le i \le 2t_{j+1}$ where $i$ is even. Then $\rho_{r_j+i} =\rho_{r_j}-\tau(r_j+i)+\tau(r_j) \le \rho_{r_j}=\rho^*_j+1$.\vspace{0.5mm}

  \noindent {\textbf{Proof.}} The first equality follows from Observation 82.3. By (81.4) \mbox{and (81.5)} we get that
$1\le r_j+i \le r_j+2t_{j+1}=r_{j+1}-1 \le r_{I+1}-1 \le J$. Then (81.3) implies  $\tau(r_j+i)\ge\tau(r_j)$ and  the inequality is also true. The last equality follows from Observation 83.1 since $j$ is odd.\vspace{2mm}

 \noindent {\textbf{Observation 83.4.}}  Suppose $0 \le j < I$ where $j$ is even, and  $0 \le i \le 2t_{j+1}+1$ where $i$ is odd. Then $\rho_{r_j+i} \le \rho^*_{j+1}+1$.\vspace{0.7mm}

  \noindent {\textbf{Proof.}} By  (81.4) and (81.5)  we get  $1\le r_j+i \le r_j+2t_{j+1}+1=r_{j+1} \le r_{I} \le J$.

  \noindent Hence,   (81.3) and Observation 82.2 imply that\vspace{0.5mm}

\noindent (83.2)  \hspace{1mm} $\tau(r_j+i)-\tau(r_j)+1 \le \tau(r_j+2t_{j+1}+1)-\tau(r_j)+1=\rho_{r_j+2t_{j+1}+1}$.\vspace{12mm}

\noindent Then  $\rho_{r_j+i} =\tau(r_j+i)-\tau(r_j)+1 \le \rho_{r_j+2t_{j+1}+1}=\rho_{r_{j+1}}=q_{j+1}^*+1$\vspace{0.3mm}
  where the first equality follows from Observation 82.2, the inequality follows \mbox{from (83.2)} and the second equality follows from (81.5). Moreover, the last equality follows from Observation 83.1 since $j+1$ is odd.\vspace{2mm}

\noindent {\large{\textbf{84. The main result.}}}\vspace{1mm}

\noindent We will prove the following result.\vspace{1mm}

\noindent (84.1) \hspace{2mm}If $V \in M^+_{p}$ where $p>0$, then $\pi(V)=(v_1^*, \cdots, v_{I+1}^*) \in M^+_{p-1}=M^+_{p^*}$\vspace{0.5mm}

\noindent where $p^*=p-1$. Suppose $V=(v_1, \cdots, v_{J+1}) \in M^+_{p}$ where $p>0$.\vspace{0.5mm}

\noindent By Observation 14.4 we get that $V \in M^*$. Hence, by Proposition 81.2 we also get that $\pi(V) \in M$. According to Section 14 it is therefore sufficient to prove that\vspace{0.5mm}

\noindent (84.2) \hspace{10mm} $\pi(V)$ has an admissible start vector with respect to $p^*$.\vspace{0.5mm}

\noindent Since $V \in M^+_{p}$, there exists an odd integer $t$ such that $1 \le t \le J$ and\vspace{0.2mm}

\noindent (84.3) \hspace{24.2mm} $\rho_{\,t} \ge p+1$ and $\rho_{\,i}>0$ for $1 \le i \le t$.\vspace{0.5mm}

\noindent  Since $1 \le t \le J$, then  by  (81.4) and (81.5)  there exists $j$ such that\vspace{0.5mm}

\noindent (84.4) \hspace{5mm} $0 \le j \le I$,  $r_j \le t < r_{j+1}$ and $t=r_j+i_*$ where $0 \le i_* \le 2t_{j+1}$.\vspace{0.5mm}

\noindent  Since $r_j \le t$, then (84.3) implies that $\rho_{i}>0$ for $1 \le i \le r_j$.
Let $\rho_0^*,  \cdots, \rho_{I+1}^*$ be the alternating parameters of $\pi(V)=(v_1^*, \cdots, v_{I+1}^*)$. Then we get by Observation 83.2 that

\noindent (84.5) \hspace{35mm} $\rho_m^*>0$ for $1 \le m \le j$.\vspace{0.5mm}

\noindent We divide the proof of (84.2) into two subcases.\vspace{1.5mm}

\noindent {\textbf{Case 1.}} Suppose $j$ is odd. By (81.7) we get that $r_j$ is odd. Since $r_j+i_*=t$ is odd, then $i_*$ is even. Then\vspace{0.5mm}

\hspace{29.8mm}$\rho^*_j  \ge \rho_{r_j+i_*}-1=\rho_{\,t}-1 \ge p=p^*+1$\vspace{0.5mm}

\noindent where the first and second inequality follow from Observation 83.3 and (84.3). The  equalities are trivial.   Hence,  we get from (84.5) that $(v_1^*, \cdots, v_j^*)$ is  an admissible start vector of $\pi(V)=(v_1^*, \cdots, v_{I+1}^*)$ with respect to $p^*$.\vspace{2mm}

\noindent {\textbf{Case 2.}} Suppose $j$ is even. Since $0 \le j \le I$ where $I$ is odd, then $0 \le j < I$. By (81.7) we get that $r_j$ is even. Since $r_j+i_*=t$ is odd, then $i_*$ is odd. Then
we get that

\hspace{27mm}$\rho^*_{j+1}  \ge \rho_{r_j+i_*}-1=\rho_{\,t}-1 \ge p=p^*+1$\vspace{0.5mm}

\noindent where the first and second inequality follow from Observation 83.4 and (84.3). The  equalities are trivial.  Hence,  we get from (84.5) that $(v_1^*, \cdots, v_{j+1}^*)$ is  an admissible start vector of $\pi(V)=(v_1^*, \cdots, v_{I+1}^*)$ with respect to $p^*$.\vspace{200mm}

                          \noindent \textbf{{\large References.}}
\begin{tabbing}
\noindent [1]\hspace{2mm}\= P.G. Drazin \& R.S. Johnson, Solitons: an
introduction, Cambridge, 1989\,.\\

\noindent [2] \> T. Helleseth, Nonlinear Shift Registers, Survey and open problems,\\

\>  Algebraic Curves and Finite Fields, De Gruyter, 2014.\\

\noindent [3] \> K. Kjeldsen, On the cycle structure of
a set of nonlinear shift registers\\

\>  with symmetric feedback functions, J.
Combinatorial Theory,\\

\> Ser. A., 20 (1976), page 154-169.\\

\noindent [4] \> J. Mykkeltveit, Nonlinear recurrences and
arithmetic Codes, Information and \\

\> control, Vol. 33 (1977), page 193-209.\\

\noindent [5] \> J. Mykkeltveit, M. K. Siu, and P. Tong,
On the cycle structure of some nonlinear\\

\> shift register sequences, Information and control,
Vol. 43 (1979), page 202-215\,.\\

\noindent [6] \> J. S{\o }reng, The periods generated by
some symmetric shift registers,\\

\> J. Combinatorial Theory, Ser. A., 21 (1976),
page 165-187.\\

[7] \> J. S{\o }reng, Symmetric shift registers, Pacific
J. Math., 85 (1979), page 201-229.\\

[8] \> J. S{\o }reng, Symmetric shift registers, Pacific
J. Math., 98 (1982), page 203-234.\\
\end{tabbing}

\newpage

         \noindent \textbf{{\large Index}}
\begin{tabbing}
\noindent \hspace{2mm}\= \hspace{53.5mm}Section \hspace{57mm}Section\\

\> $A^{\prime}$, $w(A)$ and $\overline{w}(A)$ \hspace{22mm}\= \hspace{1.3mm}3\hspace{3.5mm}
\= $A^{\,\infty}_p(Q)$ \hspace{49.1mm}\= 18\\

\> $\#V$, $sum(V)$ and $sum(V,j)$ \> \hspace{1.3mm}4 \>cyclic parameters\> 20\\

\> $V+\alpha$ and  $V-\alpha$ \> \hspace{1.3mm}4  \>$\psi$\> 20\\

\> extension of vectors \> \hspace{1.3mm}4\>dynamical parameters\> 21\\

\>$\theta$    \> \hspace{1.5mm}5\> weight parameters, $w_i, w_i^*$ \> 27 \\

\>$V(A^{\,\infty})$ \> \hspace{1.5mm}6\> shift symmetric vector, $C^{\,\infty}_p(Q)$ \> 28 \\

\>  $M$,  $M^*$ and $M_p$\>\hspace{1.3mm}7\>$\lambda_j$, $s_{j+1}$ and $e_j$\> 28\\

\> $V(A)$ and $A(V)$  \> \hspace{1.3mm}7\>complete vector\>32\\

\> $v^-$ and $\delta(V)$ \> \hspace{1.3mm}9\> $t_{max}$ and $next$\hspace{0.2mm}, infinite case \> 32\\

\> $\tau(r)$, finite distance function\> 10\> $\tau(r)$, infinite distance function\> 32\\

\> $t_{max}$ and $next$, finite case\> 11\> $r$- and $t$-indexes, infinite case\> 33\\

\> proper odd component\> 11\> $\pi(Q^{\,\infty})$, infinite contraction vector\>33\\

\> finite component decomposition\>12\> infinite component decomposition\>34\\

\>  $r$-indexes, finite case\>13\>$c$\,-\,indexes, infinite case\>35\\

\>  $\alpha$\>13\>$D(Q^{\,\infty})$, infinite distance vector\>35\\

\> $\pi(V)$,  finite contraction vector  \> 13 \> $r(D^{\,\infty},\beta)$ \>35 \\

\> $\rho_i$,  alternating parameters \> 14\> progression parameters \>35\\

\> $M^+_p$\> 14\>$y(r)$\> 35\\

\> $D(V)$, finite distance vector  \>15\>$\beta_j$ and $y_j$\> 37\\

\> $c$\,-\,indexes, finite case\> 15 \>$\omega(\alpha^*, \gamma^*, j^*, \zeta^*)$ \> 40\\

\>  $\alpha^*$ and $\gamma^*$ \> 16 \>$A_r$ \> 42 \\

\> $gcd\hspace{0.15mm}(\alpha,\gamma)$ \>16  \>$B^+_m$ and $B^-_m$\>77 \\

\end{tabbing}

\end{document}